\documentstyle[amscd,amssymb,verbatim]{amsart} 

\numberwithin{equation}{section}
\newtheorem{theorem}{\hspace*{-1.5cm}}[section]

\newenvironment{thm}{\begin{theorem}{\em{\bf%
Theorem.}\ }\ignorespaces}{\end{theorem}}
\newenvironment{prop}{\begin{theorem}{\em{\bf%
Proposition.}\ }\ignorespaces}{ \end{theorem} }
\newenvironment{lemma}{\begin{theorem}{\em{\bf%
Lemma.}\ }\ignorespaces}{ \end{theorem} }
\newenvironment{cor}{\begin{theorem}{\em{\bf%
Corollary.}\ }\ignorespaces}{ \end{theorem} }
\newenvironment{thm*}{\begin{theorem}{\em{\bf%
Theorem*.}\ }\ignorespaces}{\end{theorem}}
\newenvironment{prop*}{\begin{theorem}{\em{\bf%
Proposition*.}\ }\ignorespaces}{ \end{theorem} }
\newenvironment{lemma*}{\begin{theorem}{\em{\bf%
Lemma*.}\ }\ignorespaces}{ \end{theorem} }
\newenvironment{cor*}{\begin{theorem}{\em{\bf%
Corollary*.}\ }\ignorespaces}{ \end{theorem} }
\newenvironment{defn}{\begin{theorem}\em{{\bf%
Definition.}\ }\ignorespaces}{ \end{theorem} }

\newenvironment{conj}{\begin{theorem}{\em{\bf%
Conjecture.}\ }\ignorespaces}{ \end{theorem} }

\newenvironment{rk}{\begin{theorem}\em{{\bf%
Remark.}\ }\ignorespaces}{ \end{theorem} }
\newenvironment{rks}{\begin{theorem}\em{{\bf%
Remarks.}\ }\ignorespaces}{ \end{theorem} }

\newenvironment{exs}{\vskip1em\noindent{\bf%
Examples}.\begin{enumerate}}{\end{enumerate}}

\newenvironment{nota}{\vskip1em\noindent{\bf Notation}.}{}

\newcounter{spec}
\newenvironment{thlist}{\begin{list}{\rm{(\roman{spec})}}%
{\usecounter{spec}\labelwidth=20pt\itemindent=0pt\labelsep=10pt}}%
{\end{list}}
\renewcommand{\rm}{\normalshape}

\newcommand{\oo}{\displaystyle\operatornamewithlimits\otimes}
\newcommand{\af}{{\Bbb A}^1}
\newcommand{\tors}{{\operatorname{tors}}}
\newcommand{\cotors}{{\operatorname{cotors}}}

\newcommand{\Pic}{\operatorname{Pic}}

\newcommand{\Ker}{\operatorname{Ker}}
\newcommand{\Coker}{\operatorname{Coker}}
\newcommand{\IM}{\operatorname{Im}}
\newcommand{\car}{\operatorname{char}}

\newcommand{\Div}{\underline{\operatorname{Div}}}
\newcommand{\Codiv}{\underline{\operatorname{Codiv}}}
\renewcommand{\lim}{\varprojlim}
\newcommand{\colim}{\varinjlim}

\newcommand{\Spec}{\operatorname{Spec}}

\newcommand{\cd}{\operatorname{cd}}
\newcommand{\et}{{\operatorname{\acute{e}t}}}
\newcommand{\Et}{{\operatorname{\acute{E}t}}}
\newcommand{\cont}{{\operatorname{cont}}}

\newcommand{\ind}{{\operatorname{ind}}}
\newcommand{\ladic}{{\operatorname{\text{$l$}-adic}}}
\newcommand{\ord}{{\operatorname{ord}}}
\newcommand{\Zar}{{\operatorname{Zar}}}
\newcommand{\gr}{\operatorname{gr}}
\newcommand{\rank}{\operatorname{rank}}
\newcommand{\corank}{\operatorname{corank}}
\newcommand{\trdeg}{\operatorname{trdeg}}
\newcommand{\inj}{\hookrightarrow}
\newcommand{\Inj}{\lhook\joinrel\longrightarrow}
\newcommand{\surj}{\rightarrow\!\!\!\!\!\rightarrow}
\newcommand{\Surj}{\relbar\joinrel\surj}
\newcommand{\resp}{{\it resp.} }
\newcommand{\ie}{{\it i.e.} }
\newcommand{\eg}{{\it e.g.} }
\newcommand{\cf}{{\it cf.} }
\newcommand{\2}{_{(2)}}
\newcommand{\iso}{@>\sim>>}

\newcommand{\F}{{\Bbb{F}}}

\newcommand{\Q}{{\Bbb{Q}}}

\newcommand{\Z}{{\Bbb{Z}}}
\newcommand{\sA}{{\cal A}}
\newcommand{\sB}{{\cal B}}
\newcommand{\sC}{{\cal C}}
\newcommand{\sD}{{\cal D}}
\newcommand{\sF}{{\cal F}}
\newcommand{\sG}{{\cal G}}
\newcommand{\sH}{{\cal H}}
\newcommand{\sK}{{\cal K}}
\newcommand{\sM}{{\cal M}}
\newcommand{\sP}{{\cal P}}
\newcommand{\sO}{{\cal O}}
\newcommand{\sS}{{\cal S}}
\newcommand{\sX}{{\cal X}}
\newcommand{\bC}{{\Bbb{C}}}
\newcommand{\bF}{{\Bbb{F}}}
\newcommand{\bG}{{\Bbb{G}}}
\newcommand{\bH}{{\Bbb{H}}}
\newcommand{\bN}{{\Bbb{N}}}
\newcommand{\bP}{{\Bbb{P}}}
\newcommand{\bQ}{{\Bbb{Q}}}
\newcommand{\bR}{{\Bbb{R}}}

\newcommand{\bZ}{{\Bbb{Z}}}
\newcommand{\prf}{\noindent{\bf Proof}. }
\renewcommand{\qed}{\hfill\text{$\Box$}}

\begin{document}
\title[Tate conjecture]{A sheaf-theoretic reformulation of the Tate conjecture}
\author{Bruno Kahn}
\address{Institut de Math\'ematiques de Jussieu\\Equipe Th\'eories
G\'eom\'etriques\\Universit\'e Paris 7\\Case 7012\\75251 Paris Cedex
05\\France}
\email{kahn@@math.jussieu.fr}
\maketitle
\hfill 18-12-1997 \\

Preliminary version

\tableofcontents

\begin{abstract} We present a conjecture which unifies several conjectures on
motives in characteristic $p$.
\end{abstract}

\section*{Introduction}

In  a talk at the June 1996 Oberwolfach algebraic $K$-theory conference, I
explained that, in view of Voevodsky's proof of the Milnor conjecture, the Bass
conjecture on finite generation of $K$-groups of regular schemes of finite type
over the integers implies the rational Beilinson-Soul\'e conjecture on vanishing of
algebraic $K$-theory of low weights in characteristic $\neq 2$. This argument is
reproduced in the appendix. The next day, Thomas Geisser explained that, in
characteristic
$>0$, the Tate conjecture on surjectivity of the
$\bQ_l$-adic cycle map for smooth, projective varieties over a finite field,
together with the conjecture that rational and numerical equivalences agree for
such varieties, also implies the Beilinson-Soul\'e conjecture. His work is now
available in
\cite{geisser}. The present paper stems from an attempt to understand the
relationship between these two facts.

Let $l$ be a prime number. We present a conjecture (conjecture
\ref{co7.1}) which is equivalent to the conjunction of three well-known conjectures
on smooth, projective varieties $X$ over $\F_p$ ($p\neq l$):

\begin{enumerate}
\item {\it Tate's conjecture}: the geometric cycle map
\begin{equation*}\tag{*}\label{cycle}
CH^n(X)\otimes \Q_l\to H^{2n}(\overline X,\Q_l(n))^G
\end{equation*}
is surjective ($\overline X=X\times_{\F_p}\bar\F_p$, $G=Gal(\bar\F_p/\F_p)$).
\item {\it Partial semi-simplicity}: the characteristic subspace of $H^n(\overline
X,\Q_l(n))$ corresponding to the eigenvalue $1$ of Frobenius is semi-simple.
\item {\it rational equivalence =  homological equivalence}: the map \eqref{cycle}
is injective.
\end{enumerate}

We note that the conjunction of (1), (2) and (3) (in codimensions $n$ and $\dim
X-n$) is equivalent to the so-called strong form of the Tate conjecture: the order
of the pole of the Hasse-Weil zeta function $\zeta(X,s)$ at $s=n$ equals the rank of
the group $A^n(X)$ of codimension $n$ cycles modulo numerical equivalence, \cf
\cite[th. 2.9]{tate2}. In this light, conjecture \ref{co7.1} below can really be
thought of as a sheaf-theoretic reformulation of the (strong) Tate conjecture.

The main result of this paper is that conjecture \ref{co7.1} (and therefore the
three above together) also implies a host of other familiar conjectures:

\begin{itemize}
\item  The Tate conjecture and the partial semi-simplicity conjecture (see
definition \ref{d4.2}) for smooth, projective varieties over {\it any}
finitely generated field of characteristic $p$ (theorem
\ref{t7.2}).
\item The injectivity of the (arithmetic) $\bQ_l$-adic cycle map for {\it any}
smooth variety $X$ over $\F_p$, and a description of its
image. In particular, that rational equivalence equals homological
equivalence on $X$ (proposition \ref{p7.3} and corollary \ref{c7.2}).
\item The Bass conjecture (finite generation of the algebraic $K$-groups) after
tensoring by
$\bQ$ for such
$X$ (corollary
\ref{c7.4}).
\item The rational Beilinson-Soul\'e conjecture for such $X$ ({\it ibid.}).
\item The vanishing of $K_n^M(F)\otimes\bQ$ for $n>\trdeg(F/\bF_p)$, when
$F$ is a field of characteristic $p$ (Bass-Tate conjecture in characteristic $p$,
\cf corollary
\ref{p7.4}).
\item Soul\'e's conjecture on the order of the zeroes and poles of the zeta
function of a variety over $\bF_p$ at integers \cite{souleicm} (in fact a refined
version of it, th. \ref{t8.5}).
\item The existence of a canonical integral structure on
continuous \'etale cohomology (proposition \ref{p7.2}).
\end{itemize}

It is known that these conjectures imply other ones: for example, the Tate
conjecture plus semi-simplicity (in a range) imply that homological and
numerical equivalences agree (\cite[prop. 8.4]{milne},
\cite[(2.6)]{tate2}), that the K\"unneth components of the diagonal are
algebraic \cite[\S 3]{tate2} and that Hard Lefschetz holds for cycles
modulo numerical equivalence (see theorem \ref{t7.2} b)). On the other
hand, the injectivity of the cycle map implies the existence of
Beilinson's conjectural filtration on Chow groups (a variant of
\cite[lemma 2.7]{jannsen4}, see also \cite[lemma 2.2]{jannsen4}, \cf
theorems \ref{t7.2} b) and \ref{t7.4} here), which is in turn
equivalent to Murre's idempotents conjecture (\cite[1.4]{murre}, \cite[th.
5.2]{jannsen4},
\cf corollary \ref{c7.7} here).

One conjecture we don't know to
follow from conjecture \ref{co7.1} is the Hodge-type standard conjecture
in characteristic $p$ \cite[\S 5]{kleiman}. This is perhaps moral, since this
conjecture implies the semi-simplicity of the full Frobenius action on cohomology
groups, not just at the eigenvalue $1$ ({\it ibid.}, th. 5-6 (2)); we feel that
conjecture
\ref{co7.1} is not strong enough to imply such a fact {\it per se}. 

There is a variant of conjecture \ref{co7.1}, involving Voevodsky's
motivic cohomology (conjecture \ref{co1}). We show that, under resolution
of singularities, it is equivalent to conjecture \ref{co7.1}. It implies
Lichtenbaum's conjectures of \cite[\S 7]{licht} on the nature of the (\'etale)
motivic cohomology of a smooth projective variety over $\F_p$ and on
values of its zeta function at nonnegative integers; in fact, we extend
this consequence to arbitrary smooth varieties (proposition \ref{c5.4} and
theorem \ref{t4.4}). Under resolution of singularities, these properties
therefore follow from the Tate conjecture, the partial semi-simplicity
conjecture and the rational equivalence=homological equivalence
conjecture. 

Let us now describe conjectures \ref{co7.1} and \ref{co1}. In
\cite{jannsen}, Jannsen defines continuous \'etale cohomology as the
derived functors of the composition of left exact functors
\[\begin{CD}\sA^\bN@>\Gamma>> Ab^\bN\\ && @V{\lim}VV\\ && Ab\end{CD}\]
where $\sA$ is the category of abelian sheaves over the small \'etale site
of some scheme $X$. Our starting point, only briefly considered in
\cite[p. 219]{jannsen}, is to observe that one can equally well define
this composition by going the other way:  
\[\begin{CD}\sA^\bN\\@V{\lim}VV\\ \sA@>\Gamma>> Ab.\end{CD}\]

In general the category $\sA$ does not satisfy Grothendieck's axiom $AB4^*$
(products are not exact), but nevertheless the inverse limit functor is left
exact and preserves injectives, so it makes sense to study its higher derived
functors with values in
$\sA$ or, better, its total derived functor $R\lim$ with values in the derived
category $\sD^+(\sA)$. In particular, we can define
\[\bZ_l(n)^c_X:=R\lim\mu_{l^\nu}^{\otimes n}\in\sD^+(\sA)\]
provided the prime $l$ is invertible on $X$, a blanket assumption here.
Note that $\bZ_l(n)^c_X$ is not an $l$-adic sheaf in the sense of SGA5,
but a complex of honest sheaves (up to quasi-isomorphism).

It is more convenient to work at once over the big \'etale site of $\Spec
\bF_p$ (restricted to smooth schemes): denote by $\bZ_l(n)^c$ the
corresponding object. We prove in section \ref{5} that
$\bZ_l(n)^c=\bQ_l/\bZ_l(n)[-1]$ for $n<0$, that
\[\bZ_l(0)^c\simeq\bZ^c\otimes\bZ_l\] where $\bZ^c$ is the class of a certain
explicit complex of length $1$ coming from the small \'etale site of
$\Spec\bF_p$, such that
\[\bQ^c:=\bZ^c\otimes\bQ\simeq \bQ[0]\oplus\bQ[-1],\]
and that $\bH^i(X,\bZ_l(n)^c)$ is finite for $i\notin
[n,2n+1]$, $n>0$ and $X$ smooth over $\bF_p$. 

In section \ref{coarse}, we define a map
\[\bQ_l(0)^c\oo^L \alpha^*\sK_n^M[-n]\to
\bQ_l(n)^c\]
where $\sK_n^M$ is the sheaf of Milnor $K$ groups over the big
smooth Zariski site of $\Spec \F_p$ and $\alpha$ is the projection of the
big \'etale site onto the big Zariski site. Conjecture \ref{co7.1} says
that this map is an {\it isomorphism}.

In section \ref{6} we analogously define a map
\[\bZ_l(0)^c\oo^L\alpha^*\bZ(n)\to \bZ_l(n)^c.\]

Here, $\bZ(n)$ is the $n$-th motivic complex of Suslin-Voevodsky \cite{suvo} and
$\alpha$ is as above. Conjecture \ref{co1} also states that this
map is an isomorphism.

This paper is organised as follows. In section \ref{1} we go through basic
facts concerning inverse limits and their higher derived functors. In sections
\ref{2} and \ref{classical} we ``recall" some properties of continuous \'etale
cohomology, cohomology with proper supports and Borel-Moore homology. In section
\ref{3}, we do some computations over finite fields. In section \ref{4} we develop
some elementary useful formalism. In section \ref{5} we give our main
unconditional results. In section \ref{zeta}, we generalise the main result of
Milne \cite[th. 0.1]{milne} on the principal part of the zeta function $\zeta(X,s)$
at $s=n$ from smooth, projective varieties to arbitrary varieties over
$\bF_p$ (theorem \ref{t4.3}, theorem \ref{t6.1} and corollary \ref{c4.1});
I am particularly fond of the statement in theorem \ref{t6.1} b). In section
\ref{coarse}, we present our main conjecture and derive the consequences listed
above, and others. Perhaps the most important ones are contained in theorems
\ref{t7.2} and \ref{t7.4}. In section
\ref{6} we present the motivic variant of our main conjecture and explore some
further consequences of it. In particular, in proposition \ref{p5.7} we give as a
consequence of conjecture
\ref{co1} (and resolution of singularities) a description of the kernel and
cokernel of the integral cycle map $CH^n(X)\otimes\bZ_l\to
H^{2n}_\cont(X,\bZ_l(n))$. In section
\ref{one}, we prove the restriction of conjecture
\ref{co7.1} to curves over $\bF_p$. In section \ref{zariski} we propose a 
``Zariski" variant of conjecture \ref{co1} related to the Kato conjecture. In
section
\ref{8} we briefly discuss the much more open case of schemes of finite type
over
$\Spec\bZ[1/l]$. Here it is clear that $l$-adic cohomology does not give the full
picture; in order to explain the contribution of archimedean places, we expect that
Arakelov geometry will play its r\^ole. Finally, in the appendix, we give a proof
of our announcement that the Bass conjecture implies the Beilinson-Soul\'e
conjecture, in all characteristics.

A pendant of the present investigation for $l=p$ could doubtlessly be
developed; we don't tackle this here. Another thing which should definitely be
done is to investigate analogues of conjectures \ref{co7.1} and \ref{co1} over
$p$-adic fields, $p\neq l$ and $p=l$.

It will be obvious that this paper owes much to Uwe Jannsen's previous work on
motives. In particular, it is by reading the last chapter of \cite{jannsen2}, where
Jannsen extends the Tate and Hodge conjectures from smooth projective varieties to
arbitrary varieties, that I got the courage to look for a sheaf-theoretic version
of the Tate conjecture. But it would be unfair not to acknowledge the influence of
Vladimir Voevodsky's work as well. I also wish to thank Bernhard
Keller and Amnon Neeman for a discussion which helped clarify some confusion
related to section
\ref{3}. Finally, I would like to thank Thomas Geisser for a large number of helpful
comments.

\section{Continuous cohomology}\label{1}

\begin{nota} If $\sA$ is an abelian category and $\sD(\sA)$ is its derived
category, we usually denote by $A\mapsto A[0]$ the canonical functor $\sA\to
\sD(\sA)$.
\end{nota}

\subsection{Inverse limits}\label{1.1}

Let $\sA$ be an abelian category, and let $\sA^\bN$ be the category of functors
from $\bN$ to $\sA$, where the category structure on $\bN$ is given by its ordering.
An object in $\sA^\bN$ is the same as a projective system of objects of $\sA$.

The constant functor $\sA\to \sA^\bN$ (induced by the projection $\bN\to \{1\}$) is
exact. If $\sA$ is complete ($\iff$ $AB3^*$ holds), this functor has as a left
adjoint the left exact functor:
\begin{align*}
\lim:\sA^\bN&\to \sA\\
(A_n)&\mapsto \lim A_n.
\end{align*}
which carries injectives to injectives.

By
\cite[(1.1)]{jannsen}, $\sA^\bN$ has enough injectives if and only if $\sA$ has.
In this case, we can derive $\lim$ in 
\[\lim^{i}:\sA^\bN\to \sA\qquad (i\ge 0).\]

More globally, we can
define a total derived functor
\[R\lim:\sD^+(\sA^\bN)\to \sD^+(\sA).\]

\subsection{Extending functors in one variable}\label{1.3} Let $\sB$ be another
complete abelian category with enough injectives and
$T:\sA\to
\sB$ a left exact functor. Then $T$ induces a left exact functor from
$\sA^\bN$ to $\sB^\bN$, that we still denote by $T$. We have a natural
transformation of left exact functors
\begin{gather}\label{eq1.4}
T\circ \varprojlim\to \varprojlim\circ T\\
\intertext{hence a chain of natural transformations in $\sD(\sB)$}
\label{eq1.1}
RT\circ R\lim @<<< R(T\circ\lim)@>>> R(\lim\circ T)@>>> R\lim\circ RT.
\end{gather}

\begin{lemma}\label{l1.3}
a) The left transformation in \eqref{eq1.1} is a natural isomorphism, hence a
spectral sequence
\[R^pT(\lim^{q} A)\Rightarrow R^{p+q}(T\circ \lim)(A)
\]
for any $A\in\sA^\bN$.\\
b)  If $T$ carries
injectives to
$\lim$-acyclics, the right transformation in \eqref{eq1.1} is a natural
isomorphism, hence a spectral sequence
\[\lim^{p}R^qT(A)\Rightarrow R^{p+q}(\lim\circ T)(A).
\]
c) If $T$ commutes with products, then \eqref{eq1.4} and  the middle 
transformation in \eqref{eq1.1} are natural isomorphisms.\\
d) If $T$
has a left adjoint which respects monomorphisms, then \eqref{eq1.4} and all
transformations in \eqref{eq1.1} are natural isomorphisms.
\end{lemma}

\begin{pf} a) holds because $\lim$ carries injectives to injectives. b) and c) are
clear; the assumption in d) implies
those  of b) and c).\end{pf}

\subsection{Extending functors in two variables}\label{1.4} Let $\sA,\sB,\sC$ be
three abelian categories and let
$T:\sA\times
\sB\to
\sC$ be an biadditive, bicovariant functor. Then $T$ extends naturally to a functor
\[T:\sA^\bN\times\sB^\bN\to \sC^\bN\]
by the formula
\[T((a_\nu),(b_\nu))=(T(a_\nu,b_\nu)).\]

This formula can be derived in the usual way. For example, suppose that
$\sA=\sB=\sC$ is endowed with a right exact tensor product
\[\otimes=\sA\times \sA\to \sA\]
for which $\sA$ has enough acyclic objects. Suppose moreover that $\otimes$ has
finite homological dimension. One can then derive the natural transformation
\[\otimes \circ (\lim,\lim)\to \lim\circ \otimes\]
between functors from $\sA^\bN\times \sA^\bN$ to $\sA$ into a natural transformation
\begin{equation}\label{eq1.6}\oo^L\circ(R\lim,R\lim)\to
R\lim\circ\oo^L\end{equation}
between functors from $\sD^+(\sA^\bN)\times
\sD^+(\sA^\bN)$ to $\sD^+(\sA)$.

If $T$ is contravariant in $\sA$ and covariant in $\sB$, like $Hom$, its extension
should be defined by an end, as in \cite[ch. IX \S 5]{mcl}.

\subsection{Sheaves}\label{1.2}

Let $X$ be a site and $R$ a ring. We consider the category $\sA$ of sheaves of
$R$-modules on
$X$. Since the functor ``global sections" $\Gamma$ has as a left adjoint the exact
functor ``constant sheaf",  we can apply lemma
\ref{l1.3}. As in \cite{jannsen}, we define

\begin{defn} Let $(\sF_\nu)$ be a projective system of abelian sheaves over $X$. The
{\em continuous cohomology} of $X$ with values in $(\sF_\nu)$ is
\[H^n_\cont(X,(\sF_\nu))=R^n(\Gamma\circ\lim)((\sF_\nu))=
R^n(\lim\circ\Gamma)((\sF_\nu))\]
(see lemma \ref{l1.3} c)).
\end{defn}

Since products are exact in $R-mod$, $\lim^i=0$ for $i>1$ and the spectral
sequence of lemma \ref{l1.3} b) yields ``Milnor" exact sequences
\cite[(1.6)]{jannsen}
\begin{equation}\label{eq1.2} 0\to\lim^1H^{n-1}(X,\sF_\nu)\to
H^n_\cont(X,(\sF_\nu))\to
\lim H^n(X,\sF_\nu)\to 0.\end{equation}

But the spectral sequence of lemma \ref{l1.3} a) also says that continuous \'etale
cohomology can be computed by a hypercohomology spectral sequence:
\[H^p(X,\lim^q(\sF_\nu))\Rightarrow H^{p+q}_\cont(X,(\sF_\nu)).\]

This gives us a formula for the {\it stalks} of $\lim^q(\sF_\nu)$ (\cf
\cite[(3.12)]{jannsen}):

\begin{lemma} \label{l1.1} For $x$ a point of $X$, we have
\[(\lim^q(\sF_\nu))_x=\colim_{U\ni x} H^{q}_\cont(U,(\sF_\nu)).\qquad\qed\]
\end{lemma}

We also have:

\begin{prop}\label{p1.2} Suppose that $\cd_R(X)\le n$. Then,\\
a) For any inverse system of
$R$-sheaves $(\sF_\nu)$ and any $U\in X$, we have
$H^q_\cont(U,(\sF_\nu))=\lim^q(\sF_\nu)=0$ for
$q>n+1$.\\
b) If moreover $(\sF_\nu)$ satisfies the Mittag-Leffler condition, then 
$H^q_\cont(U,(\sF_\nu))=\lim^q(\sF_\nu)=0$ for
$q>n$.
\end{prop}

\begin{pf} a) This follows from the Milnor exact
sequence \eqref{eq1.2}.

b) By \cite[(1.14)]{jannsen}, we may assume that the transition morphisms of
$(\sF_\nu)$ are surjective. Then the same holds for the inverse system
$(H^n_\cont(U,\sF_\nu))$. The conclusion now follows again from
\eqref{eq1.2}.\end{pf}

Let $f:Y\to X$ be a morphism of sites. By lemma \ref{l1.3} c), the natural
transformation
\begin{equation}\label{eq1.5} f_*\circ\lim\to \lim\circ f_*\end{equation}
is an isomorphism; denote by
$f_*^\cont$ either side of \eqref{eq1.5}. Suppose that $f^*$ respects monomorphisms
(\eg is left exact). Applying lemma \ref{l1.3} d), we get
similarly two spectral sequences
\[\lim^p(R^qf_* \sF_n)\Rightarrow R^{p+q}f_*^\cont(\sF_n)\Leftarrow
R^pf_*\lim^q(\sF_n).\]

In general, the natural transformation
\[f^*\circ\lim\to \lim\circ f^*\]
is not an isomorphism. When $f^*$ has a left adjoint which respects
monomorphisms, we can however apply lemma \ref{l1.3} d). This is the case for
example if the associated topos $Shv(Y)$ is induced from $Shv(X)$ via
$f$ (\cf
\cite[expos\'e IV, \S 5]{sga4}).

Suppose $R$ commutative of finite Tor-dimension. Let $(\sF_\nu),(\sG_\nu)$ be
two inverse systems of sheaves. Applying
\eqref{eq1.6}, we get a morphism in $\sD^+(\sA)$
\[R\lim(\sF_\nu)\oo^LR\lim(\sG_\nu)\to R\lim(\sF_\nu[0]\oo^L\sG_\nu[0]).\]

Let $(\sH_\nu)$ be a third inverse system of sheaves and
$\alpha: (\sF_\nu)\times(\sG_\nu)\to (\sH_\nu)$ a bilinear map in $\sA^\bN$. Using
the composite
\[\sF_\nu[0]\oo^L\sG_\nu[0]\to (\sF_\nu\otimes\sG_\nu)[0]@>\alpha[0]>> \sH_\nu[0]\]
in $\sD^+(\sA)$, we get an induced morphism
\begin{equation}\label{eq1.3} \alpha:R\lim(\sF_\nu)\oo^LR\lim(\sG_\nu)\to
R\lim(\sH_\nu).\end{equation}

\section{Continuous \'etale cohomology}\label{2}

\begin{defn} Let $l$ be a prime number and $S$ be a scheme over $\Spec \bZ[1/l]$.
Let
$S_\Et$ denote the big \'etale site of
$S$, with underlying category the category of all $S$-schemes locally of finite
type. For
$n\in\bZ$, we denote by $\bZ_l(n)^c_S$ (\resp $\bQ_l(n)^c_S$) the object
$R\lim(\bZ/l^\nu(n))$ (\resp $\bZ_l(n)^c\otimes \bQ$) of
$\sD^+(Ab(S_\Et))$. Here, $\bZ/l^\nu(n)$ denotes the sheaf of $l^\nu$-th roots of
unity twisted $n$ times.
\end{defn}

For $X$ of locally finite type over $S$, we shall usually denote by
$H^i_\cont(X,\bZ_l(n))$ and $H^i_\cont(X,\bQ_l(n))$ the hypercohomology groups
$\bH^i_\et(X,\bZ_l(n)^c_S)$ and $\bH^i_\et(X,\bQ_l(n)^c_S)$.

Note that, by definition, $\bZ_l(n)^c_S$ and $\bQ_l(n)^c_S$ are complexes of honest
sheaves over
$S_\Et$ (up to quasi-isomorphism). By subsection
\ref{1.2}, their hypercohomology coincides with Jannsen's continuous
\'etale cohomology \cite{jannsen}. They are related by the obvious

\begin{lemma}\label{l2.1} There are exact triangles in $\sD^+(Ab(S_\Et))$:
\[\begin{CD}
\bZ_l(n)^c@>l^\nu>> \bZ_l(n)^c@>>> \bZ/l^\nu(n)[0]@>>> \bZ_l(n)^c[1] \\
\bZ_l(n)^c@>>> \bQ_l(n)^c@>>> \bQ_l/\bZ_l(n)[0]@>>> \bZ_l(n)^c[1] 
\end{CD}\]
where $\bQ_l/\bZ_l(n):=\colim \bZ/l^\nu(n)$. \qed
\end{lemma}

Let $X@>f>> S$ be a morphism. Since $f^*\bZ/l^\nu(n)=\bZ/l^\nu(n)$, we have
morphisms in $\sD^+(Ab(X_\Et))$
\begin{align*}
f^*\bZ_l(n)^c_S&\to \bZ_l(n)^c_X\\
f^*\bQ_l(n)^c_S&\to \bQ_l(n)^c_X.
\end{align*}

If $f$ is locally of finite type, then the topos $Shv(X_\Et)$ is
induced \cite[expos\'e IV, \S 5]{sga4} from
$Shv(S_\Et)$ via $f$ and these are isomorphisms by the remarks after proposition
\ref{p1.2}. This justifies to write
$H^*_\cont(X)$ rather than $H^*_\cont(X/S)$. In other cases, as a rule they are
(violently) {\it not} isomorphisms. For example, if
$S=\Spec k$,
$k$ a field,
$X=\Spec \bar k$,
$\bar k$ its algebraic closure, then $\bQ_l(n)^c_X=0$ while $f^*\bQ_l(n)^c_S$ is
usually nonzero. However, lemma \ref{l2.1} implies

\begin{lemma}\label{l2.0} The cone of $f^*\bZ_l(n)^c_S\to \bZ_l(n)^c_X$ has
uniquely divisible cohomology sheaves.\qed
\end{lemma}

In general, we shall denote by
\[\tilde H^*_\cont(X/S,\bZ_l(n))\]
the cohomology groups $H^*_\et(X,f^*\bZ_l(n)^c_S)$. By the above, they coincide
with $H^*_\cont(X,\bZ_l(n))$ when $f$ is locally of finite type, and
map naturally to them in any case. The same holds with $\bQ_l$ coefficients. If $S$
is affine Noetherian and
$f$ is quasi-compact and quasi-separated, then, by \cite[C.9]{TT}, $X$ is the
filtering inverse limit of a system of
$S$-schemes $X_i$ of finite presentation, with affine transition maps. Since
\'etale cohomology commutes with this type of inverse limits \cite[Expos\'e
VII, cor. 5.8]{sga4}, we have
\[\tilde H^*_\cont(X/S,\bZ_l(n))=\colim H^*(X_i,\bZ_l(n)).\]

Therefore, the groups $\tilde H^*_\cont(X/S,\bZ_l(n))$ generalise those introduced
by Jannsen in \cite[\S 11]{jannsen2}. In particular, using lemma \ref{l1.1}, we can
describe the stalks of $\bZ_l(n)^c_S$ at a geometric $S$-point $s$:
\begin{equation}\label{eq2.1}
\sH^q(\bZ_l(n)^c_S)_s=\tilde H^q_\cont(\Spec \sO_{S,s}^{sh}/S,\bZ_l(n)).
\end{equation}

Let
$m,n\in\bZ$. From
\eqref{eq1.3} we deduce canonical pairings:
\begin{equation}\label{eq2.2}
\begin{CD}
\bZ_l(m)^c_S\oo^L\bZ_l(n)^c_S@>>> \bZ_l(m+n)^c_S\\
\bQ_l(m)^c_S\oo^L\bQ_l(n)^c_S@>>> \bQ_l(m+n)^c_S.
\end{CD}
\end{equation}

\subsection*{$l$-adic sheaves} For any $\bZ[1/l]$-scheme $X$, let $\bZ_l(X_\et)$ be
the full subcategory of the category $Ab(X_\et)^\bN$ ({\em small} \'etale site)
consisting of AR-$\bZ_l$-constructible sheaves in the sense of \cite[expos\'e VI,
d\'ef. 1.5.1]{sga5} (this category is denoted by $U(X)$ in {\it loc. cit.}, not.
1.5.4). This is an abelian subcategory. Suppose $X$ connected, let $x$ be a
geometric point of $X$ and let $\Pi=\pi_1(X,x)$ be its \'etale fundamental group
at $x$. Let $\bZ_l(\Pi)_{fg}$ be the category of continuous representations of
$\Pi$ onto finitely generated $\bZ_l$-modules. Then 
\[L\mapsto (L/l^\nu)\]
defines a functor  $\bZ_l(\Pi)_{fg}\to \bZ_l(X_\et)$. We have

\begin{prop} \label{p2.5} \cite[expos\'e VI, prop. 1.2.5 and 1.5.5]{sga5} This
functor is a full embedding; moreover, any object of $\bZ_l(X_\et)$ is
AR-equivalent to the image of an object of $\bZ_l(\Pi)_{fg}$, in the sense of
\cite[expos\'e V, 2.2]{sga5}.\qed
\end{prop}

Ekedahl \cite{ek} has constructed 
a triangulated category $\sD_{\bZ_l}(X_\et)$ with a $t$-structure such that
\begin{itemize}
\item $\bZ_l(X_\et)$ is the heart of $\sD_{\bZ_l}(X_\et)$ for this $t$-structure.
\item There is a triangulated functor
\[\sD_{\bZ_l}(X_\et)\to \sD^+(Ab(X_\et)^\bN)\]
which respects the $t$-structures and whose restriction to hearts coincides with
the full embedding $\bZ_l(X_\et)\inj Ab(X_\et)^\bN$.
\end{itemize}

This is not completely true, as Ekedahl works modulo essential (ML) isomorphisms.
However one can replace his construction by that where one keeps essential
isomorphisms uninverted, so that $\sD_{\bZ_l}(X_\et)$ actually maps to
$\sD^+(Ab(X_\et)^\bN)$.

For all $n\in\bZ$, the projective system of sheaves $\bZ_l(n)=(\bZ/l^\nu(n))$
belongs to
$\bZ_l(X_\et)$. By definition, we have
\[\bZ_l(n)^c_{|X}=R\lim\bZ_l(n)\]
where $_{|X}$ denotes restriction to the small \'etale site of $X$.

Although we need $\bZ_l(n)^c$ to express our conjectures, the reader should be
aware that all serious computations are actually done with $l$-adic sheaves,
notably with the help of proposition \ref{p2.5}. Also beware of similar-looking
spectral sequences which are in fact different. For example, if $f:Y\to X$ is
a morphism, there are two spectral sequences
\[H^p_\cont(X,R^qf^\ladic_*\bZ_l(n))\Rightarrow
H^{p+q}_\cont(Y,\bZ_l(n))\Leftarrow  H^p_\et(X,R^qf_*\bZ_l(n)^c).\]

Here we have written $R^qf^\ladic_*$ for the higher derived images
of $f_*$ from
$\bZ_l(Y_\et)\to \bZ_l(X_\et)$. They are best explained by a naturally commutative
diagram of functors
\[\begin{CD}
\bZ_l(n)_Y\in&\sD_{\bZ_l}(Y_\et)@>>> \sD^+(Ab(Y_\et)^\bN)@>R\lim>>
\sD^+(Ab(Y_\et))&\ni\bZ_l(n)^c_Y \\
&@V{Rf_*^\ladic}VV @V{Rf_*}VV @V{Rf_*}VV\\
&\sD_{\bZ_l}(X_\et)@>>> \sD^+(Ab(X_\et)^\bN)@>R\lim>> \sD^+(Ab(Y_\et))
\end{CD}\]

The left spectral sequence is obtained from the composition of $Rf_*^{\text{\rm
$l$-adic}}$ with the bottom horizontal composition, applied to $\bZ_l(n)_Y$, while
the right one is simply the Leray spectral sequence for the right $Rf_*$, applied to
$\bZ_l(n)^c_Y$. These two spectral sequences will generally yield different
filtrations on $H^*_\cont(Y,\bZ_l(n))$. In practice it is the $l$-adic spectral
sequence that is the most useful.

\section{Classical theorems}\label{classical} This section is included for future
reference and can be skipped at first reading. We give a few basic results on
continuous
\'etale cohomology, which directly involve $\bZ_l(n)^c$ or variants of it. One
notable exception is Poincar\'e duality, for which we are not aware of a
formulation involving this object. 

\begin{nota} Let $\pi_X:X_\Et\to X_\et$ be the projection of the big \'etale site
of $X$ onto its small \'etale site. If $C\in \sD^+(X_\Et)$, we write
\[C_{|X}:=R(\pi_X)_*C\]
for the restriction of $C$ to $X_\et$. We simply write
\begin{align*}
\bZ_l(n)^c_{|X}&=(\bZ_l(n)^c_X)_{|X}\\
\bQ_l(n)^c_{|X}&=(\bQ_l(n)^c_X)_{|X}.
\end{align*}
\end{nota}

\begin{thm}{\rm (Proper base change)}\label{t2.-1}  Let
\[\begin{CD}
Y'@>g'>> Y\\
@V{f'}VV @V{f}VV\\
X'@>g>> X
\end{CD}\]
be a cartesian diagram of schemes, where $f$ is proper and $g$  is
locally of finite type. Then, for all $n\in\Z$, the base change morphism
\[g^*Rf_*\bZ_l(n)^c_Y\to Rf'_*g'{}^*\bZ_l(n)^c_Y\]
in $\sD^+(X'_\Et)$ is an isomorphism.
\end{thm}

\begin{pf} This follows from applying $R\lim$ to the classical proper base change
theorem \cite[Expos\'e XII, th. 5.1]{sga4}
\[g^*Rf_*\bZ/l^\nu(n)\iso Rf'_*g'{}^*\bZ/l^\nu(n)\]
extended to the big \'etale sites.
\end{pf}

\begin{cor}\label{c2.1} Let $f:\tilde X\to X$ be a
proper morphism and
$U\subset X$ be an open subset such that $f_{|f^{-1}(U)}:f^{-1}(U)\to U$ is an
isomorphism. Let $Z=X\setminus U$ and $\tilde Z=f^{-1}(Z)$, both with their
reduced structure. Then there is an exact triangle
\[\bZ_l(n)_X^c@>>>i_*\bZ_l(n)_Z^c\oplus Rf_*\bZ_l(n)_{\tilde X}^c@>>>
\tilde{\text{\i}}_*\bZ_l(n)_{\tilde Z}^c@>>>\bZ_l(n)_X^c[1]\]
where $i$ (\resp $\tilde{\text{\i}}$ is the closed immersion $Z\inj X$ (\resp
$\tilde Z\inj \tilde X$). Hence a long exact sequence
\begin{multline*}
\dots\to H^{i-1}_\cont(\tilde Z,\bZ_l(n))\to
H^i_\cont(X,\bZ_l(n))\\
\to H^i_\cont(Z,\bZ_l(n))\oplus H^i_\cont(\tilde
X,\bZ_l(n))\to H^i_\cont(\tilde Z,\bZ_l(n))\to\dots\end{multline*}
\end{cor}

\begin{pf} (compare \eg \cite[proof of prop. 2.1]{pw}) Let $j:U\to X$ be the open
immersion complementary to $i$. Apply the exact triangle of functors
$j_!j^*\to Id\to i_*i^*\to j_!j^*[1]$ to the morphism $\bZ_l(n)_X^c\to
Rf_*\bZ_l(n)_{\tilde X}^c$ to get a commutative diagram of exact triangles
\[\begin{CD}
j_!j^*\bZ_l(n)_X^c@>>> \bZ_l(n)_X^c@>>>
i_*i^*\bZ_l(n)_X^c@>>>j_!j^*\bZ_l(n)_X^c[1]\\
@VVV @VVV @VVV @VVV\\
j_!j^*Rf_*\bZ_l(n)_{\tilde X}^c@>>> Rf_*\bZ_l(n)_{\tilde X}^c@>>>
i_*i^*Rf_*\bZ_l(n)_{\tilde X}^c@>>>j_!j^*Rf_*\bZ_l(n)_{\tilde X}^c[1]
\end{CD}\]

By theorem \ref{t2.-1} and the assumption, the left vertical morphism is an
isomorphism and 
\[i_*i^*Rf_*\bZ_l(n)_{\tilde X}^c \iso i_*Rf'_*\bZ_l(n)_{\tilde Z},\]
where $f'$ is the restriction of $f$ to $\tilde Z$. Moreover,
$i_*i^*\bZ_l(n)_X^c\simeq i_*\bZ_l(n)_Z^c$. We deduce from this the exact triangle
of corollary \ref{c2.1}.
\end{pf}

\begin{thm}\label{t2.0} {\rm (Homotopy invariance)} Let $f:X\to S$ be a morphism
locally of finite type whose geometric fibres are isomorphic to affine spaces
(e.g. a torsor under a vector bundle). Then, the natural map
\[\bZ_l(n)^c_S\to Rf_*f^*\bZ_l(n)^c_S\iso Rf_*\bZ_l(n)^c_X\]
is an isomorphism.
\end{thm}

\begin{pf} This is clear from the homotopy invariance of \'etale cohomology with
coefficients
$\Z/l^\nu(n)$ \cite[Expos\'e XV, cor. 2.2]{sga4}.\end{pf}

\begin{thm}{\rm (Purity, \cf \cite[(3.17)]{jannsen})}\label{t2.1} Let $(X,Z)$ be a
regular
$S$-pair of codimension
$c$ and $i:Z\inj X$ the corresponding closed immersion.\\
a) Suppose $(X,Z)$ is a smooth $S$-pair. Then, for all
$n\in\bZ$ is an isomorphism
\[\bZ_l(n-c)^c_{|Z}[-2c]\iso Ri^!\bZ_l(n)^c_{|X} \]
with the usual naturality properties.\\
b) Assume $S$ and $X$ are quasi-compact, quasi-separated over $\Spec \bZ[1/l]$ and
that $X$ has \'etale $l$-cohomological dimension $N$. Then there are maps 
\[Ri^!\bZ_l(n)^c_{|X}@<<< M(N)\bZ_l(n-c)^c_{|Z}[-2c]@>>>
\bZ_l(n-c)^c_{|Z}[-2c]\]
whose cones are killed by $M(N)$. Here $M(N)$ is the integer described in
\cite[3.3]{th}. In particular,
\[Ri^!\bQ_l(n)^c_{|X}\simeq \bQ_l(n-c)^c_{|Z}[-2c]\]
\end{thm}

\begin{pf} a) follows from classical purity (\cite[Expos\'e XVI, cor. 3.8]{sga4} and
\cite[ch. VI, \S 6]{milneet}), since $i^!$ has as a left adjoint the exact functor
$i_*$, \cf lemma \ref{l1.3} d). b) follows from Thomason's theorem on absolute
cohomological purity
\cite[3.5]{th} and the limit result already quoted \cite[C.9]{TT}.\end{pf}

\begin{rk} O. Gabber has announced a proof of absolute cohomological purity in
general. Using this, one can get rid of the integer $M(N)$ in theorem \ref{t2.1}
b).
\end{rk}

\begin{thm} \label{t2.4} {\rm (Cohomological dimension)} If $\cd_l(X)\le
N$,
then
$H^q_\cont(X,\bZ_l(n))=0$ for $q>N$.
\end{thm}

\begin{pf} This follows from proposition \ref{p1.2}.\end{pf}

\begin{thm} \label{t2.6} {\rm (Trace maps)} Let $f:X'\to X$ be a finite
flat morphism of constant separable degree $d$.  Then there exists a morphism
\[Tr_f:f_*\bZ_l(n)^c_{X'}\to \bZ_l(n)^c_X\]
whose composition with the adjunction morphism
\[\bZ_l(n)^c_X\to f_*\bZ_l(n)^c_{X'}\]
is multiplication by $d$. This morphism commutes with base change. If
$g:X''\to X'$
is another such morphism, then
\[Tr_{f\circ g}=Tr_f\circ f_*Tr_g.\]
Moreover, if $Z@>i>> X$ is a smooth $S$-pair of codimension $c$, $f:X'\to X$ is
finite and flat and
$Z'=f^{-1}(Z)$, then the diagram
\[\begin{CD}
\bZ_l(n-c)^c_{|Z'}[-2c]\iso Ri^{'!}\bZ_l(n)^c_{|X'} \\
@V{Tr_{f'}}VV @V{Tr_f}VV\\
\bZ_l(n-c)^c_{|Z}[-2c]\iso Ri^!\bZ_l(n)^c_{|X} 
\end{CD}\]
commutes, where $i'$ is the closed immersion $Z'\inj X'$ and $f'$ is the
restriction of $f$ to $Z'$.
\end{thm}

\begin{pf} Except for the last claim, this has nothing to do with continuous
cohomology: the way the trace is described in
\cite[exp. IX, \S 5.1]{sga4} shows that it defines a natural 
transformation
\[Tr_f:f_*f^*\to Id\]
on $\sD^b(X_\Et)$, which commutes with base change. In fact, \cite[exp. IX, \S
5.1]{sga4} only deals with the case when $f$ is
\'etale, but the extension to a finite morphism reduces by the same method to the
case of a radicial morphism, when $f^*$ is an equivalence of categories.

The last claim is local for the \'etale topology and to prove it we reduce as
usual to the case when $c=1$ and $Z$ is defined say, by a
rational function $a\in \Gamma(X-Z,\bG_m)$. Then the isomorphism
\[i_*\bZ/l^\nu(n-1)_{Z}\iso \sH^2_Z(\bZ/l^\nu(n))\]
is given for all $\nu$ by the ``extraordinary cup-product" \cite[9.14]{iversen} by
the Kummer class of $a$ in $H^2_Z(X,\bZ/l^\nu(1))$, and similarly for $(X',Z')$. 
The diagram of sheaves
\[\begin{CD}
f_*i'_*\bZ/l^\nu(n-1)_{Z'}@>(a)>> f_*\sH^2_{Z'}(\bZ/l^\nu(n))\\
@V{f_*i'_*Tr_{f'}}VV @V{Tr_f}VV\\
i_*\bZ/l^\nu(n-1)_{Z}@>(a)>> \sH^2_Z(\bZ/l^\nu(n))
\end{CD}\]
clearly commutes for all $\nu\ge 1$; applying $Ri^!R\lim$ to this
commutative diagram of $l$-adic sheaves yields the result.
\end{pf}

\begin{thm} \label{t2.7} {\rm (Gersten's conjecture)} Let $k$ be a field of
characteristic
$\neq l$ and $X$ a smooth affine variety over $k$. Let $\{x_1,\dots, x_r\}$ be a
finite set of points of $X$ and $Y=\Spec \sO_{X,x_1,\dots,x_r}$. Then the Gersten
complex
\[0\to \tilde H^i_\cont(Y/k,\bZ_l(n))\to \coprod_{y\in Y^{(0)}}
\tilde H^i_\cont(k(y)/k,\bZ_l(n))\to \coprod_{y\in Y^{(1)}}
\tilde H^{i-1}_\cont(k(y)/k,\bZ_l(n-1))\to\dots\]
is universally exact in the sense of Grayson \cite{grayson} for all $i$. If
$X'@>f>> X$ is finite and flat, the diagram
\[\begin{CD}
0&\to& \tilde H^i_\cont(Y'/k,\bZ_l(n))&\to& \displaystyle\coprod_{y\in Y^{'(0)}}
\tilde H^i_\cont(k(y)/k,\bZ_l(n))&\to& \displaystyle\coprod_{y\in Y^{'(1)}}
\tilde H^{i-1}_\cont(k(y)/k,\bZ_l(n-1))&\to&\dots\\
&& @V{f_*}VV @V{f_*}VV @V{f_*}VV \\
0&\to& \tilde H^i_\cont(Y/k,\bZ_l(n))&\to& \displaystyle\coprod_{y\in Y^{(0)}}
\tilde H^i_\cont(k(y)/k,\bZ_l(n))&\to& \displaystyle\coprod_{y\in Y^{(1)}}
\tilde H^{i-1}_\cont(k(y)/k,\bZ_l(n-1))&\to&\dots
\end{CD}\]
commutes, where $Y'=Y\times_X X'$.
\end{thm}

\begin{pf} The first claim follows from \cite[cor. 5.1.11, prop. 5.3.2 and ex.
7.1.8]{cthk}, theorem
\ref{t2.0} and theorem \ref{t2.1}, plus \cite[th. 6.2.5]{cthk} and theorem
\ref{t2.6} in case $k$ is finite. More precisely, by homotopy invariance we get
universally exact Cousin complexes
\[0\to \tilde H^i_\cont(Y/k,\bZ_l(n))\to \coprod_{y\in Y^{(0)}}
\tilde H^i_{y,\cont}(Y/k,\bZ_l(n))\to \coprod_{y\in Y^{(1)}}
\tilde H^{i+1}_{y,\cont}(Y/k,\bZ_l(n-1))\to\dots\]
and then purity allows us to rewrite them as Gersten complexes. The second claim
is clear for the Cousin complexes, and follows for Gersten complexes from the last
part of  theorem
\ref{t2.6}.
\end{pf}

\begin{rk} With more effort one could check that
\[K\mapsto \tilde H^*_\cont(K/k,\bZ_l(n+*)),\]
where $k$ is a base field and $K$ runs through finitely generated extensions of
$K$, defines a cycle module in the sense of Rost \cite{rost}. We skip this since
it will not be needed here.
\end{rk}

\subsection{Cohomology with proper supports} For $X@>f>> S$ separated
of finite type, we define {\it continuous cohomology with proper supports}
as
\begin{equation}\label{eq2.4} H^i_{\cont,c}(X/S,\bZ_l(n)):= H^i_\et(S,R\lim
Rf_!\bZ/l^\nu(n))\end{equation}
where $Rf_!$ is the functor ``higher direct image with
proper support" from \cite[appendix to expos\'e XVII]{sga4}. By abuse of notation,
we shall denote the object $R\lim
Rf_!\bZ/l^\nu(n)$ by $Rf_!\bZ_l(n)^c_X$. We have

\begin{lemma}\label{l2.3} There are natural homomorphisms
\[H^i_{\cont,c}(X/S,\bZ_l(n))\to H^i_\cont(X,\bZ_l(n))\]
which are isomorphisms when $X/S$ is proper.
\end{lemma}

\begin{pf} This follows from the same fact with finite coefficients.\end{pf}

Using \cite[exp. XVII, (6.2.7.3)]{sga4}, one can define trace maps for
finite flat morphisms in continuous cohomology with proper supports,
verifying the same properties as in cohomology without supports (see
theorem \ref{t2.6}). We shall not need them here, so skip the details. 

\subsection{(Borel-Moore) \'etale homology} For any $X@>f>> S$ locally of finite
type, define 
\[L_l(n)^c_{X/S}:= R\lim Rf^!\bZ/l^\nu(n)\]
where $Rf^!$ is
the extraordinary inverse image of \cite[Expos\'e XVIII]{sga4}. We define
{\it continuous \'etale homology} as
\begin{equation}\label{eq2.5}
H_i^\cont(X/S,\bZ_l(n)):=H^{-i}_\et(X,L_l(-n)^c_{X/S}).\end{equation}

This coincides with Laumon's definition of $l$-adic \'etale homology \cite{laumon}
in the case $S=\Spec k$, $k$ an algebraically closed field.

\begin{prop}\label{p2.2} a) Let $i:Z\inj X$ be a closed immersion and
$j:U\inj X$ be the complementary open immersion. Then there are long exact
sequences
\[\dots\to H_i^\cont(Z/S,\bZ_l(n))@>i_*>> H_i^\cont(X/S,\bZ_l(n))@>j^!>>
H_i^\cont(U/S,\bZ_l(n))\to H_{i-1}^\cont(Z/S,\bZ_l(n))\to\dots\]
b) Continuous \'etale homology is contravariant for \'etale
$S$-morphisms or finite flat
$S$-morphisms and covariant for finite flat $S$-morphisms; if $u:X'\to X$ is
finite and flat, of constant separable degree $d$, then the composition
\[H_*^\cont(X/S,\bZ_l(n)) @>u^*>> H_*^\cont(X'/S,\bZ_l(n)) @>u_*>>
H_*^\cont(X/S,\bZ_l(n))\]
is multiplication by $d$.\\
c) If $X$ is smooth compactifiable of pure dimension $d$, there are canonical
isomorphisms
\[H_i^\cont(X/S,\bZ_l(n))\simeq H^{2d-i}_\cont(X,\bZ_l(d-n)).\]
\end{prop}

\begin{pf} a) This follows from applying $R\lim$ to the classical exact
triangles
\[i_*Ri^!Rf^!\bZ/l^\nu(n)\to Rf^!\bZ/l^\nu(n)\to
Rj_*j^!Rf^!\bZ/l^\nu(n)\to i_*Ri^!Rf^!\bZ/l^\nu(n)[1].\]

b) If $u:X'\to X$ is \'etale or finite and flat, the map
\[u^*:H_*^\cont(X/S,\bZ_l(n))\to H_*^\cont(X'/S,\bZ_l(n))\]
is defined through the unit morphism
\[L_l(-n)^c_{X/S}\to
Ru_*u^*L_l(-n)^c_{X/S}=Ru_*u^!L_l(-n)^c_{X/S}=Ru_*L_l(-n)^c_{X'/S}.\]

If $u$ is finite and flat, the map
\[u_*:H_*^\cont(X'/S,\bZ_l(n))\to H_*^\cont(X/S,\bZ_l(n))\]
is defined through the trace map (see proof of theorem \ref{t2.6})
\[u_*L_l(-n)^c_{X'/S}=u_*u^*L_l(-n)^c_{X/S}@>Tr_u>> L_l(-n)^c_{X/S}.\]

The formula $u_*u^*=d$ is then clear.

c) This follows from applying $R\lim$ to the isomorphisms of
\cite[expos\'e XVIII, th. 3.2.5]{sga4}
\[Rf^!\bZ/l^\nu(n)\simeq f^*\bZ/l^\nu(n+d)[2d].\]
\end{pf}

\subsection{Duality} Recall \cite[expos\'e XVIII, th. 3.1.4]{sga4} that, for
$f:X\to S$ separated of finite type, the functors
$Rf_!$ and $Rf^!$ are (partially) adjoint on the derived categories of sheaves
of $\bZ/l^\nu$-Modules. Let $T_f:Rf_!Rf^!\to Id$ be the counit of this
adjunction. It induces an isomorphism
\[Rf_*RHom_{\sD^b(X_\et,\bZ/l^\nu)}(A,Rf^!B)@>(T_f)_*>>
RHom_{\sD^b(S_\et,\bZ/l^\nu)}(Rf_!A,B)\]
for any $A,B$. Taking $A=\bZ/l^\nu, B=\bZ/l^\nu(-n)$, this translates as
\[Rf_*Rf^!\bZ/l^\nu(-n)\iso RHom_S(Rf_!\bZ/l^\nu,\bZ/l^\nu(-n))\simeq
RHom_S(Rf_!\bZ/l^\nu(n),\bZ/l^\nu).\]

These isomorphisms are compatible when $\nu$ varies in the sense that,
for
$\mu\ge \nu$, the diagram of pairings
\[\begin{CD}
Rf_*Rf^!\bZ/l^\mu(-n)\oo^L Rf_!\bZ/l^\mu(n)@>>> \bZ/l^\mu\\
@VVV @VVV\\
Rf_*Rf^!\bZ/l^\nu(-n)\oo^L Rf_!\bZ/l^\nu(n)@>>> \bZ/l^\nu
\end{CD}\]
commutes, where the left map is the tensor product of the two projections.
This translates as a pairing of $\bZ_l$-sheaves, in $\sD_{\bZ_l}(S_\et)$
\[Rf^\ladic_*L_l(-n)_{X/S}\oo^L Rf^\ladic_!\bZ_l(n)_X@>>> \bZ_l\]
whose adjunction
\begin{equation}\label{eq2.6}
Rf^\ladic_*L_l(-n)_{X/S}@>>>
R\underline{Hom}_{\sD_{\bZ_l}(S_\et)}(Rf^\ladic_!\bZ_l(n)_X,\bZ_l)
\end{equation}
is an isomorphism in $\sD_{\bZ_l}(S_\et)$.

Note that, more generally, there are pairings of $l$-adic sheaves
\[Rf^\ladic_*L_l(-p)_{X/S}\oo^L Rf^\ladic_!\bZ_l(q)_X@>>> \bZ_l(q-p)\]
which induce products in cohomology
\begin{equation}\label{eq3.2} H_i^\cont(X/S,\bZ_l(p))\times
H^j_{\cont,c}(X/S,\bZ_l(q))\to H^{j-i}_\cont(S,\bZ_l(q-p)).
\end{equation}

\subsection{Schemes over $\bF_p$} In this paper, we shall mainly concentrate on the
case
$S=\Spec\bF_p$, where $p$ is a prime number $\neq l$. The next results are special
to this case.

\begin{lemma} \label{l2.2} Let $X$ be a scheme of finite type over $\Spec \bF_p$.
Then for all
$n\in\bZ$,
$i\ge 0$,\\
a) $H^i_\cont(X,\bZ_l(n))\iso\lim H^i_\et(X,\bZ/l^\nu(n))$.\\
b) $H^i_\cont(X,\bZ_l(n))$ is a finitely generated $\bZ_l$-module.\\
The same assertions hold for \'etale cohomology with proper supports and \'etale
homology. 
\end{lemma}

\begin{pf} a) follows from Deligne's theorem that the $H^i_\et(X,\bZ/l^\nu(n))$ are
finite \cite[Th. finitude]{sga}. b) By a), this $\bZ_l$-module is compact. 
On the other hand, the exact sequence
\[H^i_\cont(X,\bZ_l(n))@>l>> H^i_\cont(X,\bZ_l(n))\to
H^i_\et(X,\bZ/l(n))\]
shows that it is finite modulo $l$. Therefore it is finitely generated.
Similarly for the other cases.\end{pf}

\begin{prop}\label{p2.6} Let $X$ be a (not necessarily smooth) variety over
a finite field $k$, with $\car k\neq l$. Let $\bar k$ be the algebraic
closure of $k$, $G_k=Gal(\bar k/k)$ and $\overline X=\bar k\otimes_k X$. Then, for
all
$n\in\Z$ and $q\ge 0$, there is an exact sequence
\[0\to H^{q-1}_\cont(\overline X,\bZ_l(n))_{G_k}\to H^q_\cont(X,\bZ_l(n))\to
H^q_\cont(\overline X,\bZ_l(n))^{G_k}\to 0.\]
There are similar exact sequences for continuous cohomology with proper supports
and continuous homology.
\end{prop}

\begin{pf} (Note how we use essentially the $l$-adic sheaf definition of continuous
\'etale cohomology, and not the complexes $\bZ_l(n)^c$.) We do it only for
continuous cohomology, the other cases being similar. By
\cite[(3.4)]{jannsen}, there is a ``Hochschild-Serre" spectral sequence
\[E_2^{p,q}=H^p_\cont(k,H^q_\cont(\overline X,\bZ_l(n)))\Rightarrow
H^{p+q}_\cont(X,\bZ_l(n)).\]

Moreover, by \cite[\S 2]{jannsen}, we can identify the $E_2$-terms to continuous
cohomology of the profinite group $G_k=Gal(\bar k/k)$ in the sense of Tate
\cite{tate}.

Since $k$ is finite, we have
\[E_2^{p,q}=\begin{cases} 0&\text{if $p\neq 0,1$}\\
H^q_\cont(\overline X,\bZ_l(n))^{G_k}&\text{if $p=0$}\\
H^q_\cont(\overline X,\bZ_l(n))_{G_k}&\text{if $p=1$.}
\end{cases}\]

This concludes the proof.\end{pf}

\begin{cor}\label{l6.1} With $X$ as in proposition \ref{p2.6}, let
$d_i(n)=\dim_{\bQ_l} H^i_\cont(X,\bQ_l(n))$. Then we have
\[\dim H^i_\cont(\overline X,\bQ_l(n))^{G_k}=\dim H^i_\cont(\overline
X,\bQ_l(n))_{G_k}=d_i(n)-d_{i-1}(n)+\dots.\] There are similar identities for
continuous cohomology with proper supports and continuous homology.
\end{cor}

\begin{pf} Since $H^i_\cont(\overline
X,\bQ_l(n))$ is a finite-dimensional $\bQ_l$-vector space, we have
\[\dim_{\bQ_l}H^i_\cont(\overline X,\bQ_l(n))^{G_k}=\dim_{\bQ_l}H^i_\cont(\overline
X,\bQ_l(n))_{G_k}\]
for all $i$, as follows from the exact sequence
\[0\to H^i_\cont(\overline X,\bQ_l(n))^{G_k}\to H^i_\cont(\overline
X,\bQ_l(n))@>1-F>> H^i_\cont(\overline X,\bQ_l(n))\to H^i_\cont(\overline
X,\bQ_l(n))_{G_k}\to 0.\]

The result then follows from proposition \ref{p2.6} which shows that
\[d_i(n)=\dim H^i_\cont(\overline X,\bQ_l(n))^{G_k}+\dim H^{i-1}_\cont(\overline
X,\bQ_l(n))_{G_k}\]\end{pf}

\begin{cor}\label{p2.1} Let $X$ be a (not necessarily smooth) variety over
a finite field $k$, with $\car k\neq l$. Then, for all $n\in\Z$,
\[\sum_{i\ge 0}(-1)^i\dim_{\bQ_l} H^i_\cont(X,\bQ_l(n))=0.\]
The same holds for continuous cohomology with proper supports and for continuous
\'etale homology.\qed
\end{cor}

\subsection{Duality again} Let $B\mapsto
B^\#[1]$ be the functor from
$\sD^b(\bZ/l^\nu-mod)$ to $\sD^b((\Spec \bF_p)_\et,\bZ/l^\nu)$ which
associates to a complex of $\bZ/l^\nu$-modules the corresponding complex of constant
$(\Spec \bF_p)_\et$-sheaves, shifted once to the left. This is right adjoint to
$R\Gamma$, with counit 
\[R\Gamma(\bF_p,B^\#[1])\to B\]
induced by the canonical isomorphism
\[H^1(\bF_p,\bZ/l^\nu) \iso \bZ/l^\nu\]
given by evaluation on the Frobenius. Applying $R\Gamma(\bF_p,-)$ to \eqref{eq2.6}
and then applying this adjunction, we get a chain of adjunction isomorphisms in
$\sD^b_{\bZ_l}(*)$:
\begin{multline*}
R\Gamma(X,L_l(-n))\simeq R\Gamma(\bF_p,Rf^\ladic_*L_l(-n)_{X/\Spec \bF_p})\\
\simeq
R\Gamma(\bF_p,R\underline{Hom}_{\sD_{\bZ_l}((\Spec
\bF_p)_\et)}(Rf^\ladic_!\bZ_l(n)_X,\bZ_l))\\
\simeq  R{Hom}_{\sD_{\bZ_l}((\Spec
\bF_p)_\et)}(Rf^\ladic_!\bZ_l(n)_X,\bZ_l))\\
\simeq RHom_{\sD^b_{\bZ_l}(*)}(R\Gamma(\bF_p,Rf^\ladic_!\bZ_l(n)_X),\bZ_l)[1].
\end{multline*}

By \cite[th. 7.2]{ek}, the two sides can be identified with objects of the derived
category of finitely generated $\bZ_l$-modules. We then get isomorphisms
\[H_{i}^{\cont}(X,\bZ_l(n))\iso
H^{i+1}(RHom_{\sD^b(\bZ_l-mod)}(R\Gamma(\bF_p,Rf^\ladic_!\bZ_l(n)_X),\bZ_l)).\]

By the standard Ext spectral sequence, this yields in a short exact sequence
\begin{equation}\label{eq2.7} 0\to Ext(H^{i+2}_{c,\cont}(X,\bZ_l(n)),\Z_l) \to
H_{i}^{\cont}(X,\bZ_l(n))\to Hom(H^{i+1}_{c,\cont}(X,\bZ_l(n)),\Z_l)
\to 0.\end{equation}

For any
$\bZ_l$-module
$A$, let
$\bar A$ denote $A/\tors$. 

\begin{lemma}\label{l2.7} Let $A$ be a finitely generated $\bZ_l$-module. Then
there is a canonical isomorphism
\[Ext(A,\bZ_l)\simeq Hom(A_\tors,\bQ_l/\bZ_l).\]
\end{lemma}

\begin{pf} This follows from the commutative diagram with exact rows and columns
\[\begin{CD}
0@>>>0@>>>0@>>>0\\
@AAA @AAA @AAA @AAA\\
0@>>>Hom(A_\tors,\bQ_l/\bZ_l) @>>> Ext(A_\tors,\bZ_l) @>>> 0\\
@AAA @AAA @AAA @AAA\\
Hom(A,\bQ_l) @>>>Hom(A,\bQ_l/\bZ_l) @>>> Ext(A,\bZ_l) @>>> Ext(A,\bQ_l)\\
@AAA @AAA @AAA @AAA\\
Hom(\bar A,\bQ_l) @>>>Hom(\bar A,\bQ_l/\bZ_l) @>>> 0 @>>> 0.
\end{CD}\]
(We used that $A_\tors$ has finite exponent and that $\bar
A$ is free). \end{pf}

In view of lemma \ref{l2.7}, the first group in the exact sequence \eqref{eq2.7} is
isomorphic to \[Hom(H^{i+2}_{c,\cont}(X,\bZ_l(n))_\tors,\Q_l/\Z_l)\] while the
last one is obviously isomorphic to  \[Hom(\bar
H^{i+1}_{c,\cont}(X,\bZ_l(n)),\Z_l).\]
We finally get the following theorem:

\begin{thm} {\rm (\cf \cite[lemma 5.3]{milne2} in the smooth, proper
case)}\label{t2.3} For any $X$ separated of finite type over $\bF_p$, the geometric
adjunction between
$Rf_!$ and $Rf^!$ and the arithmetic adjunction over $\Spec \bF_p$ induce perfect
pairings
\[\bar H_{i}^{\cont}(X,\bZ_l(n))\times \bar H^{i+1}_{c,\cont}(X,\bZ_l(n))\to
\bZ_l\]
given by the products \eqref{eq3.2} followed by the trace,
as well as isomorphisms
\[Hom(H^{i+2}_{c,\cont}(X,\bZ_l(n))_\tors,\Q_l/\Z_l)\iso
H_i^\cont(X,\bZ_l(n))_\tors.\]
In particular, if $X$ is smooth and proper of pure
dimension $d$,
$H^{2d+1}_\cont(X,\bZ_l(d))$ is canonically isomorphic to $\bZ_l$ and the
resulting pairings induced by cup-product
\[H^i_\cont(X,\bZ_l(n))\times H^{2d+1-i}_\cont(X,\bZ_l(d-n))\to \bZ_l\]
are perfect modulo torsion.\qed
\end{thm}

\section{Restriction to $\bF_p$}\label{3}

We now start our investigation. Let $p\neq l$. We shall begin by
computing the restriction of
$\bZ_l(n)^c_{\F_p}$ and
$\bQ_l(n)^c_{\F_p}$ to the {\it small}
\'etale site of $\Spec {\F_p}$.

Let $\hat\bZ=\lim\bZ/n$. Note that the discrete abelian group $\bQ/\bZ$ is
naturally a module over this profinite ring, with continuous action. Consider
$\hat\bZ$ as an additive profinite group. We introduce the discrete topological
$\hat\bZ$-module 
\[\tilde M=\bQ/\bZ\times \bQ/\bZ\]
provided with the following action: for $(a,r,s)\in\hat\bZ\times \bQ/\bZ\times
\bQ/\bZ$, $a(r,s)=(r,ar+s)$. So we have a short exact sequence
\begin{equation}\label{eq17}0\to \bQ/\bZ\to\tilde M\to \bQ/\bZ\to 0\end{equation}
whose kernel and cokernel are trivial $\hat\bZ$-modules.

Let $M$ be the pull-back of this extension under the projection map
$\bQ\to\bQ/\bZ$, so that we have another short exact sequence of $\hat\bZ$-modules
\begin{equation}\label{eq3.1}0\to\bQ/\bZ\to M\to \bQ\to 0.\end{equation}

\begin{defn}\label{d2.1} We denote by $\bZ^c$ the object of
$\sD^b(Ab(\hat\bZ))$ represented by the complex of length $1$
\begin{equation}\label{eq0}\bQ@>\gamma>> M\end{equation}
where $M$ is as above and $\gamma$ is the composition $\bQ\surj\bQ/\bZ\inj M$.
($\bQ$ sits in degree $0$ and $M$ in degree $1$.)\\
We set $\bQ^c=\bZ^c\otimes \bQ$.
\end{defn}

\begin{prop}\label{p2.4}We have
\[\sH^i(\bZ^c)=
\begin{cases}
\bZ &\text{if $i=0$}\\
\bQ &\text{if $i=1$}\\
0 &\text{otherwise,}
\end{cases}
\]
\[\bH^i(\hat\bZ,\bZ^c)=
\begin{cases}
\bZ &\text{if $i=0$}\\
\bZ &\text{if $i=1$}\\
0 &\text{otherwise.}
\end{cases}
\]
\end{prop}

\begin{pf} The first formulas are obvious in view of the definition of
$\bZ^c$. For the next ones, we need: 

\begin{lemma}\label{l2.4} $H^i(\hat\bZ,M)=0$ for $i>0$; there is an exact
sequence 
\[0\to \bQ/\bZ\to H^0(\hat\bZ,M)\to \bZ\to 0.\]
\end{lemma}

Indeed, $H^2(\hat\bZ,M)=0$ since $M$ is divisible, and the exact sequence 
\eqref{eq3.1} yields a long cohomology exact sequence
\[
0\to H^0(\hat\bZ,\bQ/\bZ)\to H^0(\hat\bZ,M)\to H^0(\hat\bZ,\bQ)@>\delta>>
H^1(\hat\bZ,\bQ/\bZ)\\
\to H^1(\hat\bZ,M)\to H^1(\hat\bZ,\bQ).
\]

The definition of the action of $\hat{\bZ}$ on $\tilde M$ shows that $\delta$ is
surjective, which gives the claim.

To compute the hypercohomology of $\bZ^c$, we use lemma \ref{l2.4} and the long
exact sequence
\begin{multline*}
0\to \bH^0(\hat\bZ,\bZ^c)\to H^0(\hat\bZ,\bQ)@>\gamma_*>> H^0(\hat\bZ,M)
\to\bH^1(\hat\bZ,\bZ^c)\\
\to H^1(\hat\bZ,\bQ)@>\gamma_*>> H^1(\hat\bZ,M)
\to
\bH^2(\hat\bZ,\bZ^c)\to 0.\end{multline*}

We have $H^i(\hat\bZ,\bQ)=H^i(\hat\bZ,M)=0$ for $i>0$ and the map
\[H^0(\hat\bZ,\bQ)@>\gamma_*>> H^0(\hat\bZ,M)\]
sends $\bQ$ onto the kernel $\bQ/\bZ$ of the extension in lemma \ref{l2.4}. Hence
the claims on $\bH^*(\hat\bZ,\bZ^c)$.\end{pf}

The canonical generator $e$ of $H^1_\cont(\hat\bZ,\hat\bZ)$ represented by the
continuous homomorphism
\[e:\hat\bZ@>Id>> \hat\bZ\]
will play an important r\^ole in the sequel. We now
use it to describe the boundary homomorphism
associated to $\bZ^c$:

\begin{prop}\label{p2.3} Let $\partial$ be the morphism defined by the exact
triangle
\[\bZ\to \bZ^c\to
\bQ[-1]@>\partial>>\bZ[1]\]
stemming from proposition \ref{p2.4}. Then the diagram
\[\begin{CD}
\bQ[-1]@>\partial>>\bZ[1]\\
@VVV @A{\beta}AA\\
\bQ/\bZ[-1] @>\cdot e>>\bQ/\bZ[0]
\end{CD}\]
commutes, where $\cdot e$ is cup-product by $e$, $\beta$ is the Bockstein and
the left vertical map is induced by the projection $\bQ\to
\bQ/\bZ$.
\end{prop}

\begin{pf} This follows from the two commutative diagrams of exact triangles
\[\begin{CD}
&&\bQ[-1]@= \bQ[-1]\\
&&@VVV @V{\partial}VV\\
\bQ[0]@>>>\bQ/\bZ[0]@>\beta>> \bZ[1]@>>>\bQ[1]\\
||&& @VVV @VVV ||\\
\bQ[0]@>\gamma>> M[0]@>>> \bZ^c[1]@>>>\bQ[1]\\
&&@VVV @VVV\\
&&\bQ[0]@= \bQ[0]
\end{CD}\qquad\qquad
\begin{CD}
\bQ[-1]@>>> \bQ/\bZ[-1]\\
@VVV @VVV\\
\bQ/\bZ[0]@=\bQ/\bZ[0]\\
@VVV @VVV\\
M[0]@>>> \tilde M[0]\\
@VVV @VVV\\
\bQ[0]@>>> \bQ/\bZ[0]
\end{CD}\]
and the fact that, by definition of $\tilde M$, the boundary map $\bQ/\bZ[-1]\to
\bQ/\bZ[0]$ corresponding to
$\tilde M$ is given by cup-product by $e$.\end{pf}

\begin{cor}\label{c2.3} $\bQ^c\simeq \bQ[0]\oplus\bQ[-1]$.
\end{cor}

\begin{pf} Indeed, proposition \ref{p2.3} shows that $\partial\otimes
\bQ=0$.\end{pf}

We can consider $\bZ^c$ as a complex of sheaves over $\Spec {\F_p}$, by
identifying $G$ with $\hat\bZ$ by means of the absolute Frobenius.

\begin{thm}\label{t2.2} Denote by
$\bZ_l(n)^c_{|{\F_p}}$ (\resp
$\bQ_l(n)^c_{|{\F_p}}$) the object
$R\pi_*\bZ_l(n)^c_{\F_p}$ (\resp $R\pi_*\bQ_l(n)^c_{\F_p}$) of the derived
category
$\sD^+(Ab(\Spec {\F_p})_\et)$, where $\pi:(\Spec {\F_p})_\Et\to (\Spec
{\F_p})_\et$ is the natural projection. Then\\
a) $\bZ_l(n)^c_{|{\F_p}}=\bQ_l/\bZ_l(n)[-1]$ if $n\neq 0$.\\
b) $\bZ_l(0)^c_{|{\F_p}}=\bZ^c\otimes\bZ_l$.
\end{thm}

\begin{pf} We first compute the
cohomology sheaves of $\bZ_l(n)^c_{|\bF_p}$ for $n\neq 0$. Let $k$ be a finite field
of characteristic $p$. We note that the groups
$H^0(k,\bZ/l^\nu(n))$ and $H^1(k,\bZ/l^\nu(n))$ have order bounded independently of
$\nu$. Using lemma \ref{l2.2}, it follows that
\begin{align*}
&H^0_\cont(k,\bZ_l(n))=0\\
&H^1_\cont(k,\bZ_l(n))\quad\text{is finite}\\
&H^q_\cont(k,\bZ_l(n))=0\quad\text{for $q>1$.}
\end{align*}

It follows from this that $H^q_\cont(k,\bQ_l(n))=0$ for all $q$, hence that
$\bQ_l(n)^c_{|\bF_p}=0$. We conclude by using lemma
\ref{l2.1}.

To prove b), we note that, after proposition \ref{p2.4},
$\bZ^c\oo^L\bZ/l^\nu\allowbreak\simeq \bZ/l^\nu[0]$ for all $\nu\ge 1$. From this
we deduce a morphism
$\bZ^c\otimes\bZ_l\to \bZ_l(0)^c_{|{\F_p}}$ out of the composite morphisms
\[\bZ^c\otimes\bZ_l\to \bZ^c\oo^L\bZ/l^\nu\to \bZ/l^\nu[0].\]

We claim that this morphism is an isomorphism. To see
this, we note that it inserts into an exact triangle
\[C@>>> \bZ^c\otimes\bZ_l@>>> \bZ_l(0)^c_{|{\F_p}}@>>> C[1]\]
where $C=R\lim(\bZ^c\otimes\bZ_l,l)$. The claim now follows from

\begin{lemma} $C=0$.
\end{lemma}

Indeed, we shall show that the map $R\lim(\bQ_l,l)@>R\lim\gamma\otimes 1>>
R\lim(M\otimes\bZ_l,l)$ is an isomorphism. We have evidently
$R\lim(\bQ_l,l)=\bQ_l[0]$. On the other hand, lemma
\ref{l2.4} and a little bit of computation show that
\begin{align*}
\bH^0(k,R\lim(M\otimes\bZ_l,l))&=\bQ_l\\
\bH^q(k,R\lim(M\otimes\bZ_l,l))&=0\quad\text{for $q>0$}\\
\intertext{ for any finite extension $k/\bF_p$, that is}
R\lim(M\otimes\bZ_l,l)&\simeq \bQ_l[0].
\end{align*}

To conclude, it suffices to notice that the map $\bQ_l[0]\to \bQ_l[0]$ defined by
$\gamma\otimes 1$ and the above computation is the inverse limit
of the projections $\bQ_l\to \bQ_l/\bZ_l$ (for multiplication by $l$), hence the
identity.\end{pf}

\begin{cor}\label{c2.2} The canonical generator $e$ of $H^1_\cont(\F_p,\Z_l)$ is
the image of the canonical generator of $\bH^1(\F_p,\Z^c)$ (also denoted by
$e$). The composition
\[\bQ^c\to \bQ[-1]\to \bQ^c[-1]\]
is given by cup-product by $e$, where the left map is ``projection onto $H^1$" and
the right one is ``inclusion of $H^0$".\qed
\end{cor}

\section{Cohomology and homology theories}\label{4}

Let $k$ be a field. We shall only need a basic concept of a ``pure
cohomology theory":
 
\begin{defn}\label{d3.1} Let $Sm/k$ denote the category of smooth $k$-schemes of
finite type, and let
$\sA$ be an abelian category. A {\em pure cohomology theory} with values in $\sA$
is an assignment
\[h=(h^n):Ob(Sm/k)\to
Ob(\sA^\bZ)\] 
such that, for any closed immersion $Z\inj X$ of pure codimension $c$ of
smooth schemes, with complementary open set $U$, there is a long exact
sequence
\[\dots\to h^{i-2c}(Z)\to h^i(X)\to h^i(U)\to h^{i-2c+1}(Z)\to\dots\]
b) $h$ {\em has transfers} if, for any finite morphism $f:Y\to X$ in
$Sm/k$, of pure degree $d$, there are maps
\[h^i(X)@>f^*>> h^i(X),\quad h^i(Y)@>f_*>> h^i(X)\quad i\in \bZ\]
such that $f_*\circ f^*$ is multiplication by $d$.
\end{defn}

\begin{thm}\label{t3.1} Let $\sB$ be a thick subcategory of $\sA$. Let $N>0$, and
let
$h$ be a pure cohomology theory such that
$h^*(X)\in \sB$ for any smooth, projective $X$ with $\dim X\le N$.\\
a) Suppose $\car k=0$. Then $h^*(X)\in \sB$ for any $X\in Sm/k$ with $\dim X\le
N$.\\
b) Suppose $\car k>0$. Then the conclusion still holds if $h$ has transfers
and
$\sA/\sB$ is a $\bQ$-linear category.
\end{thm}

\begin{pf} Passing to the quotient category $\sA/\sB$, we reduce to $\sB=0$. We show that
$h^*(X)=0$ for all equidimensional
$X\in
Sm/k$ by induction on
$n=\dim X$. The case $n=0$ follows from the assumption. Suppose $n>0$ and the
claim known in dimensions $<n$. Let $X\in Sm/k$ of dimension $n$. We first prove
that, for any proper closed subset of $X$, $h^*(X)\iso h^*(X-Z)$. Suppose
first $Z$ smooth. Then we get the
result by purity and the induction hypothesis. In general, let
$Z'\subseteq Z$ be
the singular locus of $Z$. If $\car k=0$, then $Z'\neq Z$; if 
$\car k>0$, then $Z'\neq Z$ after a suitable finite radicial
extension of $k$, which is allowable in view of the additional
hypotheses. Then $Z\setminus Z'$ is smooth in $X\setminus Z'$.
Writing $X\setminus Z=(X\setminus Z')\setminus (Z\setminus Z')$, we get
$h^*(X\setminus Z')\iso h^*(X\setminus Z)$. We conclude by Noetherian induction on
$Z$.

If $\car k=0$, then, by resolution of singularities \cite{hiro}, a suitable open
subset of
$X$ embeds into a smooth, projective $k$-variety and the proof is complete. If $\car
k>0$, then by de Jong's theorem \cite{dJ}, a suitable open subset of $X$ is
finite-covered by an open subset of a smooth, projective $k$-variety. We now
conclude by a transfer argument.\end{pf}

We have the following variants of definition \ref{d3.1} and theorem
\ref{t3.1}.

\begin{defn} \label{d3.2} a) Let $\sA$ be an abelian category. A {\em pure graded
cohomology theory} with values in $\sA$ is a pure cohomology theory
\begin{align*}
h&=(h_n)_{n\in\bZ}\\
X&\mapsto (h^*_n(X)):=(h^*(X,n))
\end{align*}
with values
in $\sA^\bN$, satisfying the
following condition: for $(X,Z,U)$ as in definition \ref{d3.1} a), the exact
sequences of {\it loc. cit.} split into exact sequences
\[\dots\to h^{i-2c}(Z,n-c)\to h^i(X,n)\to h^i(U,n)\to
h^{i-2c+1}(Z,n-c)\to\dots\]
c) $h$ {\em has transfers} if all $h_n$ have transfers.
\end{defn}

\begin{thm}\label{t3.2} Let $\sB$ be a thick subcategory of $\sA$. Let $N>0$,
$n\in\bZ$ and let
$h$ be a pure graded cohomology theory such that
$h^i(X,m)\in\sB$ for all $i\le m\le n$ and any smooth, projective $X$ with $\dim
X\le d$.\\
a) Suppose $\car k=0$. Then $h^i(X,m)\in\sB$ for any $i\le m\le n$ and any $X\in
Sm/k$ with
$\dim X\le d$.\\
b) Suppose $\car k>0$. Then the conclusion still holds if $h$ has
transfers and
$\sA/\sB$ is a $\bQ$-linear category.
\end{thm}

\begin{pf} The same as that of theorem \ref{t3.1}, noting that if $i\le m$ and
$c>0$, then $i-2c+1\le m-c$.\end{pf}

\begin{defn}\label{d3.3} Let $Sch/k$ be the category of schemes of finite type
over $k$ and $\sA$ an abelian category.\\
a) A {\em homology theory} with values in
$\sA$ is an assignment
\[h=(h_n)_{n\in\bZ}:Ob(Sch/k)\to Ob(\sA^\bZ)\]
such that, for any closed immersion $Z\inj X$ with complementary open set
$U$, there is a long exact sequence
\[\dots\to h_{i}(Z)\to h_i(X)\to  h_i(U)\to h_{i-1}(Z)\to\dots\]
b) A homology theory $h$ {\em has transfers} if, for any
finite morphism $f:Y\to X$ in
$Sm/k$, of pure degree $d$, there are maps
\[h_i(X)@>f^*>> h_i(Y),\quad h_i(Y)@>f_*>> h_i(X),\quad i\in \bZ\]
such that $f_*\circ f^*$ is multiplication by $d$.
\end{defn}

\begin{thm} \label{t3.3} Let $h$ be a homology theory with values in
$\sA$, and let $\sB$ be a thick subcategory of $\sA$ such that
$h_i(X)\in\sB$ for
all
$i<m$ for all smooth, projective $X$ of dimension $\le N$.\\
a) Suppose $\car k=0$. Then  $h_i(X)\in\sB$ for
all
$i<m$ for all $X\in Sch/k$ of dimension $\le N$.\\
b) Suppose $\car k>0$. Then the conclusion still holds if $h$ has transfers and
$\sA/\sB$ is a
$\bQ$-linear category.
\end{thm}

\begin{pf} By induction on $N$, we have $h_i(X)\iso h_i(U)$ for any open
subscheme $U$
of $X$ and $i< m$. If $\car k=0$, we can choose such an $U$ smooth and
embedded
in a smooth, projective scheme of dimension $N$ (by resolution of singularities
\cite{hiro}). If $\car k>0$, the extra assumptions allow us to pass to a finite,
radicial extension of $k$ if necessary. Then $X$ gets a nonempty smooth open
subscheme, and we can conclude as before, using de Jong's theorem
\cite{dJ}.\end{pf}

\section{Back to $\bZ_l(n)^c$ and $\bQ_l(n)^c$}\label{5}

In this section, we fix a prime number $p\neq l$ and denote by $\sS$ the category
of smooth schemes of finite type over $\bF_p$. We consider the restriction of
$\bZ_l(n)^c_{\bF_p}$ and $\bQ_l(n)^c_{\bF_p}$ to $\sS$, considered as a
subsite of the big
\'etale site of $\Spec \bF_p$ (the {\it smooth big \'etale site} of $\Spec
\bF_p$); for
simplicity, we denote these objects by
$\bZ_l(n)^c$ and $\bQ_l(n)^c$. Let $G$ still denote the absolute Galois group of
$\bF_p$.

\begin{prop}\label{p4.1} Let $X$ be a smooth, projective variety over $\bF_p$. Then,
for
$n\in\bZ$, 
\[H^i_\cont(X,\bQ_l(n))=\begin{cases}0&\text{if $i\neq 2n, 2n+1$}\\
H^{2n}(\overline X,\bQ_l(n))^{G}&\text{if $i=2n$}\\
H^{2n}(\overline X,\bQ_l(n))_{G}&\text{if $i=2n+1$.}
\end{cases}\]
\end{prop}

\begin{pf} By Deligne's proof of the Weil conjectures
\cite{deligne}, the eigenvalues of the Frobenius action on $H^q(\overline X,\bQ_l)$
are algebraic integers whose infinite absolute values are all equal to $q/2$. The
proposition then follows from proposition \ref{p2.6}.\end{pf}

\begin{cor}\label{c4.2} $H^i_\cont(X,\bZ_l(n))$ and $H^i_\et(X,\bQ_l/\bZ_l(n))$ are
finite for $i\neq 2n,2n+1$.
\end{cor}

\begin{pf} The first claim follows from proposition \ref{p4.1} and lemma
\ref{l2.2}; the second one follows from proposition \ref{p4.1}, lemma \ref{l2.2}
and the long exact sequence
\[\dots \to H^i_\cont(X,\bZ_l(n))\to H^i_\cont(X,\bQ_l(n))\to
H^i_\et(X,\bQ_l/\bZ_l(n))\to H^{i+1}_\cont(X,\bZ_l(n))\to \dots\]
stemming from lemma \ref{l2.1}.\end{pf}

\begin{thm}\label{t4.1} Let $\pi$ be the
projection of $\sS$ onto the small \'etale site of $\Spec \bF_p$.
Then the natural map
\[\pi^*\bZ_l(n)^c_{|\bF_p}=\pi^*\pi_*\bZ_l(n)^c@>\alpha>> \bZ_l(n)^c\]
is an isomorphism for $n\le 0$. ($\bZ_l(n)^c$ is ``locally constructible".)
\end{thm}

\begin{pf} Let
$K(n)$ be the cone of $\alpha$. We have to prove that
$H^*(X,K(n))=0$ for any $X\in \sS$. By lemma \ref{l2.1}, $K(n)$ has uniquely
divisible cohomology sheaves. By theorem \ref{t2.6} (or rather its proof), $X\mapsto
H^*(X,K(n))$ therefore defines a pure graded cohomology theory with transfers with
values in
$\bQ$-vector spaces, in the sense of definition \ref{d3.2}. (To see that it is pure,
it is enough to see that the complexes of sheaves
$\pi^*\bQ_l(n)^c_{|\bF_p}$ verify purity for $n\le 0$, which is trivial from theorem
\ref{t2.2}).  By theorem
\ref{t3.1} b), it is enough to prove the claim for $X$ smooth and projective. In
this case, it follows immediately from proposition \ref{p4.1} and theorem
\ref{t2.2}.\end{pf}

\begin{cor}\label{c4.7} Over $\sS$, we have
\[\bQ_l(0)^c\simeq \bQ_l[0]\oplus \bQ_l[-1].\]
\end{cor}

\begin{pf} This follows from corollary \ref{c2.3}, theorem \ref{t2.2} and theorem
\ref{t4.1}.
\end{pf}

Recall from section \ref{3} the canonical generator $e\in
H^1_\cont(\bF_p,\bZ_l)$. Cup-product by $e$ induces endomorphisms of degree $1$ and
square
$0$
\begin{equation}\label{eq10}\bZ_l(n)^c@>e>> \bZ_l(n)^c[1].\end{equation}

\begin{prop}{\rm (\cf \cite[prop. 6.5]{milne})}\label{p5.8} For any scheme $X$ of
finite type over
$\bF_p$, the diagram
\[\begin{CD}
H^q_\cont(\overline X,\bZ_l(n))^G@>>>H^q_\cont(\overline X,\bZ_l(n))_G\\
@AAA @VVV\\
H^q_\cont(X,\bZ_l(n)) @>\cdot e>> H^{q+1}_\cont(X,\bZ_l(n))
\end{CD}\]
commutes, where the top horizontal map is the natural one and the vertical maps come
from the Hochschild-Serre spectral sequence. The same holds for continuous
cohomology with proper supports and continuous homology.
\end{prop}

\begin{pf} This is an immediate consequence of the multiplicativity of the
Hochschild-Serre spectral sequence (view $e$ as an element of
$H^1(\bF_p,H^0_\cont(\overline X,\bZ_l))\subseteq H^1_\cont(X,\bZ_l)$) and the
following lemma:

\begin{lemma} \label{l4.1} Let $A$ be a continuous $G$-module. Then the
diagram
\[\begin{CD}
H^0(G,A) @>\cdot e>> H^1(G,A)\\
|| && @V{\wr}VV\\
A^G@>>> A_G
\end{CD}\]
commutes. Here, the bottom horizontal map is the composition $A^G\inj A\surj A_G$.
\end{lemma}

Indeed, the right vertical isomorphism is induced by the map sending a continuous
cocycle $c$ to $c(F)$. On the other hand, cup-product by $e$ is defined at the
level of cocycles by the formula $(a\cdot e)(g)=e(g)a$.\end{pf}

Recall condition
$S^n$ from \cite{tate2}:

\begin{defn}\label{d4.2} a) Let $G$ be a profinite group and $A$ a topological
$G$-module. We say that $A$ is {\em semi-simple at $1$} if the composition
\[A^G\Inj A \Surj A_G\]
is bijective.\\
b) Let $k$ be a field and $X/k$ be smooth and projective. Let $k_s$ be a separable
closure of $k$, $G_k=Gal(k_s/k)$, $\overline
X=X\otimes_{k} \bar k_s$ and $n\in \bZ$. We say that $X$ {\em satisfies condition
$S^n(X)$} if the topological $G_k$-module
$H^{2n}_\cont(\overline X,\bQ_l(n))$  is semi-simple at $1$.
\end{defn}

\begin{cor}\label{c4.4} Let $X$ be smooth, projective over $\bF_p$. Then, condition
$S^n(X)$ is equivalent to the bijectivity of $H^{2n}_\cont(X,\bQ_l(n))\allowbreak
@>\cdot e>> H^{2n+1}_\cont(X,\bQ_l(n))$.\qed
\end{cor}

In the remainder of this section, we shall freely use the notion of a $\bQ_l$-sheaf
appearing in Deligne \cite{weilII}; we refer to this article for a precise
definition.

\begin{thm} \label{t4.2} Let $\sF$ be a $\bQ_l$-sheaf over $\bF_p$, pure of weight
$m$. Consider $\sF$ as a $\bQ_l$-sheaf over the big \'etale site of $\bF_p$ by
pull-back. Then,\\
a) For any smooth variety
$X/\bF_p$ and
$q\ge 0$, $H^q_\cont(\overline X,\sF)$ is mixed of weights between $m+q$ and
$m+2q$.\\
b) For any variety $X/\bF_p$ of dimension $d$ and $q\ge 0$, $H^q_\cont(\overline
X,\sF)$ is mixed of weights between $m+2q-2d$ and
$m+2q$.\\
c) If $\sF\spcheck(-m)$ is entire, then for any $X/\bF_p$ of dimension $d$ and
$q\ge 0$,
$H^q_\cont(\overline X,\sF)$ is mixed of weights $\le 2m+2d$.\end{thm}

In c), $\sF\spcheck$ denotes the dual of $\sF$ and ``entire" means that the
eigenvalues of the action of Frobenius are algebraic integers.

\begin{pf} a) The bound $m+q$ follows from \cite[cor. 3.3.5]{weilII}. To prove the
other, we consider the graded pure cohomology theory with transfers (definition
\ref{d3.2})
\[h^i(X,n)=H^i_\cont(\overline X,\sF(n))\]
with values in the category $\sA$ of $\bQ_l$-sheaves over $\bF_p$. Let $\sB$ be
the thick subcategory consisting of sheaves of weights $\le  m$. We want to prove
that $h^i(X,n)\in \sB$ for $i\le n$. By theorem \ref{t3.2}, we reduce to the case
where $X$ is projective; then it follows from the main result of
\cite{deligne}.

b) The proof should be along the same lines, but we could not find a suitable
formalisation. Suppose first
$X$ smooth. By the known $l$-cohomological dimension of $\overline X$, we may
assume $q\le 2d$. Then the lower bound follows from a) since
$m+2q-2d\le m+q$. 

In general, we argue by
induction on
$\dim X$. Recall \cite{dJ} that an {\it alteration} is a proper morphism
$f:X'\to X$ such that
$f^{-1}(U)\to U$ is finite and flat for some nonempty open subset $U$ of $X$. Let
$f:X'\to X$ be an alteration, with $X'$
smooth over $\bF_p$: it exists by \cite{dJ}. Let $U$ be the maximal open subset of
$X$ such that $U'=f^{-1}(U)$ is flat over $U$, and $Z=X\setminus U$ (the ``center"
of the alteration). Set $Z'=f^{-1}(Z)$. So we have a diagram
\[\begin{CD}
Z'@>i'>> X'@<j'<< U'\\
@V{f'}VV @V{f}VV @V{f''}VV\\
Z@>i>> X@<j<< U.
\end{CD}\]

We argue as in the
proof of corollary
\ref{c2.1}. By proper base change, we have a commutative diagram of long exact
sequences
\[\begin{CD}
\dots &\to& H^{q-1}_\cont(\overline Z,\sF)&\to& H^q(\overline
X,j_!\sF_U)&\to& H^q_\cont(\overline X,\sF)&\to&
H^q_\cont(\overline Z,\sF)&\to&\dots\\
&&@VVV @VVV @VVV @VVV\\
\dots &\to& H^{q-1}_\cont(\overline{Z'},\sF)&\to& 
H^q(\overline X,j_!Rf''_*{f''}^{*}\sF_{U'})&\to&
H^q_\cont(\overline{X'},\sF)&\to& H^q_\cont(\overline{Z'},\sF)&\to&\dots
\end{CD}\]

By induction and the smooth case, the sheaves
$H^q_\cont(\overline{X'},\sF)$, $H^q_\cont(\overline Z,\sF)$ and
$H^{q-1}_\cont(\overline{Z'},\sF)$ are all of weights $\le m+2q$. Therefore
$H^q(\overline X,j_!Rf''_*{f''}^{*}\sF_{U'})$ is of weights $\le m+2q$ as
well. On the other hand, since
$f''$ is finite and flat, there is a trace morphism $Rf''_*{f''}^{*}\sF_{U'}\to
\sF_U$ whose composition on the left with the pull-back map is
multiplication by $\deg f''$. Therefore $H^q(X,j_!\sF_U)$ is of weights $\le m+2q$
as well and finally $H^q_\cont(X,\sF)$ is of weights $\le m+2q$, as desired.
For the lower bound we argue similarly, noting that $m+2(q-1)-2\dim Z'\ge m+2q-2d$.

c) If $X$ is smooth, this follows from \cite[cor. 3.3.3]{weilII} and Poincar\'e
duality (compare \cite[proof of th. 2]{jannsen4}); in general, we proceed as in
the proof of b).
\end{pf}

\begin{cor}\label{c4.3} 
a) The sheaf $\sH^i(\bQ_l(n)^c)$ is $0$ for $i<n$.\\
b) For all $X\in Sch/\F_p$,
the groups
$H^i_\cont(X,\bZ_l(n))$ are finite for
$i\notin [n,n+d+1]$, and for $i\notin [n,2n+1]$ if $X$ is smooth.\\
c) For all $X\in Sch/\bF_p$, the groups $H^i_\cont(X,\bQ_l/\bZ_l(n))$ are finite
for
$i\notin [n,n+d+1]$, and for $i\notin [n,2n+1]$ if $X$ is smooth.
\end{cor}

\begin{pf} This follows from theorem \ref{t4.2}, applied to
$\sF=\bQ_l(n)$, proposition \ref{p2.6} and lemma \ref{l2.2} b).
\end{pf}

\begin{defn}\label{d6.1}
For $d\ge 0$, let $Sch_d$ be the category of schemes of finite type of
dimension $\le d$ over $\F_p$ and $\pi_d:(\Spec \bF_p)_\Et=(Sch/\bF_p)_\et\to
(Sch_d)_\et$ the projection of the big \'etale site of $\bF_p$ over $Sch_d$
provided with the \'etale topology. If $C\in \sD((\Spec \bF_p)_\Et)$, we denote
$R(\pi_d)_*C$ simply by $C_{|Sch_d}$.
\end{defn}

\begin{cor}\label{c4.5}  For $n>d$, $\bQ_l(n)^c_{|Sch_d}=0$ and
$\bZ_l(n)^c_{|Sch_d}=\bQ_l/\bZ_l(n)[-1]$.
\end{cor}

\begin{pf} By theorem \ref{t4.2} c), applied to $\sF=\bQ_l$, for any $X\in Sch_d$
and any
$q\ge 0$,
$H^q_\cont(\overline X,\bQ_l(n))$ is mixed of weights $\le 2d-2n$. By proposition
\ref{p2.6}, this implies that $H^*_\cont(X,\bQ_l(n))=0$ for $n>d$. The
statement for
$\bZ_l(n)^c$ follows from this by lemma \ref{l2.1}. \end{pf}

\begin{cor}\label{c4.6} For any $X\in Sch/\bF_p$, one has
$H_*^\cont(X,\bQ_l(n))=H^*_{c,\cont}(X,\bQ_l(n))=0$ for any $n<0$.
\end{cor}

\begin{pf} If $X$ is smooth, this follows from corollary \ref{c4.5} and
proposition \ref{p2.2} c). In general, it follows from this, theorem
\ref{t3.3} (for homology) and theorem \ref{t2.3} (for cohomology with proper
supports).
\end{pf}

\begin{cor}\label{t5.3} Let $f:X\to U$ be a smooth, projective morphism of
varieties over
$\bF_p$. Let $\delta=\dim U$ and $d$ be the relative dimension of $f$. Then we have
\[H^p_\cont(U,R^q f_*^\ladic \bQ_l(n))=0\quad\text{for}\quad
p+q<n,p+d<n,  d+\delta<n\quad\text{or} \quad p+q>2n+1.\]
If $U$ is smooth, this group vanishes for $\delta+q<n$ as well. Finally, if $f$ is
``defined over
$\bF_p$",
\ie
$X=X_0\times_{\bF_p}U$ with
$X_0$ smooth and projective over $\bF_p$, it also vanishes for $2p+q<2n$.
\end{cor}

\begin{pf}  By Hard Lefschetz \cite[th. 4.1.1]{weilII}, smooth and proper base
change and Deligne's degeneracy criterion for spectral sequences, the Leray
spectral sequence degenerates at $E_2$, hence the group
$H^p_\cont(U,R^q f_*^\ladic\bQ_l(n))$ is a direct summand of $H^{p+q}(X,\bQ_l(n))$
(compare \cite[\S 4]{deligne2}). The vanishings for $p+q<n$ and $p+q>2n+1$ then
follow from corollary \ref{c4.3}. The vanishing for $d+\delta<n$ follows from
corollary \ref{c4.5}. The vanishing for $p+d<n$ follows from Hard Lefschetz,
which implies that $R^q f_*^\ladic \bQ_l(n)\simeq
R^{2d-q}f_*^\ladic\bQ_l(n+d-q)$ (apply the inequality $p+q<n$ with these new
values of
$p,q,n$). For the bound
$\delta+q<n$, we note that the result of theorem \ref{t4.2} c) in the smooth case
does not necessitate
$\sF$ to be defined over
$\bF_p$. Finally, the bound $2p+q<2n$ follows from theorem
\ref{t4.2} b), this time applied with
$\sF=H^q_\cont(\overline X_0,\bQ_l(n))$, which is pure of weight $q-2n$ by the main
result of \cite{deligne}, and proposition \ref{p2.6}.
\end{pf}

\begin{rk} 
This extends Jannsen's inequalities in
\cite[th. 1]{jannsen5}: {\it loc. cit.} says that, for $U$ a smooth curve over
$\bF_p$ and
$X$ as above, $H^2_\cont(\pi_1(U),H^q(\overline X,\bQ_l(n)))=0$ whenever
$q+2\le n$ or
$q+2>2n+1$. 
\end{rk}

Since continuous \'etale cohomology acts by products on continuous \'etale
cohomology with supports and \'etale homology, for any $X\in Sch/\F_p$
we have complexes
\begin{equation}\label{eq4.1}\begin{CD}
\dots &\to& H^{i-1}_\cont(X,\bQ_l(n))@>\cdot e>>
H^i_\cont(X,\bQ_l(n))@>\cdot e>>
H^{i+1}_\cont(X,\bQ_l(n))&\to&\dots\\
\dots &\to& H^{i-1}_{c,\cont}(X,\bQ_l(n))@>\cdot e>>
H^i_{c,\cont}(X,\bQ_l(n))@>\cdot e>>
H^{i+1}_{c,\cont}(X,\bQ_l(n))&\to&\dots\\
\dots &\to& H_{i+1}^\cont(X,\bQ_l(n))@>\cdot e>>
H_i^\cont(X,\bQ_l(n))@>\cdot e>>
H_{i-1}^\cont(X,\bQ_l(n))&\to&\dots
\end{CD}\end{equation}

\begin{thm}\label{t4.6} a) If $S^n$ holds for any smooth, projective
variety over $\F_p$, then the last two complexes of \eqref{eq4.1} are
acyclic for any $X\in Sch/\F_p$, and the first one is acyclic for $X$
smooth of pure dimension $d$ of one replaces $n$ by $d-n$.\\
b) If $S^n$ holds for $n\le d$ and any smooth, projective variety of
dimension $\le d$, then the first complex of \eqref{eq4.1} is acyclic for
$n\le d$ and any $X\in Sch_d$. \end{thm}

\begin{pf} a) First assume that $X$ is smooth projective. Then the first
sequence of \eqref{eq4.1} is exact by proposition \ref{p4.1} and corollary
\ref{c4.4}; so is the second one by lemma \ref{l2.3} and so is also the
third one by proposition \ref{p2.2} c) and theorem \ref{t2.3}. 

In general, we first deal with continuous homology. Let $X\in Sch/\F_p$,
$Z$ a closed subset and $U$ the open complement. Proposition \ref{p2.2} a)
and a little inductive argument involving a big diagram chase shows that,
if theorem \ref{t4.6} holds for two among $X,Z,U$, then it holds for the
third. Moreover, if $f:U'\to U$ is a finite and flat morphism and theorem
\ref{t4.6} holds for $U'$, then it holds for $U$ by a transfer argument.
An argument as in the proof of theorem \ref{t3.3} now yields the
conclusion. 

The case of cohomology with proper supports follows from this and theorem
\ref{t2.3}; the case of cohomology for $X$ smooth follows from this and
proposition \ref{p2.2} c). 

For the sake of the proof of b), we note that the above argument works if
we restrict to varieties of dimension $\le d$ for some $d$.

b) By a) and the assumption, the first sequence of \eqref{eq4.1} is exact
for any $X$ smooth of dimension $\le d$. For an arbitrary $X$ of dimension
$\le d$, we can argue as in the proof of theorem \ref{t4.2} b), using the
ideas in the second paragraph of the proof of a).\end{pf}

\begin{cor}\label{c5.3} If $S^n$ holds for any
smooth, projective variety over
$\F_p$, then, for  any
$X\in Sch/\bF_p$ and any $i\in\bZ$,
$H^i_{c,\cont}(\overline X,\Q_l(n))$ and $H_i^\cont(\overline X,\Q_l(n))$ are
semi-simple at $1$. If $S^n$ holds for all $n\le d$ and all smooth, projective
varieties of dimension $\le d$, then
$H^i_\cont(\overline X,\Q_l(n))$ is also semi-simple at $1$ for all $n\le d$ and
$X\in Sch_d$.
\end{cor}

\begin{pf} This follows from theorem \ref{t4.6} and proposition \ref{p5.8}.\end{pf}

\begin{rk} One should compare this with \cite[th. 12.7]{jannsen2}, where a similar
result is proven with heavier assumptions.
\end{rk}

\section{Values of the zeta function}\label{zeta}

In this section, we generalise the $l$-primary part of \cite[th. 0.1]{milne} to
arbitrary varieties.

\begin{defn}\label{d4.1} a) Let $u:A\to B$ be a homomorphism of 
$\bZ_l$-modules with finite kernel and cokernel. The {\em index} of
$u$ is
\[\ind(u)={|\Ker u|\over |\Coker u|}.\]
b) Let $C$ be a bounded cochain complex of $\bZ_l$-modules with finite cohomology
groups. We set
\[\chi(C)=\prod_{i\in\bZ}|H^i(C)|^{(-1)^i}.\]
\end{defn}

Note that a) is a special case of b).

Let $X\in Sch/\bF_p$ of dimension $d$ and $\zeta(X,s)$ be its zeta
function. By \cite[expos\'e XV]{sga5}, we have $\zeta(X,s)=Z(X,p^{-s})$, where
\begin{equation}\label{eq20}
Z(X,t)=\prod_{i=0}^{2d} \det(1-tF|H^i_c(\overline
X,\bQ_l))^{(-1)^{i+1}}\end{equation}
where $F$ denotes the geometric Frobenius acting on the $\bQ_l$-adic cohomology
with proper supports of
$\overline X:=X\otimes_{\bF_p} \bar \bF_p$\footnote{Throughout this section, we drop
the index cont from continuous cohomology groups, since we shall not consider any
other cohomology.}. In this section, we shall prove:

\begin{thm}\label{t4.3} Let $r_i(n)=\dim_{\bQ_l} H^i_c(X,\bQ_l(n))$ and $a_i(n)$ be
the order of the zero of 
\linebreak[3] $\det(1-tF|H^i_c(\overline
X,\bQ_l))$ at
$t=p^{-n}$. Assume that $S^n(X)$ holds for any smooth, projective $X/\bF_p$ (\cf
definition \ref{d4.2}). Then, for any $X\in Sch/\bF_p$ and $i\ge 0$:\\
a) $a_i(n)=r_i(n)-r_{i-1}(n)+\dots$\\
b) Let
$\det(1-tF|H^i_c(\overline
X,\bQ_l))=(1-p^nt)^{a_i(n)}f_i(t)$. Then
\[|f_i(p^{-n})|_l=\left|\ind(H^i_c(\overline X,\bZ_l(n))^G\to H^i_c(\overline
X,\bZ_l(n))_G)\right|^{-1}_l\]
where $||_l$ is the $l$-adic absolute value.
\end{thm}

\begin{rk} Condition $S^n$ is necessary in theorem \ref{t4.3}, as is easily
seen from the special case where $X$ is smooth and projective.
\end{rk}

\begin{pf} a) The statement is clear from corollaries \ref{c5.3} and
\ref{l6.1}. The proof of b) will be divided into a series
of lemmas.

\begin{lemma}\label{l4.3}
a) Let $0\to C'\to C\to C''\to 0$ be a short exact
sequence of complexes satisfying the assumption of definition
\ref{d4.1} b). Then
\[\chi(C)=\chi(C')\chi(C'').\]
b) Let $A@>u>> B@>v>> C$ be a chain of
$\bZ_l$-homomorphisms such that $u$ and $v$ have finite kernel and cokernel. Then
\[\ind(v\circ u)=\ind(v)\ind(u).\]
\end{lemma}

\begin{pf} a) is classical; b) follows from the exact sequence
\[0\to\Ker u\to \Ker vu\to\Ker v\to \Coker u\to\Coker vu\to \Coker v\to 0.\]
\end{pf}

\begin{lemma}\label{l4.2} If $A=B$ in definition \ref{d4.1}, then
$|\ind(u)|_l=|\det(u\otimes \bQ_l)^{-1}|_l$.
\end{lemma}

\begin{pf} Lemma \ref{l4.3} a) shows that, if $A'\subset A$ is stable under $u$, we have
\[\ind(u)=\ind(u_{|A'})\ind(\bar u)\]
where $\bar u$ is the induced endomorphism of $A/A'$. Therefore we reduce to the
cases where $A$ is finite or torsion-free. In the first case, the claim is clear
in view of the exact sequence
\[0\to\Ker u\to A@>u>> A\to \Coker u\to 0.\]

In the second case, it follows \eg from \cite[ch. III, prop. 2]{CL}.\end{pf}

\begin{lemma}\label{l4.5} Under $S^n$, we have $\det(1-tF|H^i_c(\overline
X,\bQ_l))=(1-p^{n}t)^{a_i(n)}f_i(t)$, with
\[|f_i(p^{-n})|_l=\left|\ind\left({H^{i}_c(\overline X,\bZ_l(n))\over
H^i_c(\overline X,\bZ_l(n))^G}@>1-F>> {H^i_c(\overline X,\bZ_l(n))\over
H^i_c(\overline X,\bZ_l(n))^G}\right)^{-1}\right|_l.\]
\end{lemma}

\begin{pf} For simplicity, set $H=\displaystyle{H^i_c(\overline X,\bZ_l(n))\over
H^i_c(\overline X,\bZ_l(n))^G}$. The claim follows from lemma \ref{l4.2} and the
exact sequence of
$G$-modules
\begin{equation}\label{eq21}0\to H^i_c(\overline X,\bZ_l(n))^G\to H^i_c(\overline
X,\bZ_l(n))\to H\to 0.\end{equation}
\end{pf}

Theorem \ref{t4.3} b) now follows from lemma \ref{l4.5} by taking the cohomology of
\eqref{eq21}, which gives an exact sequence
\[0\to H^G\to H^i_c(\overline X,\bZ_l(n))^G\to H^i_c(\overline
X,\bZ_l(n))_G\to H_G\to 0.\]
\end{pf}

We shall need the following corollary in section \ref{6}:

\begin{cor}\label{c6.2} With notation as in theorem \ref{t4.3}, we also have
\[|f_i(p^{-n})|_l=\left|\ind(H_i(\overline X,\bQ_l/\bZ_l(n))^G\to H_i(\overline
X,\bQ_l/\bZ_l(n))_G)\right|_l.\]
If $X$ is smooth of pure dimension $d$, we have
\[|f_i(p^{-n})|_l=\left|\ind(H^{2d-i}(\overline X,\bQ_l/\bZ_l(d-n))^G\to
H^{2d-i}(\overline X,\bQ_l/\bZ_l(d-n))_G)\right|_l.\]
\end{cor}

\begin{pf} The first formula follows from the duality between $H^i_c(\overline
X,\bZ_l(n))$ and
$H_i(\overline X,\bQ_l/\bZ_l(n))$, while the second one follows from the
isomorphism between \'etale homology and \'etale cohomology for smooth
varieties.
\end{pf}

\begin{thm}\label{t6.1} With assumptions and notation as in theorem \ref{t4.3}, we
have\\ 
a) $a(n):=\ord_{s=n}\zeta(X,s)=\sum_{i\ge 0}(-1)^{i+1}(2d+1-i)r_i(n)$, where
$d=\dim X$.\\
b) Let
$\zeta(X,s)=(1-p^{n-s})^{a(n)}\varphi(s)$. Then
\[|\varphi(n)|_l=\left|\chi\left(H^*_c(X,\bZ_l(n)),e\right)\right|_l.\]
\end{thm}

\begin{pf} a) follows immediately from theorem \ref{t4.3} a). b) follows from
theorem \ref{t4.3} b) and a little diagram chase in the commutative diagram of exact
sequences (see proposition
\ref{p5.8})
\[\begin{CD}
&&0&&0&&0\\
&&@VVV @AAA @VVV\\
&\to&H^{i-2}_c(\overline X,\bZ_l(n))_G &&H^i_c(\overline X,\bZ_l(n))^G@>>>
H^i_c(\overline X,\bZ_l(n))_G\\ &&@VVV @AAA @VVV\\
\dots &\to& H^{i-1}_c(X,\bZ_l(n))@>\cdot e>>H^i_c(X,\bZ_l(n))@>\cdot e>>
H^{i+1}_c(X,\bZ_l(n))&\to&\dots\\
&&@VVV @AAA @VVV\\
&&H^{i-1}_c(\overline X,\bZ_l(n))^G@>>> H^{i-1}_c(\overline X,\bZ_l(n))_G&&
H^{i+1}_c(\overline X,\bZ_l(n))^G &\to&\\
&&@VVV @AAA @VVV&\qquad &\qquad &\\
&&0&&0&&0
\end{CD}\]
\end{pf}

\begin{cor}\label{c6.1} With assumptions and notation as in theorem \ref{t6.1}, we
have\\
\[|\varphi(n)|_l=\left|\prod_{i\ge 0}
|H^i_c(X,\bZ_l)_\tors|^{(-1)^i}\chi\left(\overline H^*_c(X,\bZ_l(n)),\bar
e\right)\right|_l\]
where $\overline H^i_c(X,\bZ_l(n))=H^i_c(X,\bZ_l(n))/\tors$ and $\bar e$ is the
induced homomorphism.
\end{cor}

\begin{pf} This follows from theorem \ref{t6.1} and lemma \ref{l4.3} a).\end{pf}

\begin{cor}\label{c4.1} Keep assumptions and notation as in corollary
\ref{c6.1}. Assume $X$ smooth and projective. Let $R_n(X)$ be the determinant of the
pairing
\[\overline H^{2n}(X,\bZ_l(n))\times \overline H^{2d-2n}(X,\bZ_l(d-n))\to
\overline H^{2d}(X,\bZ_l(d))\iso H^{2d}(\overline X,\bZ_l(d))\simeq \bZ_l\]
with respect to any bases of $\overline H^{2n}(X,\bZ_l(n))$ and $\overline
H^{2d-2n}(X,\bZ_l(d-n))$ (this is well-defined up to a unit of $\bZ_l$). Then
\[|\varphi(n)|_l=\left|\prod_{i\neq 2n,2n+1}
 |H^i(X,\bZ_l(n))|^{(-1)^{i}}\cdot
{|H^{2n}(X,\bZ_l(n))_\tors|\over
|H^{2n+1}(X,\bZ_l(n))_\tors|R_n(X)}\right|_l.\]
\end{cor}

\begin{pf} We have a commutative diagram
\[\begin{CD}
\overline H^i_c(X,\bZ_l(n))@>\bar e>> \overline H^{i+1}_c(X,\bZ_l(n))\\
@V{\alpha}VV @V{\beta}VV\\
Hom(\overline H^{2d-2n}(X,\bZ_l(d-n)),\overline H^{2d})@>Hom(Id,\bar
e)>> Hom(\overline H^{2d-2n}(X,\bZ_l(d-n)),\overline H^{2d+1})
\end{CD}\]
where $\overline H^{2d}$ and $\overline H^{2d+1}$ respectively stand for
$H^{2d}(X,\bZ_l(d))$ and $H^{2d+1}(X,\bZ_l(d))$ (in order for the diagram to fit
in the page). By lemma \ref{l4.2}, we have $|\ind(\alpha)|_l=|R_n(X)^{-1}|_l$ and
$\beta$ is an isomorphism by theorem
\ref{t2.3}, hence $\ind(\beta)=1$. Finally, $\overline
H^{2d}@>\bar e>> \overline H^{2d+1}$ is also an isomorphism by proposition
\ref{p5.8}. Applying lemma \ref{l4.3} b) once again, we get
\[|\ind(\bar e)|_l=|R_n(X)^{-1}|_l\]
and the result now follows from corollary \ref{c6.1}.\end{pf}

Corollary \ref{c4.1} is a variant of \cite[th. 0.1]{milne}.

\section{The conjecture}\label{coarse}

In this section and the next one, it will be occasionally convenient to use the
following definition:

\begin{defn}\label{d7.1}For any integer $d\ge 0$, we denote by $\sS_d$ the
full subcategory of
$Sm/\bF_p$ formed of smooth varieties of dimension $\le d$; we consider $\sS_d$ as a
subsite of
$Sm/\bF_p$ (for the \'etale or Zariski topology).
\end{defn}

\subsection{Milnor's $K$-theory} For any commutative ring $R$, denote by $K_n^M(R)$
the group defined by generators and relations as in the case of a field
\cite{milnor} (to be on the safe side, include the relation $\{x,-x\}=0$). The
Zariski sheaf
$\sK_n^M$ associated to the presheaf $\Spec R\mapsto K_n^M(R)$ is called the {\it
sheaf of Milnor
$K$-groups}: it is a sheaf on the big Zariski site of $\Spec \bZ$. We have

\begin{thm}\label{t5.1} Suppose $X$ is a smooth variety over a field $k$. Then\\
a) There is a complex of Zariski sheaves over $X$:
\[0\to \sK_n^M\to \coprod_{x\in X^{(0)}} (i_x)_*
K_n^M(k(x))\to
\coprod_{x\in X^{(1)}} (i_x)_* K_{n-1}^M(k(x))\to\dots  \]
which is exact, except perhaps at $\sK_n^M$, where its kernel is killed by
$(n-1)!$. Here $X^{(p)}$ denotes the set of points of $X$ of codimension $p$,
$i_x:\{x\}\inj X$ is the natural immersion and $K_n^M(k(x))$ is considered as a
constant sheaf on $x_\Zar$. In particular the sequence
\[0\to \sK_n^M \otimes\bQ\to \coprod_{x\in X^{(0)}} (i_x)_*
K_n^M(k(x))\otimes\bQ\to
\coprod_{x\in X^{(1)}} (i_x)_* K_{n-1}^M(k(x))\otimes\bQ\to\dots  \]
is exact.\\
b) We have
\[H^p_\Zar(X,\sK_n^M)\otimes \bQ=
\begin{cases}
0& \text{for $p>n$}\\
CH^n(X)\otimes \bQ&\text{for $p=n$.}
\end{cases}\]
\end{thm}

\begin{pf} a) The complex was defined by Kato \cite{kato}. By Rost
\cite[th. 6.1]{rost}, it is acyclic, except perhaps at $\sK_n^M$ and
$\coprod_{x\in X^{(0)}} (i_x)_* K_n^M(k(x))$. By a result of Gabber
(unpublished), it is exact at $\coprod_{x\in X^{(0)}} (i_x)_* K_n^M(k(x))$
as well. Finally, $\sK_n^M$ maps to the sheaf $\sK_n$ of algebraic
$K$-groups with kernel killed by $(n-1)!$, by results of Suslin \cite{suslin} and
Guin \cite{guin}. The claim on $\Ker(\sK_n^M\to \coprod_{x\in X^{(0)}} (i_x)_*
K_n^M(k(x)))$ follows from this and Gersten's conjecture for the algebraic
$K$-theory of
$X$ (Quillen \cite[th. 7.5.11]{quillen}). The statement after tensoring by $\bQ$,
hence b), readily follow.\end{pf}

\begin{rk} If one is interested only in $\sK_*^M\otimes\bQ$, one can get another
proof of a) by using the fact that this complex is the weight $n$-part for the Adams
operations of the corresponding Gersten complex for algebraic $K$-theory, \cf
\cite{soule}. Moreover, b) holds without tensoring by $\bQ$ for $n=\dim X$, see
\cite{kato}.
\end{rk}

\begin{cor}\label{c7.8} Let $k$ be a field and $f:Y\to X$ a finite flat
morphism of smooth $k$-schemes,
of constant separable degree $d$. Then there exists a morphism
\[Tr_f:f_*\sK_n^M \otimes\bQ\to \sK_n^M \otimes\bQ\]
whose left composition with the natural morphism $\sK_n^M \otimes\bQ\to
f_*\sK_n^M \otimes\bQ$ is multiplication by $d$. If $g:Z\to Y$ is another such
morphism, we have
\[Tr_{f\circ g}= Tr_f\circ f_* Tr_g.\]
\end{cor}

\begin{pf} This follows from the result in the case of a finite field extension
(Bass-Tate \cite[\S 5]{bt}, Kato \cite[prop. 5]{kato2}) and theorem \ref{t5.1} b).
One should use \cite[lemma 16]{kato2}, which implies that the Milnor $K$-theory
transfer is compatible with the Gersten resolution of theorem \ref{t5.1} a).
\end{pf}

\begin{cor}\label{c5.1} Let $(X,Z)$ be a smooth pair of codimension $c$ of
$k$-schemes, where
$k$ is a field, and
$i:Z\inj X$ the corresponding closed immersion. Then
\[Ri^!(\sK_n^M)_X \otimes\bQ\simeq (\sK_{n-c}^M)_Z \otimes\bQ[-c].\]
\end{cor}

\begin{pf} Applying the functor $i^!$ to the flasque resolution of
$(\sK_n^M)_X \otimes\bQ$ given by theorem \ref{t5.1} a), we get the
corresponding resolution of $(\sK_{n-c}^M)_Z \otimes\bQ$ shifted by $c$ to
the right.\end{pf}

\begin{cor}\label{c5.2} If $k$ is a field, the map over the big smooth Zariski
site of
$\Spec k$
\[\sK_n^M \otimes\bQ\to R\alpha_*\alpha^*\sK_n^M \otimes\bQ\]
is an isomorphism, where $\alpha$ is the projection of the big \'etale site
onto the big Zariski site.
\end{cor}

\begin{pf} Over the small Zariski site of $\Spec k$, we have
$(R^q\alpha_*\alpha^*\sK_n^M \otimes\bQ)(k)= H^q(k, K_n^M(k_s)
\otimes\bQ)=0$ for $q>0$ and $\sK_n^M(k) \otimes\bQ\iso
(\alpha_*\alpha^*\sK_n^M \otimes\bQ)(k)=(K_n^M(k_s) \otimes\bQ)^{G_k}$
by a transfer argument, where $k_s$ is a separable closure of $k$ and
$G_k=Gal(k_s/k)$. Over the small \'etale site of a smooth $k$-scheme $X$,
the isomorphism of corollary \ref{c5.2} follows from this (applied to all residue
fields of $X$) and theorem
\ref{t5.1} a).  \end{pf}

Corollary \ref{c5.2} implies that the map
\[H^i_\Zar(X,\sK_n)\otimes \bQ\to H^i_\et(X,\sK_n)\otimes \bQ\]
is an isomorphism for all $i,n$; we shall take this opportunity to write both
groups $H^i(X,\sK_n)\otimes \bQ$.

\begin{cor}\label{c7.5} If $k$ is a field, over the big smooth Zariski site of
$\Spec k$, we have
\[\sK_n^M\otimes\Q\iso R\pi_*\pi^*\sK_n^M\otimes\Q\]
where $\pi$ is the projection $\af_k\to\Spec k$.
\end{cor}

\begin{pf} This follows from \cite[prop. 8.6]{rost}.\end{pf}

\subsection{$\sK^M$-homology} Let $Z$ be an excellent scheme. For all $n\in \bZ$, we
denote by $\sM_{n,Z}$ the (chain) complex of Zariski sheaves
\[\dots @>>> \coprod_{x\in Z_{(i)}}(i_x)_*K_{n+i}^M(\kappa(x))@>>> \coprod_{x\in
Z_{(i-1)}}(i_x)_*K_{n+i-1}^M(\kappa(x))@>>>\dots\]
defined by Kato \cite{kato}. Since its terms are flasque, its hypercohomology is
computed by the cohomology of its global sections\footnote{Recall that the
hypercohomology of a chain complex $C_i$ is the hypercohomology of the cochain
complex $C^i$ defined by $C^i=C_{-i}$.}. In other terms,
\[\bH^{-i}_\Zar(Z,\sM_{n,Z})=A_i(Z,K^M_*,n)\]
where the right hand group is the one defined in \cite[\S 5]{rost}. For the
purpose of future compatibility with \'etale homology (see subsection \ref{8.4}), we
set:

\begin{defn}\label{d8.1} $H_i(Z,\sK_n^M)=\bH^{-i}_\Zar(Z,\sM_{-n,Z})$.
\end{defn}

We collect in the following proposition some properties of $\sK^M$-homology.

\begin{prop}\label{p8.3}
a) If $X$ is smooth over a field of pure dimension $d$, one has a quasi-isomorphism
\[\sM_{n,X}\otimes\bQ\simeq \sK_{n+d,X}^M\otimes \bQ[d].\]
In particular,
\[H_i(X,\sK_n^M)\otimes \bQ\simeq H^{d-i}_\Zar(X,\sK_{d-n}^M)\otimes \bQ.\]
b) Let $i:Z\inj Y$ be a closed immersion. One has
\begin{equation}\label{eq8.2} Ri^!\sM_{n,Y}=\sM_{n,Z}.
\end{equation}
c) The assignment $Z\mapsto \sM_{n,Z}$ is covariant for proper morphisms and
contravariant for flat, equidimensional morphisms of schemes of finite type over a
field.
\end{prop}

\begin{pf} a)  follows from theorem \ref{t5.1}; b) is obvious and c) follows from 
\cite[(4.5)]{rost}.
\end{pf}

\subsection{The conjecture} For $n=1$, $\alpha^*\sK_n^M=\bG_m$; the Kummer exact
sequences
\[1\to \mu_{l^\nu}\to \bG_m@>l^\nu>> \bG_m\to 1\quad (\nu\ge 1)\]
induce a ``Kummer" morphism
\[\bG_m[0]\to \bZ_l(1)^c[1]\]
in the derived category, and in particular a morphism of sheaves
\[u_1:\alpha^*\sK_1^M@>>> \sH^1(\bZ_l(1)^c).\]

We denote by $r\mapsto (r)$ the induced morphism $R^*\to H^1_\cont(\Spec
R,\bZ_l(1))$ for any ring $R$ containing $1/l$.

\begin{prop}\label{p5.10} a) The morphism $u_1$ and cup-product define morphisms
\[u_n:\alpha^*\sK_n^M@>>> \sH^n(\bZ_l(n)^c).\]
b) $u_n$ is compatible with the transfer of corollary \ref{c7.8} and the direct
image morphism of theorem \ref{t2.6} in continuous \'etale cohomology.
\end{prop}

\begin{pf} a) It is enough to show that, for any affine scheme
$U=\Spec R$ and any $r\in R^*$, the cup-product
\[(r)\cdot (1-r)\]
is $0$ in $H^2_\cont(U,\bZ_l(2))$. To do this we just mimic Tate's classical
argument \cite[th. 3.1]{tate}. We use the \'etale covering
$U_1=\Spec R_1$, with $R_1=R[t]/(t^l-r)$, noting that, if $r_1$ is the image of $t$
in
$R_1$, $r_1^l=r$ and $N_{R_1/R}(1-r_1)=1-r$, hence
\[(r)\cdot (1-r)=(r)\cdot f_*(1-r_1)=f_*(f^*(r)\cdot
(1-r_1))=lf_*((r_1\cdot(1-r_1))\]
where $f$ is the projection $U_1\to U$. This shows that $(r)\cdot (1-r)$ is
contained in a divisible subgroup of
$H^2_\cont(U,\bZ_l(2))$; by \cite[(4.9)]{jannsen}, we conclude that $(r)\cdot
(1-r)=0$.

b) This follows from theorem \ref{t2.7} and the construction of the transfer on
the Milnor $K$-sheaves (corollary \ref{c7.8}).
\end{pf}

By corollary \ref{c4.3} a), the morphism of proposition \ref{p5.10}
tensored with $\bQ$ defines a morphism
\[\alpha^*\sK_n^M[-n]\otimes\bQ\to\bQ_l(n)^c\]

Tensoring this with\
$\bQ_l(0)^c$ and using the pairing
\[\bQ_l(0)^c\oo^L\bQ_l(n)^c\to \bQ_l(n)^c\]
(\cf \eqref{eq2.2}), we get a morphism in
$\sD^+(Ab(\Spec \bF_p)_\Et)$
\begin{equation}\label{eq33}\bQ_l(0)^c\oo^L \alpha^*\sK_n^M[-n]\to
\bQ_l(n)^c.\end{equation}

Let $X$ be a smooth variety. Applying $H^{2n}_\et(X,-)$ to \eqref{eq33} we get a map
\[H^n_\et(X,\alpha^*\sK_n^M)\otimes\bQ_l\to H^{2n}_\cont(X,\bQ_l(n)).\]

By theorem \ref{t5.1} b) and corollary \ref{c5.2}, the left hand side is
isomorphic to $CH^n(X)\otimes\bQ_l$, hence we get a map
\[CH^n(X)\otimes\bQ_l\to H^{2n}_\cont(X,\bQ_l(n)).\]

\begin{lemma}\label{l7.1} This map is the $l$-adic cycle map of
\cite[(6.14)]{jannsen} tensored with $\bQ$.
\end{lemma}

\begin{pf} In view of the way the cycle class is defined (\cite[cycle 2.2.8]{sga},
\cite[(3.23)]{jannsen}), we reduce to $n=1$, in which case the result is
trivial.\end{pf}

\begin{conj}\label{co7.1} The restriction of \eqref{eq33} to the smooth big \'etale
site is an isomorphism.
\end{conj}

Conjecture \ref{co7.1} is trivially true for $n=0$. In more
concrete terms, it implies:

\begin{prop}\label{co7.2} If conjecture \ref{co7.1} holds, then there is an exact
triangle
\[\bZ_l(n)^c\to \alpha^*\sK_n^M\otimes\bQ_l[-n]\oplus
\alpha^*\sK_n^M\otimes\bQ_l[-n-1]\to
\bQ_l/\bZ_l[0]\to \bZ_l(n)^c[1].\]
In particular, for any
$X\in
\sS$, there is a long exact sequence
\begin{multline*}\dots\to H^i_\cont(X,\bZ_l(n))\to  H^{i-n}(X,\sK_n^M)\otimes
\bQ_l\oplus H^{i-n-1}(X,\sK_n^M)\otimes
\bQ_l\\
\to H^i_\et(X,\bQ_l/\bZ_l(n))\to H^{i+1}_\cont(X,\bZ_l(n))\to\dots\end{multline*}
\end{prop}

\begin{pf} This follows from lemma \ref{l2.1}, corollary \ref{c4.7} and corollary
\ref{c5.2}.
\end{pf}

\begin{rk} Conjecture \ref{co7.1} extended to all of $Sch/\bF_p$ is
false. For example, continuous \'etale cohomology is homotopy invariant on
$Sch/\bF_p$, while $\sK^M$-cohomology is not. See subsection \ref{8.4} for a {\em
homological} formulation.
\end{rk}

\begin{prop}\label{p7.1} Suppose conjecture \ref{co7.1} holds. Then,
for any smooth, projective variety $X/\bF_p$,
\begin{itemize}
\item $H^i(X,\sK_n^M)$ is torsion for $i<n$.
\item The cycle map $CH^n(X)\otimes\bQ_l\to
H^{2n}_\cont(X,\bQ_l(n))$ is an isomorphism.
\item Cup-product by $e$: $H^{2n}_\cont(X,\bQ_l(n))\to H^{2n+1}_\cont(X,\bQ_l(n))$
is an isomorphism. 
\end{itemize}
\end{prop}

\begin{pf} The first two claims directly follow from proposition \ref{p4.1}
and lemma
\ref{l7.1}. The last one follows from the same (obvious) fact with $\bQ_l(n)^c$
replaced by
$\bQ_l(0)^c\oo^L
\alpha^*\sK_n^M[-n]$ (\cf corollary \ref{c2.2}).\end{pf}

\begin{cor}\label{c7.1} Conjecture \ref{co7.1} implies the Tate conjecture
in codimension $n$, the equality of rational and homological equivalences
and condition $S^n(X)$ for any smooth, projective variety $X$ over
$\bF_p$.  \end{cor}

\begin{pf} Follows from proposition \ref{p7.1}, proposition \ref{p4.1} and
corollary \ref{c4.4}.\end{pf}

Conversely:

\begin{thm}\label{t7.1} Suppose that the Tate conjecture in codimension
$n$, the equality of rational and homological equivalences and condition
$S^n(X)$ hold for $n\le d$ and any smooth, projective variety $X$ over
$\bF_p$ of dimension $\le d$. Then conjecture \ref{co7.1}
holds for $n\le d$ on $\sS_d$ (see definition \ref{d7.1}).  \end{thm}

\begin{pf} Fix $n\le d$. Let $K(n)$ be the cone of \eqref{eq33}: we need
to show that the restriction of $K(n)$ to $\sS_d$ is $0$. We note that
$X\mapsto H^*(X,K(n))$ defines (for varying $n$) a graded pure cohomology
theory with transfers with values in $\bQ$-vector spaces in the sense of
definition \ref{d3.2}:  this follows from theorem \ref{t2.1}, corollary
\ref{c5.1} and proposition \ref{p5.10} b). Applying theorem \ref{t3.1} b),
we see that it is enough to prove that $\bH^*(X,\eqref{eq33})$ is an isomorphism for
any smooth, projective variety $X$ of dimension $\le d$. By the above arguments,
the only thing which is left to prove is the vanishing of $H^i(X,\sK_n^M)\otimes
\bQ$ for $i<n$. 

For this, we observe that the assumptions imply that numerical and
homological equivalences agree over $X$ (\cite[prop. 8.4]{milne},
\cite[(2.6)]{tate2}).  Therefore, the category of Chow motives of
dimension $\le d$ over $\bF_p$ is abelian and semi-simple \cite{jannsen3}
and we can apply the Soul\'e-Geisser argument \cite{soule2}, \cite[proof of th.  3.3
b)]{geisser}. We have to prove that, for any simple motive $M$ of
dimension $\le d$, we have \[H^i(M,\sK_n^M)\otimes
\bQ=0\quad\text{for}\quad i<n\quad\text{and}\quad n\le d.\]

Suppose first that $M$ is the Lefschetz motive $L^n$. Then $M$ is a direct
summand of $\bP^n_{\bF_p}$ and one sees easily that the conclusion is
valid. Suppose now that $M\neq L^n$. Let $P$ be the minimum polynomial of
the geometric Frobenius $F_M$ of $M$, viewed as an element of the (finite
dimensional $\bQ$-algebra) $End(M)$. By \cite[th. 2.8]{geisser}, we have
$P(p^m)\neq 0$. On the other hand, $F_M$ acts on $H^i(M,\sK_n^M)\otimes
\bQ$ by multiplication by $p^n$ \cite[prop. 2 (iv)]{soule2}.  So
$H^i(M,\sK_n^M)\otimes \bQ$ is killed by multiplication by a nonzero
rational number, and therefore is $0$.\end{pf}

In the opposite direction:

\begin{thm} \label{t7.3} Conjecture \ref{co7.1} holds if and only if, for any
finitely generated field $F$ of characteristic $p$,
\[\tilde H^i_\cont(F/\bF_p,\bQ_l(n))=
\begin{cases} 
0& \text{if $i\neq n,n+1$}\\
K_n^M(F)\otimes \bQ_l &\text{if $i=n,n+1$.}
\end{cases}\]
\end{thm}

\begin{pf} Necessity is clear. For sufficiency, let $K(n)$ be as in the proof of
theorem
\ref{t7.1}. The functor $(X,Z)\mapsto H^*_Z(X,K(n))$ defines a cohomology
theory with supports in the sense of \cite[def. 5.1.1]{cthk}, which verifies
axioms COH1 (\'etale excision), COH3 (homotopy invariance) and COH6 (transfers) 
of {\it loc. cit.}, 5.1, 5.3 and 6.2. By {\it loc. cit.}, theorem 6.2.5, for any
local ring $R$ of a smooth variety over $\bF_p$, with field of fractions $F$, the
map
\[H^*(R,K(n))\to H^*(F,K(n))\]
is injective. This concludes the proof.\end{pf}

We now give some further consequences of conjecture \ref{co7.1}.

\begin{prop}\label{p7.2} If conjecture \ref{co7.1} holds, then, for all
$X\in \sS$, $H^*_\cont(X,\bZ_l(n))$ carries a canonical (\ie independent
of $l$) integral structure.  \end{prop}

\begin{pf} Consider the morphism \[\bQ_l(0)^c\otimes\alpha^*\sK_n^M[-n]
\to \bQ_l/\bZ_l(n)[0]\] obtained by composing \eqref{eq33} with the
morphism $\bQ_l(n)^c\to \bQ_l/\bZ_l(n)[0]$. Writing $\bQ_l/\bZ_l(n)$ as a
direct summand of $(\bQ/\bZ)'(n):=\coprod_{l\neq p}\bQ_l/\bZ_l(n)$, it
factors as 
\[\bQ_l(0)^c\otimes\alpha^*\sK_n^M[-n]\to (\Q/\Z)'(n)[0] \to
\bQ_l/\bZ_l(n)[0].\]

Let $\bQ^c$ be as in definition \ref{d2.1}. The way
\eqref{eq33} is defined show that the composite morphism
\[\bQ^c\otimes\alpha^*\sK_n^M[-n]@>\otimes
\bQ_l>>\bQ_l(0)^c\otimes\alpha^*\sK_n^M[-n]\to (\Q/\Z)'(n)[0]\] does not depend on
$l\neq p$. If we denote by $C$ the fibre of this morphism, then conjecture
\ref{co7.1} says that there is an isomorphism \[\bZ_l(n)^c\iso \bZ_l\otimes C.\]
\end{pf}

\begin{prop}\label{p7.3} Suppose conjecture \ref{co7.1} holds. Let $X$ be
a smooth variety. Then\\
a) We have canonical isomorphisms for all $i\in\bZ$: 
\[H^i_\cont(X,\bQ_l(n))\simeq H^{i-n}_\Zar(X,\sK_n^M)\otimes\bQ_l\oplus
H^{i-n-1}_\Zar(X,\sK_n^M)\otimes\bQ_l.\] 
b) Under this decomposition of $H^i_\cont(X,\bQ_l(n))$, cup-product by
$e$ has matrix
\[\begin{pmatrix}0&Id\\0&0\end{pmatrix}.\] \end{prop}

\begin{pf} a) is clear from the exact sequence in proposition \ref{co7.2};
b) follows immediately from a).\end{pf}

\begin{rk}
We recover the conclusion of theorem \ref{t4.6} in a more concrete way.
\end{rk}

\begin{cor}\label{c7.2} If conjecture \ref{co7.1} holds, then the $\bQ_l$-adic
cycle map of \cite[(6.14)]{jannsen}
\[CH^n(X)\otimes\bQ_l\to H^{2n}_\cont(X,\bQ_l(n))\] is injective for any
smooth variety $X$ over $\bF_p$, and the composition
\[CH^n(X)\otimes\bQ_l\to H^{2n}_\cont(X,\bQ_l(n))@>\cdot e>>
H^{2n+1}_\cont(X,\bQ_l(n))\] is bijective.\qed \end{cor}

\begin{cor} \label{c7.6} Conjecture \ref{co7.1} implies the existence of
Beilinson's conjectural filtration in its weak form (\cf \cite[conj. 
2.1]{jannsen4}) on $CH^n(X)$ for any smooth, projective variety $X$ over a
field $F$ of characteristic $p$.  \end{cor}

\begin{pf} We may assume $F$ finitely generated. The proof is then a
variant of \cite[lemma 2.7]{jannsen4}, by applying corollary \ref{c7.2} to
a extension of $X$ to a smooth, projective scheme over a smooth model of
$F$ over $\bF_p$ (in other words, we use $\tilde
H^{2n}_\cont(X/\bF_p,\bQ_l(n))$ instead of
$H^{2n}_\cont(X,\bQ_l(n))$ as in {\it loc. cit.}).\end{pf}

\begin{cor}\label{c8.5} Assume conjecture \ref{co7.1}. Then, for any smooth variety
$X$ and any $i$, the composition
\[H^{i-n}(X,\sK_n^M)\otimes\bQ_l\to H^i_\cont(X,\bQ_l(n))\to H^i_\cont(\overline
X,\bQ_l(n))^G\]
is an isomorphism.
\end{cor}

\begin{pf} This follows from propositions \ref{p7.3} and \ref{p5.8}.\end{pf}

\begin{rk} Corollary \ref{c8.5} extends the Tate conjecture.
\end{rk}

\begin{cor}\label{c7.3} For any smooth variety $X$ over $\bF_p$, let
$d_i(n)=\dim_{\bQ_l} H^i_\cont(X,\bQ_l(n))$. Assume conjecture \ref{co7.1}.
Then, for all $i\le n$, $H^i(X,\sK_n^M)\otimes \bQ$ is a
finite-dimensional vector space of dimension $d_{i+n}(n)-d_{i+n-1}(n)+\dots$.
\end{cor}

\begin{pf} This is clear from corollary
\ref{c8.5} and proposition
\ref{l6.1}.\end{pf}

\begin{cor}\label{c7.3bis} Assume conjecture \ref{co7.1}.
Let $X/\bF_p$ be smooth of pure dimension $d$. Then, with notation as in
\eqref{eq20}, the order of $\det(1-tF|H^i_c(\overline X,\bQ_l))$ at $t=p^{n-d}$ is
$\dim_\Q H^{2d-n-i}(X,\sK_n^M)\otimes \bQ$.  \end{cor}

\begin{pf}  Since conjecture \ref{co7.1} implies condition $S^n(X)$, this order
equals
$\dim_{\bQ_l}H^i_c(\overline X,\bQ_l(d-n))^G=\dim_{\bQ_l}H^i_c(\overline
X,\bQ_l(d-n))_G$. Since $X$ is smooth, by (geometric) Poincar\'e duality, there is
an isomorphism
\[H^i_c(\overline X,\bQ_l(d-n))_G\simeq (H^{2d-i}(\overline X,\bQ_l(n))^G)^*.\]

The claim now follows from corollary
\ref{c8.5}.\end{pf}

\begin{cor}\label{c7.4} If conjecture \ref{co7.1} holds for all $n\ge 0$,
then 
\begin{thlist} 
\item For $X$ smooth projective over $\F_p$, $K_i(X)$
is torsion for $i>0$ (Parshin's conjecture).  
\item For $X$ smooth over $\F_p$, there are isomorphisms for $i,n\ge 0$
\[K_i(X)^{(n)}_\Q\simeq H^{n-i}(X,\sK_n)\otimes\bQ\simeq
H^{n-i}(X,\sK_n^M)\otimes\bQ.\] 
\item The Beilinson-Soul\'e conjecture holds rationally: for any smooth
$X$, $K_i(X)^{(n)}_\Q=0$ for $i\ge 2n$.  
\item The Bass conjecture holds rationally: for any smooth variety $X$
over $\bF_p$, the groups $K_i(X)\otimes \bQ$ are finite dimensional vector
spaces.  
\item For any field $F$ of characteristic $p$ and any $n\ge 0$, the map
$K_n^M(F)\otimes \bQ\to K_n(F)\otimes\bQ$ is an isomorphism.
\end{thlist}
\end{cor}

\begin{pf} (i) follows from corollary \ref{c7.1} and Geisser's theorem
\cite[th. 3.3 b)]{geisser}. (ii) also follows from results of Geisser
\cite[th. 3.4 (ii) and cor. 3.6]{geisser}. This implies (iii) trivially,
and (iv) follows from (ii) and corollary \ref{c7.3} Finally, (v) also
follows from \cite[th. 3.4 (ii)]{geisser}.\end{pf}

\begin{cor}\label{c7.10} Assume conjecture \ref{co7.1}. If $X/\bF_p$ is smooth
affine of dimension
$d$, then
$H^i(X,\sK_n^M)$ is torsion for $i>d-n$. In particular, $CH^n(X)$ is torsion if
$d<2n$.
\end{cor}

\begin{pf} This follows from
proposition \ref{p7.3} a), theorem \ref{t2.4} and
the cohomological dimension of affine schemes \cite[expos\'e XIX]{sga4}.\end{pf}

\begin{cor}\label{p7.4} {\rm (Bass-Tate conjecture, \cf \cite[question p.
390]{bt})} If conjecture
\ref{co7.1} holds in weight
$n$, then
$K_n^M(F)\otimes\bQ=0$ for any extension $F$ of
$\bF_p$, of transcendence degree $< n$.
\end{cor}

\begin{pf} We may assume $F$ finitely generated; then this follows from
corollary \ref{c7.10} (or from corollary \ref{c4.5}).\end{pf}

\begin{rks}
\begin{enumerate}
\item The bound is sharp by \cite[prop. 5.10]{bt}.
\item This also follows from corollary \ref{c7.1} and
Geisser's
theorem
\cite[th. 3.4]{geisser}, but the above reason seems more enlightening.
\end{enumerate}
\end{rks}

\begin{thm}\label{t7.2} Let $F$ be a field of characteristic $p$ and let  
$X_F$ be a smooth, projective variety of dimension $d$ over $F$.\\
a) Conjecture \ref{co7.1} in weight $n$ implies the Tate conjecture in
codimension $n$ and condition $S^n$ for $X_F$.\\
b) If Conjecture \ref{co7.1} holds for $n\le d$, it implies the algebraicity of the
K\"unneth components of the diagonal, Hard Lefschetz for cycles modulo numerical
equivalence and the strong version of Beilinson's conjectural filtration (\cf 
\cite[strong conj. 2.1]{jannsen4}) on all
Chow groups of $X_F$.
\end{thm}

\begin{pf} a) As in the proof of corollary \ref{c7.6}, we
assume $F$ finitely generated over $\bF_p$ and extend
$X_F$ to a smooth, projective morphism
\[f:X\to U\]
over a smooth model $U$ of $F$. We have the Leray spectral sequence for
$l$-adic cohomology ($l$-adic sheaves)
\begin{equation}\label{eq7.1}
E_2^{p,q}=H^p_\cont(U,R^qf_*^\ladic\bQ_l(n))\Rightarrow
H^{p+q}_\cont(X,\bQ_l(n)).\end{equation}

As seen in the proof of corollary \ref{t5.3}, this spectral sequence degenerates at
$E_2$
\cite{deligne2}. Using proposition \ref{p7.3}, we therefore get a commutative
diagram
\begin{equation}\label{eq8.1}\begin{CD}
H^{n-2}_\Zar(X,\sK_n^M)\otimes \bQ_l\oplus H^{n-1}_\Zar(X,\sK_n^M)\otimes
\bQ_l\iso H^{2n-1}_\cont(X,\bQ_l(n))\iso \bigoplus_{p+q=2n} E_2^{p-1,q}\\
@V{\left(\begin{smallmatrix}0&1\\0&0\end{smallmatrix}\right)}VV @V{\cdot
e}VV @V{\cdot e}VV\\
H^{n-1}_\Zar(X,\sK_n^M)\otimes \bQ_l\oplus CH^n(X)\otimes
\bQ_l\iso H^{2n}_\cont(X,\bQ_l(n))\iso \bigoplus_{p+q=2n} E_2^{p,q}.
\end{CD}\end{equation}
Here, cup-product by $e$ maps $E_2^{p-1,q}$ to $E_2^{p,q}$. We therefore
get an isomorphism
\[CH^n(X)\otimes\bQ_l\iso \bigoplus_{p+q=2n} E_2^{p,q}/e\cdot
E_2^{p-1,q}.\]

In particular, the right
hand side surjects onto 
\[E_2^{0,2n}=H^0_\cont(U,R^{2n}f_*^\ladic\bQ_l(n))=H^{2n}_\cont(\overline
X,\bQ_l(n))^{\pi_1(U)}.\]

Passing to the limit over $U$, we get the Tate conjecture and some information in
passing on the conjectural filtration.

As for $S^n(X_F)$, the degeneration of the spectral
sequence, lemma \ref{l4.1} and theorem \ref{t4.6} yield a
commutative diagram
\[\begin{CD}
H^{0}_\cont(\overline U,R^{2n} f_*^\ladic\bQ_l(n))^G@>>>
H^{0}_\cont(\overline U,R^{2n} f_*^\ladic\bQ_l(n))_G\\
@A{\wr}AA @VVV\\
H^0_\cont(U,R^{2n} f_*^\ladic\bQ_l(n))@>\cdot 
e>\sim>H^1_\cont(U,R^{2n} f_*^\ladic\bQ_l(n))
\end{CD}\]
in which the right vertical map is injective. Here $\overline U=U\otimes_{\bF_p}\overline \bF_p$. This yields
\begin{equation}\label{eq15} H^{0}_\cont(\overline
U,R^{2n}f_*^\ladic\bQ_l(n))^G\iso H^{0}_\cont(\overline
U,R^{2n}f_*^\ladic\bQ_l(n))_G.\end{equation}

On the other hand, by
\cite[cor. 3.4.13]{weilII}, the $\pi_1(\overline U)$-module
$R^{2n}f_*^\ladic\bQ_l(n)$ is semi-simple, hence
\begin{equation}\label{eq16}H^{0}_\cont(\overline
U,R^{2n}f_*^\ladic\bQ_l(n))=R^{2n} f_*^\ladic\bQ_l(n)^{\pi_1(\overline U)}\iso
R^{2n} f_*^\ladic\bQ_l(n)_{\pi_1(\overline U)}\end{equation} and then \eqref{eq15}
and \eqref{eq16} yield
\[R^{2n} f_*^\ladic\bQ_l(n)^{\pi_1(U)}\iso R^{2n} f_*^\ladic\bQ_l(n)_{\pi_1(U)}.\]

Passing back to the generic fibre, we get $S^n(X_F)$.

b) The algebraicity of the K\"unneth components follows from \cite[\S
3]{tate2}. We prove the assertion on Hard Lefschetz, since we didn't trace
a proof in the literature. We first note that a) and \cite[prop.
2.6]{tate2} imply that homological equivalence equals numerical
equivalence on $X_F$. Let $h$ be the cycle class of a hyperplane section of
$X_F$ (relative to some projective embedding), and let us denote by
$A^n(X_F)$ the group of cycles of codimension $n$ on $X_F$ modulo
numerical equivalence. Let $d=\dim X_F$. We have to show that, for all
$n<d/2$, the product by $h^{d-2n}$ \[A^n(X_F)\otimes\Q@>\cdot h^{d-2n}>>
A^{d-n}(X_F)\otimes\Q\] is bijective.  By the Tate conjecture and the
equality of homological and numerical equivalences, this translates as the
bijectivity of \[H^{2n}(\overline X_F,\Q_l(n))^{G_F}@>\cdot h^{d-2n}>>
H^{2d-2n}(\overline X_F,\Q_l(d-n))^{G_F}.\]

By the geometric Hard Lefschetz, this map is injective. On the other hand,
Poincar\'e duality gives an isomorphism
\[H^{2n}(\overline X_F,\Q_l(n))^{G_F}\iso Hom(H^{2d-2n}(\overline
X_F,\Q_l(d-n))_{G_F},\Q_l).\]

Hence the two $\Q_l$-vector spaces $H^{2n}(\overline X_F,\Q_l(n))^{G_F}$
and $H^{2d-2n}(\overline X_F,\Q_l(d-n))_{G_F}$ have the same dimension.
But $\dim H^{2d-2n}(\overline X_F,\Q_l(d-n))_{G_F}=\dim
H^{2d-2n}(\overline X_F,\Q_l(d-n))^{G_F}$ by $S^n(X_F)$ and cup-product by
$h^{d-2n}$ must be bijective. 

The strong form of Beilinson's conjecture follows from this, corollary
\ref{c7.6} and \cite[lemma 2.2]{jannsen}.  \end{pf}

\begin{cor} \label{c7.7} With notation as in theorem \ref{t7.2} and under
the assumptions of theorem \ref{t7.2} b), Murre's conjecture
\cite[1.4]{murre} holds for $X_F$.  
\end{cor}

\begin{pf} This follows from theorem \ref{t7.2} b) and \cite[th.
5.2]{jannsen4}.\end{pf}

\begin{rks}
\begin{enumerate}
\item Using $H^{2n+1}_\cont(X,\bQ_l(n))$ yields a different
presentation of the filtration on $CH^n(X_F)\otimes\bQ_l$.
\item One can mimic the proof of \cite[lemma 3.1]{jannsen4}, thereby obtaining
extra information on the direct limit of the spectral sequences \eqref{eq7.1} and a
direct proof of the strong form of Beilinson's filtration conjecture. More
precisely we note that, for
$(U,f,X)$ be as in the proof of theorem
\ref{t7.2}, when $U$ shrinks to smaller and smaller neighbourhoods of $\Spec F$,
the spectral sequences \eqref{eq7.1} have a direct limit
\begin{equation}\label{eq7.2}\tilde E_2^{p,q}=\tilde
H^p_\cont(F/\bF_p,H^q(\overline X_F,\bQ_l(n)))\Rightarrow \tilde
H^{p+q}_\cont(X_F/\bF_p,\bQ_l(n))
\end{equation}
which degenerates and only depends on $X_F$. 
\end{enumerate}
\end{rks}

\begin{cor}\label{c7.9} Let $F,X_F$ be as in theorem \ref{t7.2}, and let
$\delta=\trdeg(F/\bF_p)$. If conjecture
\ref{co7.1} holds for all $n\le d$, then, in the spectral sequence \eqref{eq7.2}, 
\begin{thlist}
\item $\tilde E_2^{p,q}=0$ for $p+q<n$, $p+d<n$, $\delta+q<n$, $\delta+d<n$ or
$p>n+1$;
\item For all $q$, $\tilde E_2^{n,q}@>\cdot e>> \tilde E_2^{n+1,q}$ is surjective.
\end{thlist}
Moreover, if $X$ is defined over a finite field, then $\tilde E_2^{p,q}=0$ for
$2p+q<2n$.
\end{cor}

\begin{pf} All inequalities of (i), except the last one, and the last statement
follow from corollary
\ref{t5.3}. For the last inequality of (i) and for (ii), let us write more precisely
$\tilde E_2^{p,q}(n)$ in order to keep track of the Tate twist $n$. By Hard
Lefschetz for
$l$-adic cohomology, we have an isomorphism
\[\tilde E_2^{p,q}(n)\iso \tilde E_2^{p,2d-q}(n+d-q).\]

The right hand side group is a direct summand of $\tilde
H^{p+2d-q}_\cont(X_F/\bF_p,\bQ_l(n+d-q))$; under conjecture \ref{co7.1} (for $d-n$),
this group is $0$ for $p>n+1$ by proposition \ref{p7.3} a). On the other hand, by
proposition \ref{p7.3} b), cup-product by $e$
\[\tilde H^{n+2d-q}_\cont(X_F/\bF_p,\bQ_l(n+d-q))@>\cdot e>> \tilde
H^{n+2d-q+1}_\cont(X_F/\bF_p,\bQ_l(n+d-q))\]
is surjective, which proves the second assertion of corollary \ref{c7.9}.
\end{pf}

We use this to get much more precise information on the conjectural Beilinson
filtration, not only on Chow groups but also on $\sK^M$-cohomology groups. Note how
the structure obtained is reminiscent of a pure Hodge structure.

\begin{thm}\label{t7.4} With notation and assumptions as in corollary \ref{c7.9},
there is for all $i$ a filtration on $H^i(X_F,\sK_n^M)\otimes\bQ$
\[F^p H^i(X_F,\sK_n^M)_\bQ,\quad p\ge 0\]
such that
\begin{thlist}
\item The action of correspondences on the associated graded $\gr^p
H^i(X_F,\sK_n^M)_\bQ$ factors through numerical equivalence.
\item $\gr^p
H^i(X_F,\sK_n^M)_\bQ=0$ for $p>n$, $p+d<n$, $p>\delta+i$ or $n>\delta+d$.
\item Usual naturality properties (contravariance, compatibility with products).
\item There are ``higher Abel-Jacobi maps"
\[\gr^p H^i(X_F,\sK_n^M)_\bQ\to \tilde H^p_\cont(F/\bF_p,H^{i+n-p}_\cont(\overline
X_F,\bQ_l(n)))/e\cdot H^{p-1}_\cont(F/\bF_p,H^{i+n-p}_\cont(\overline
X_F,\bQ_l(n)))\]
such that the induced maps on 
$\gr^p H^i(X_F,\sK_n^M)_\bQ\otimes_\bQ \bQ_l$ are isomorphisms.
\end{thlist}
In particular, there is a canonical isomorphism
\[H^i(X_F,\sK_n^M)\otimes\bQ_l\simeq \coprod_{p\le n}\gr^p
H^i(X_F,\sK_n^M)_\bQ\otimes_\bQ\bQ_l.\]
Moreover, if $X$ is defined over a finite field, then $\gr^p
H^i(X_F,\sK_n^M)_\bQ=0$ for $p<n-i$.
\end{thm}

\begin{pf} Using a diagram analogous to \eqref{eq8.1}, we get an isomorphism
\[H^i(X_F,\sK_n^M)\otimes\bQ_l\simeq \coprod_{p+q=i+n}\tilde E_2^{p,q}(n)/e\cdot
\tilde E_2^{p-1,q}(n).\]

This isomorphism depends on the Hard Lefschetz theorem; it can be made canonical
by \cite{deligne2}. By corollary
\ref{c7.9}, we have 
\[\tilde E_2^{p,q}(n)/e\cdot
\tilde E_2^{p-1,q}(n)=0\quad \text{for}\quad p\notin [n-i,n].\]

Define $F^p H^i(X_F,\sK_n^M)\otimes\bQ_l$ as the sub-vector space
$\coprod_{r\ge p}\tilde E_2^{r,i+n-r}(n)/e\cdot \tilde E_2^{r-1,i+n-r}(n)$.

\begin{lemma}\label{l8.4} This filtration is $\bQ$-rational; the corresponding
filtration $F^p H^i(X_F,\sK_n^M)_\bQ$ on $H^i(X_F,\sK_n^M)\otimes\bQ$ does not
depend on the choice of $l$.
\end{lemma}

\begin{pf} Choose $f:X\to U$ as in the proof of theorem \ref{t7.2}. Applying
$Rf_*$ to the isomorphism \eqref{eq33} restricted to $X$, we get an isomorphism
\[Rf_*\bQ_l(0)^c\oo^L \alpha^*\sK_n^M[-n]\iso
Rf_*\bQ_l(n)^c\]
and the rational filtration comes from the Leray spectral sequence
\[R^pf_*(\sH^q(\bQ_l(0)^c\oo^L \alpha^*\sK_n^M[-n]))\Rightarrow
R^{p+q}f_*(\bQ_l(0)^c\oo^L \alpha^*\sK_n^M[-n]).\]
\end{pf}

Lemma \ref{l8.4} shows that $F^p H^i(X_F,\sK_n^M)_\bQ\otimes_\bQ \bQ_l\iso F^p
H^i(X_F,\sK_n^M)\otimes\bQ_l$.
\end{pf}

\subsection{Homological transformation}\label{8.4}

Let $i:Z\inj X$ be a closed injection in $Sch/\bF_p$, with $X$ smooth of pure
dimension $d$. Applying
$Ri^!$ to \eqref{eq33} restricted to $X$, we get a morphism in $\sD^+(Ab(Z_\et))$:
\[Ri^!(\bQ_l(0)^c\oo^L\alpha^*\sK_{n,X}^M[-n])\to Ri^!\bQ_l(n)^c_X.\]

By proposition \ref{p2.2} c), the right hand side can be rewritten as
$L_l(n-d)_Z[-2d]\otimes \bQ_l$. On the other hand, since
$\bQ_l(0)^c=\bQ_l[0]\oplus\bQ_l[-1]$, the left hand side can be rewritten
\[\bQ_l(0)^c\oo^L Ri^!\alpha^*\sK_{n,X}^M[-n].\]

Using the natural transformation $\alpha^*Ri^!\to Ri^!\alpha^*$, we get a
composition
\[\bQ_l(0)^c\oo^L \alpha^*Ri^!\sK_{n,X}^M[-n]\to L_l(n-d)_Z[-2d]\otimes \bQ_l.\]

By proposition \ref{p8.3}, the left hand side is isomorphic to
$\bQ_l(0)^c\oo^L\alpha^*\sM_{n-d,Z}[-n-d]$. After shifting and changing $n$ into
$n-d$, we therefore get a natural transformation
\begin{equation}\label{eq8.3}
\bQ_l(0)^c\oo^L\alpha^*\sM_{n,Z}[-n]\to L_l(n)_Z\otimes\bQ_l.
\end{equation}

It can be shown by the method of \cite[proof of th. 4.1]{gillet} that for $Z$
quasi-projective \eqref{eq8.3} does not depend on the choice of the closed embedding
$i$ and has the usual naturality properties; we shall not need this here, so we skip
the proof.

\begin{lemma}\label{l8.5} The natural map \[\alpha^*Ri^!\sK_{n,X}^M\otimes\bQ\to
Ri^!\alpha^*\sK_{n,X}^M\otimes\bQ\]
is an isomorphism.
\end{lemma}

\begin{pf} This follows from corollary \ref{c5.2} and the classical isomorphism of
functors
\[Ri^!R\alpha_*\iso R\alpha_*Ri^!.\]
\end{pf}

\begin{thm} \label{t8.4} Conjecture \ref{co7.1} for all $n$ is equivalent to the
following:
\eqref{eq8.3} is an isomorphism for all $n$ and all $Z\in Sch/\bF_p$. In
particular, under conjecture \ref{co7.1}, one has isomorphisms
\[H_i^\cont(Z,\bQ_l(n))\simeq H_{i-n}(Z,\sK_n^M)\otimes\bQ_l\oplus
H_{i-n+1}(Z,\sK_n^M)\otimes\bQ_l.\]
\end{thm}

\begin{pf} This follows from the construction of \eqref{eq8.3} and lemma \ref{l8.5}.
\end{pf}

\begin{prop} \label{p8.4} Under conjecture \ref{co7.1}, there are isomorphisms for
all quasi-projective
$X\in Sch/\bF_p$:
\[H_i(X,\sK_n^M)\otimes \bQ\iso H_i(X,\sK_n)\otimes \bQ\simeq \gr_n
K'_{i-n}(X)\otimes\bQ\]
where the second group is defined analogously to the first
(with $K$-groups) and the third one is defined in \cite[th. 7 and 7.4]{soule}.
\end{prop}

\begin{pf} The first isomorphism follows from corollary \ref{c7.4} (v); the second
follows from the spectral sequence of \cite[th. 8 (iv)]{soule} and corollary
\ref{c7.4} (v) again, which implies that this spectral sequence degenerates at
$E^2$.
\end{pf}

With the help of theorem \ref{t8.4}, we can now extend corollary \ref{c7.3bis} to
singular varieties:

\begin{thm}\label{t8.5} Assume conjecture \ref{co7.1}.
Let $X$ be a quasi-projective variety over $\bF_p$. Then, with notation as in
\eqref{eq20}, the order of $\det(1-tF|H^i_c(\overline X,\bQ_l))$ at $t=p^{-n}$ is
$\dim_\bQ H_{i-n}(X,\sK_n^M)\otimes \bQ=\dim_\bQ \gr_n K'_{i-2n}(X)\otimes\bQ$.
\end{thm}

\begin{pf} This follows from theorem \ref{t4.3}, theorem \ref{t8.4} and
proposition \ref{p8.4}.
\end{pf}

Theorem \ref{t8.5} implies (and precises) Soul\'e's conjecture that
\[\ord_{s=n}\zeta(X,s)=\sum_i (-1)^{i+1}\dim_\bQ \gr_n K'_i(X)\otimes\bQ\]
\cite[conj. 2.2]{souleicm}.

\section{Integral refinement}\label{6}

In all this section, we assume resolution of singularities for varieties over
$\bF_p$. We stress this assumption occasionally. The reader who does not want to
make such an assumption is advised to skip this section.

\subsection{Review of motivic cohomology}

Recall
the motivic complexes
$\bZ(n)$ of Suslin and Voevodsky
\cite{suvo}. We shall denote by $H^*(X,\bZ(n))$ (\resp $H^*_\et(X,\bZ(n))$)
the groups respectively denoted by $H^*_B(X,\bZ(n))$ and $H^*_L(X,\bZ(n))$ in
\cite{voem}, and similarly when replacing $\bZ(n)$ by other $A(n):= A\oo^L\bZ(n)$,
where $A$ is an abelian group. Occasionally we shall write $H^*_M$ instead of
$H^*$.

\begin{prop}\label{p8.2} We have quasi-isomorphisms
\begin{equation}\label{eq22} \bQ(n)\iso R\alpha_*\alpha^*\bQ(n) 
\end{equation} and \begin{equation}\label{eq11} \alpha^*\bZ(n)\oo^L
\bZ/l^\nu\iso \mu_{l^\nu}^{\otimes n}[0]\end{equation} where $\alpha$ is
the projection of the big \'etale site of $\Spec\bF_p$ onto its big
Zariski site; in \eqref{eq22}, $\bQ(n)$ denotes $\bZ(n)\otimes\bQ$.
\end{prop}

\begin{pf} This follows from \cite[th. 2.5 and 2.6]{voem}.
\end{pf}

From this we get some information on the \'etale motivic cohomology groups
of a smooth scheme $X$: 

\begin{prop}\label{p5.3} Let $\bZ_{(l)}(n)$ denote $\bZ(n)\otimes \bZ_{(l)}$,
where $\bZ_{(l)}$ is the localisation of $\bZ$ at $l$. Then, for any smooth scheme
$X$ over $\bF_p$, $H^i_\et(X,\bZ_{(l)}(n))$ is
uniquely divisible-by-finite for $i<n$, torsion for
$i=2n+1$, torsion of cofinite type for $i=2n+2$ and finite for
$i>2n+2$. If moreover $X$ is projective, then $H^i_\et(X,\bZ_{(l)}(n))$ is
\begin{thlist}
\item uniquely divisible-by-finite for $i<2n$;
\item with finite torsion for $i=2n$.
\end{thlist}
\end{prop}

\begin{pf} By \cite[cor. 2.3]{voem}, we have
$H^i_\Zar(X,\bZ(n))\allowbreak =0$ for $i>2n$. It follows from this and
\eqref{eq22} that $H^i_\et(X,\bZ(n))$ is torsion
for $i>2n$. The claims now follow from proposition
\ref{p4.1}, corollaries \ref{c4.2} and \ref{c4.3} and the exact sequence
\begin{equation}\label{eq23}
\dots \to H^{i-1}_\et(X,\bQ_l/\bZ_l(n))@>\beta>>
H^i_\et(X,\bZ_{(l)}(n))
\to
H^i_\et(X,\bQ(n))\to
H^i_\et(X,\bQ_l/\bZ_l(n))\to\dots\end{equation}
\end{pf}

\begin{prop}\label{p8.1} Let $X$ be a smooth variety over $\bF_p$. Then,
\begin{thlist}
\item $\sH^i(\bZ(n))=0$ for $i>n$.
\item $\sH^n(\bQ(n))\simeq \sK_n^M\otimes\bQ$
\item $H^p(X,\bZ(n))=0$ for $p>2n$
\item $H^{2n}(X,\bZ(n))$ is canonically isomorphic to $CH^n(X)$.
\end{thlist}
\end{prop}

\begin{pf} See \cite[cor. 2.3 and cor. 2.4]{voem}.\end{pf}

For any $\bF_p$-scheme $X$, let $F_X$ denote the geometric Frobenius of $X$
\cite[expos\'e XV]{sga5}. We note that $F_X$ defines a finite correspondence from
$X$ to itself in the sense of \cite{vtr}. It is easy to
check that $X\mapsto F_X$ ``extends" to an endofunctor $F$ of the category of
effective geometrical motives
$DM_{gm}^{eff}(\bF_p)$. It is clear that, for $M,M'\in DM_{gm}^{eff}(\bF_p)$,
\[F_M\otimes F_{M'}\mapsto F_{M\otimes M'}\]
under the natural map
\[End_{DM}(M)\otimes_\bZ End_{DM}(M')\to End_{DM}(M\otimes M').\]

\begin{prop}\label{p5.4} For any smooth scheme $X$ over $\bF_p$, the
geometric
Frobenius of $X$ acts on $H^*(X,\bZ(n))$ by multiplication by $p^n$.
\end{prop}

\begin{pf} We interpret motivic cohomology groups as Homs in
$DM_{gm}^{eff}(\bF_p)$. So we have
\[H^i(X,\bZ(n))=Hom_{DM}(M(X),\bZ(n)[i])\]

By naturality of the geometric Frobenius, its action on the latter groups
is the same as the action induced by the geometric Frobenius of $\bZ(n)$. 
Therefore it suffices to show that its class in $End_{DM}(\bZ(n))$ is
$p^n$ times the identity. Since $\bZ(n)=\bZ(1)^{\otimes n}$, we are
reduced to the case $n=1$ by the remarks just before the statement of
proposition \ref{p5.4}. Then $\bZ(1)$ is a shift of the Lefschetz motive
$L$. By \cite[prop.  2.1.4]{vtr}, we can compute $F_{L}$ in the category
of effective Chow motives;  then the result is classical \cite{soule2}.
\end{pf}

\subsection{The conjecture}

From \eqref{eq11} we get compatible maps ($\nu\ge 0$)
\[\alpha^*\bZ(n)\to \mu_{l^\nu}^{\otimes n}[0].\]

Tensoring it with $\bZ_l(0)^c$, we get compatible maps
\[
\bZ_l(0)^c\oo^L\alpha^*\bZ(n)\to\mu_{l^\nu}^{\otimes n}[0].\]
which yield a morphism
\begin{equation}\label{eq3}
\bZ_l(0)^c\oo^L\alpha^*\bZ(n)\to\bZ_l(n)^c.\end{equation}

The following lemma is obvious from \eqref{eq11} and lemma \ref{l2.1}.

\begin{lemma}\label{l5.1} The morphism \eqref{eq3}$\oo^L\bZ/l$ is an
isomorphism.\qed
\end{lemma}

\begin{conj}\label{co1} The morphism \eqref{eq3} is an isomorphism.
\end{conj}

In view of lemma \ref{l5.1}, conjecture \ref{co1} is
equivalent to

\begin{conj}\label{co2} The morphism \eqref{eq3}$\oo\bQ$ is an isomorphism.
\end{conj}

Assume resolution of singularities. By proposition \ref{p8.1} and corollary
\ref{c4.7}, the morphism
\eqref{eq3}$\otimes
\bQ$ induces a homomorphism of
\'etale sheaves
\begin{equation}\label{eq2.3}\alpha^*\sK_n^M\otimes \bQ_l@>u_n>>
\sH^n(\bQ_l(n)^c).\end{equation}

It is easy to check that this morphism coincides with that of proposition
\ref{p5.10} (tensored with $\Q$).

\begin{prop} \label{p5.-1} Under resolution of singularities, conjectures
\ref{co1} and \ref{co2} are equivalent to the following:\\
\begin{thlist}
\item $\sH^i(\bZ(n))\otimes \bQ=0$ for $i<n$.
\item The morphism $u_n$ from \eqref{eq2.3} is an isomorphism.
\item The morphism of \'etale sheaves
\[\sH^n(\bQ_l(n)^c)@>e>> \sH^{n+1}(\bQ_l(n)^c)\]
given by cup-product by $e$ is an isomorphism.
\item $\sH^i(\bQ_l(n)^c)=0$ for $i>n+1$.
\end{thlist}
\end{prop}

\begin{pf} Assume conjecture \ref{co2}. Then (i) follows from theorem
\ref{t4.2} and \cite[th. 2.5]{voem}, (ii) is clear, (iii) follows from the
fact that it is true for $n=0$ and (iv) follows from the fact that
$\bZ(n)$ is acyclic in degrees $>n$ \cite{suvo}. Conversely, (i)--(iv)
imply that \eqref{eq3}$\otimes \bQ$ induces an isomorphism on cohomology
sheaves, hence is an isomorphism.\end{pf}

\begin{cor}\label{c8.3} Under conjecture \ref{co2}, $H^i(X,\bQ(n))\simeq
H^{i-n}(X,\sK_n^M)\otimes \bQ$.\qed
\end{cor}

\begin{thm}\label{t8.1} Under resolution of singularities, conjectures
\ref{co7.1} and \ref{co1} are equivalent.
\end{thm}

\begin{pf} From proposition \ref{p5.-1} it easily follows that conjecture
\ref{co1} implies conjecture \ref{co7.1}. Conversely, conjecture
\ref{co7.1} implies (ii)--(iv) in proposition \ref{p5.-1} and we are left
to see that it also implies (i). By theorem \ref{t3.2} and the purity
theorem for motivic cohomology, we have to show that
\[H^i(X,\Q(n))=0\quad\text{for}\quad i<n\] for any smooth, projective $X$.
The proof is similar to that of theorem \ref{t7.1}, using corollary
\ref{c7.1} and proposition \ref{p5.4}.\end{pf}

The following proposition is trivial from theorem \ref{t2.2} b).

\begin{prop}\label{p5.9} If conjecture \ref{co1} holds, then
$\bZ_l(n)^c\simeq \bZ^c\oo^L\alpha^*\bZ(n)\oo^L\bZ_l$.\qed \end{prop}

In particular,
we recover proposition \ref{p7.2} in a more explicit way.

\begin{prop}\label{p5.0a} Suppose conjecture \ref{co1} holds. Let $X$ be a
smooth variety. Then there is a long exact sequence 
\begin{multline*}
\dots \to H^{i-2}_\et(X,\bZ(n))\otimes \bQ_l@>\partial>>
H^i_\et(X,\bZ(n))\otimes \bZ_l\\ \to H^i_\cont(X,\bZ_l(n))\to
H^{i-1}_\et(X,\bZ(n))\otimes \bQ_l\to\dots
\end{multline*} 
Moreover, the connecting homomorphism $\partial$ is given by tensoring by
$\bZ_l$ the composition 
\[\begin{CD} H^{i-2}_\et(X,\bQ(n))@>\partial>>
H^i_\et(X,\bZ_{(l)}(n))\\ @VVV @A{\beta}AA\\
H^{i-2}_\et(X,\bQ_l/\bZ_l(n))@>\cdot e>> H^{i-1}_\et(X,\bQ_l/\bZ_l(n))
\end{CD}\] 
where $\beta$ is the Bockstein map of \eqref{eq23}. In
particular, the image of $\partial$ is torsion.  
\end{prop}

\begin{pf} The exact sequence is clear; the computation of $\partial$
follows from proposition \ref{p2.3}.\end{pf}

\begin{prop}\label{c5.4} Keep the assumptions and notation as in proposition
\ref{p5.0a}. Then\\
a) For any $X$,
\begin{thlist}
\item The map $H^i_\et(X,\bZ(n))\otimes \bZ_l\to H^i_\cont(X,\bZ_l(n))$
is bijective for $i\le n$ and injective for $i=n+1$.
\item As a module over $\bZ_{(l)}$, $H^i_\et(X,\bZ_{(l)}(n))$ is finite for $i<n$,
finitely generated for
$i=n$ and an extension of a finitely generated module by a torsion
divisible module of finite corank for $i>n$. We have
\[\corank_{\bQ_l/\bZ_l}H^i_\et(X,\bZ(n))=\rank_\bZ H^{i-2}_\et(X,\bZ(n)).\]
\end{thlist}
b) If
$X$ is moreover projective, then
\begin{thlist}
\item The map $H^i_\et(X,\bZ(n))\otimes \bZ_l\to
H^i_\cont(X,\bZ_l(n))$ is bijective for $i\neq 2n+1, 2n+2$ and injective for
$i=2n+1$.
\item As a module over $\bZ_{(l)}$, $H^i_\et(X,\bZ_{(l)}(n))$ is finite for $i\neq
2n, 2n+2$, finitely generated for $i=2n$ and cofinitely generated for $i=2n+2$.
\item Suppose conjecture \ref{co1} holds for all $n\le d$. If $d=\dim X$,
the group $H^{2d+2}_\et(X,\bZ_{(l)}(d))$ is canonically isomorphic to
$\bQ_l/\bZ_l$ and the pairings
\[H^i_\et(X,\bZ_{(l)}(n))\times H^{2d+2-i}_\et(X,\bZ_{(l)}(d-n))\to
H^{2d+2}_\et(X,\bZ_{(l)}(d))\]
are perfect.
\item We have isomorphisms
\begin{align*}
H^{2n}_\et(X,\bZ_{(l)}(n))_\tors&\iso H^{2n}_\cont(X,\bZ_l(n))_\tors\\
H^{2n+1}_\et(X,\bZ_{(l)}(n))&\iso H^{2n+1}_\cont(X,\bZ_l(n))_\tors\\
H^{2n+2}_\et(X,\bZ_{(l)}(n))_\cotors&\iso H^{2n+2}_\cont(X,\bZ_l(n)).
\end{align*}
\end{thlist}
\end{prop}

\begin{pf} Everything follows easily from propositions \ref{p5.3} and
\ref{p5.0a} plus the fact that the ring extension $\bZ_l/\bZ_{(l)}$ is
faithfully flat, except for the corank computation. We have exact
sequences \[H^{i-1}_\cont(X,\bZ_l(n))\to H^{i-2}_\et(X,\bZ(n))\otimes
\bQ\to Div(H^i_\et(X,\bZ(n))\otimes\bZ_l)\to 0.\] Since the last group is
torsion and the first one is finitely generated, the image of the left map
must be a lattice in $H^{i-2}_\et(X,\bZ(n))\otimes \bQ$. This gives the
claim.\end{pf}

\begin{rk} In particular, propsition \ref{c5.4} b) shows that conjecture
\ref{co1} implies S. Lichtenbaum's conjectures 1) to 7) of \cite[\S
7]{licht} on the values of \'etale motivic cohomology, after localisation
at $l$. By \cite{milne2}, this implies conjecture 8) of {\it loc. cit.} on
the values of the zeta function of $X$ at nonnegative integers (note that
there is an obvious misprint in the formula of {\it loc. cit.}). We now
extend this to get a formula for the principal part of $\zeta(X,s)$ at
$s=n$ in terms of motivic cohomology for {\it any} smooth variety,
thereby completing the result of corollary \ref{c7.3bis}.  \end{rk}

\begin{lemma}\label{l8.1} The homomorphism
$\partial$ of proposition
\ref{p5.0a} induces a homomorphism 
\[\bar\partial^{i-2}:H^{i-2}_\et(X,\bZ(n))\otimes \bQ_l/\bZ_l@>>>
H^i_\et(X,\bZ_{(l)}(n))_{\tors}\]
with finite kernel and cokernel.
\end{lemma}

\begin{pf} The existence of $\bar\partial^{i-2}$ is clear from the
factorisation of $\partial$ given in proposition \ref{p5.0a}. By the exact
sequence in this proposition, $\Ker\partial$ is a lattice in
$H^{i-2}_\et(X,\bQ(n))$ and $\Coker\bar\partial^{i-2}\iso
H^i_\cont(X,\bZ_l(n))_\tors$, which is finite.  \end{pf}

\begin{thm}\label{t4.4} Let $X$ be a smooth variety of pure dimension
$d$ over $\bF_p$. If conjecture \ref{co1} holds, then:
\begin{thlist}
\item the order of $\det(1-tF|H^i_c(\overline X,\bQ_l))$ at $t=p^{n-d}$ is
$\rank H^{2d-i}(X,\bZ(n))$.
\item $\zeta(X,s)=(1-p^{n-d-s})^{a_{d-n}}\varphi(s)$, with $a_{d-n}=\sum
(-1)^{i+1}\rank
H^{i}(X,\bZ(n))$  and
\item $|\varphi(d-n)|_l=
\left|\prod_{i=0}^{2d}\ind(\bar\partial^{i})^{(-1)^{i+1}}\right|_l$.
\end{thlist}
\end{thm}

Note that we don't get a separate formula for the principal part of
$\det(1-tF|H^i_c(\overline X,\bQ_l))$ at $t=p^{-n}$ in terms of motivic
cohomology. See remark \ref{r9.1} for more on this point.

\begin{pf} (i) and (ii) follow from corollaries \ref{c7.3bis} and \ref{c8.3}.
(iii) With the
help of lemma \ref{l8.1}, we extend the commutative
square of proposition \ref{p5.0a} into a bigger commutative diagram:
\[\begin{CD}
H^{i-2}_\et(X,\bZ(n))\otimes \bQ_l/\bZ_l@>\bar\partial^{i-2}>>
H^i_\et(X,\bZ_{(l)}(n))_{\tors}\\
@VVV @A{\beta}AA\\
H^{i-2}_\et(X,\bQ_l/\bZ_l(n))@>\cdot e>> H^{i-1}_\et(X,\bQ_l/\bZ_l(n))\\
@VVV @AAA\\
H^{i-2}_\et(\overline X,\bQ_l/\bZ_l(n))^G @>\rho^{i-2}>> H^{i-2}_\et(\overline
X,\bQ_l/\bZ_l(n))_G
\end{CD}\]
which follows from lemma \ref{l4.1} as in the proof of proposition
\ref{p5.8}. By corollary \ref{c6.2}, the bottom horizontal map has finite
kernel and cokernel, and so has the top one by lemma \ref{l8.1}. By the
cross of exact sequences
\[\begin{CD}
&&0\\
&&@VVV\\
&&H^{i-1}(X,\bZ(n))\otimes \Q_l/\Z_l\\
&&\Big\downarrow &\searrow\scriptstyle \alpha^{i-1}\\
0\to H^{i-2}(\overline X,\Q_l/\Z_l(n))_G@>>> H^{i-1}(X,\Q_l/\Z_l(n))@>>>
H^{i-1}(\overline X,\Q_l/\Z_l(n))^{G}\to 0\\
&\searrow\scriptstyle\beta^{i-2} & @VVV\\
&&H^{i}_\et(X,\Z_{(l)}(n))_\tors\\
&&@VVV\\
&& 0
\end{CD}\]
the compositions
\[\begin{CD}
H^{i-1}(X,\bZ(n))\otimes \Q_l/\Z_l@>\alpha^{i-1}>>H^{i-1}(\overline
X,\Q_l/\Z_l(n))^{G}\\
H^{i-2}(\overline X,\Q_l/\Z_l(n))_G@>\beta^{i-2}>> H^{i}_\et(X,\Z_{(l)}(n))_\tors
\end{CD}\]
have same kernel and cokernel. Since, by the above, $\Ker \alpha^{i-2}$ and
$\Coker \beta^{i-2}$ are finite for any $i$, this proves

\begin{lemma} \label{l8.2} For all $i$, $\alpha^i$ and $\beta^i$ have finite
kernel and cokernel; moreover $\ind(\alpha^{i-1})=\ind(\beta^{i-2})$.\qed
\end{lemma}

By lemmas \ref{l8.2}, \ref{l4.3} b) and the above commutative
diagram, we have
\[\ind(\bar\partial^{i-2})=\ind(\alpha^{i-2})\ind(\rho^{i-2})\ind(\beta^{i-2})=
\ind(\alpha^{i-2})\ind(\rho^{i-2})\ind(\alpha^{i-1}).\]

By corollary \ref{c6.2}, $|f_i(p^{n-d})|_l=|\ind(\rho^{2d-i})|_l$. It follows
that
\[|f(p^{n-d})|_l=\left|\prod_{i=0}^{2d}f_i(p^{n-d})^{(-1)^{i+1}}\right|_l
=\left|\prod_{i=0}^{2d}\ind(\bar\partial^{i})^{(-1)^{i+1}}\right|_l\]
as desired.
\end{pf}

\begin{rk}\label{r9.1} Theorem \ref{t4.4} (i) becomes much more suggestive if one
compares the group $H^{2d-i}(X,\bZ(n))$ to
\[H^i_c(X,\bZ(d-n)),\] 
motivic cohomology with compact supports \cite{vtr}. It is likely that these two
groups are dual to each other after tensoring by $\bQ$. More precisely, by
\cite[\S 9]{frivoe} we have
\[H^i_c(X,\bZ(d-n))\simeq H_{2d-i}(X,\bZ(n))\]
since $X$ is smooth, where the right hand group is motivic homology, and there is a
pairing
\begin{multline*}H_{2d-i}(X,\bZ(n))\times
H^{2d-i}(X,\bZ(n))\\
=Hom_{DM}(\bZ(n)[2d-i],M(X))\times Hom_{DM}(M(X),\bZ(n)[2d-i])\\
\to End_{DM}(\bZ(n)[2d-i])=\bZ.\end{multline*}

Here the motivic groups are Zariski, but they do or should coincide with the
\'etale groups after tensoring by $\bQ$. This strongly hints that one should
formulate a version of conjecture
\ref{co1} for continuous \'etale homology analogous to theorem \ref{t8.4}, involving
this time Borel-Moore motivic homology. This would presumably allow one to remove
the smoothness assumption in theorem
\ref{t4.4}. In this generality, there is a pairing between \'etale motivic
cohomology with compact supports and
\'etale motivic Borel-Moore homology, analogous to the above, that one can
conjecture to be perfect after tensoring by $\bQ$. What would then be the
relationship between the principal part of $\det(1-tF|H^i_c(\overline X,\bQ_l))$ at
$t=p^{n-d}$ and the determinant of this pairing (together with the torsion of the
two motivic groups)?

We even suspect that this conjectural picture formally follows from conjecture
\ref{co1}.
\end{rk}

Specialising to the case where $X$ is projective and taking account of proposition
\ref{c5.4} and the functional equation of $\zeta(X,s)$, we get:

\begin{cor}\label{c8.2} With the assumptions and notation of theorem \ref{t4.4},
if $X$ is moreover projective, then
$\zeta(X,s)=(1-p^{n-s})^{a_n}\varphi(s)$, with
$a_n=\rank  H^{2n}(X,\bZ(n))$  and
\[|\varphi(n)|_l=
\left|\prod_{i\neq 2n}
|H^i_\et(X,\bZ(n))|^{(-1)^i}\ind\left(H^{2n}(X,\bZ_{(l)}(n))\otimes
\bQ_l/\bZ_l@>\bar\partial^{2n}>> H^{2n+2}(X,\bZ(n))\right)^{-1}\right|_l.\]
\end{cor}

Alternatively, using corollary \ref{c4.1} and proposition \ref{c5.4} b) (iv), we
get the following formula, which is the $l$-primary part of \cite[\S 7, conj.
8]{licht} (after correcting the weight $n+2$ to $n$ in {\it loc. cit.}):

\begin{thm}\label{t8.3} With the assumptions and notation of corollary \ref{c8.2},
we have
\[|\varphi(n)|_l=\left|\prod_{i\neq 2n,2n+2} 
|H^i_\et(X,\bZ_{(l)}(n))|^{(-1)^{i}}\cdot
{|H^{2n}_\et(X,\bZ_{(l)}(n))_\tors||H^{2n+2}_\et(X,\bZ_{(l)}(n))_\cotors|\over
|R_n(X)|_l}\right|_l\]
where $H^{2n+2}_\et(X,\bZ_{(l)}(n))_\cotors$ is the quotient of
$H^{2n+2}_\et(X,\bZ_{(l)}(n))$ by its maximal divisible subgroup (a finite group
by corollary \ref{c4.2} and proposition \ref{c5.4} b) (iv)) and $R_n(X)$ is the
discriminant of the pairing
\[H^{2n}_\et(X,\bZ_{(l)}(n))/\tors\times
H^{2d-2n}_\et(X,\bZ_{(l)}(d-n))/\tors\to
H^{2d}_\et(X,\bZ_{(l)}(d))/\tors\simeq\bZ_{(l)}\] with respect to any bases of
$H^{2n}_\et(X,\bZ_{(l)}(n))/\tors$ and
$H^{2d-2n}_\et(X,\bZ_{(l)}(d-n))/\tors$ (this is well-defined up to a unit
of $\bZ_{(l)}$).\qed
\end{thm}

Proposition \ref{c5.4} shows that, under resolution of singularities and
conjecture \ref{co1}, the \'etale motivic cohomology groups
$H^i_\et(X,\bZ_{(l)}(n))$ are finitely generated over $\bZ_{(l)}$ for $i\le 2n+1$
and any smooth, projective variety $X/\bF_p$. What about the converse? The answer
is surprisingly negative and is given by theorem \ref{p5.5} and corollary
\ref{c8.4} below.

\begin{prop}\label{p5.1} Conjecture \ref{co1} holds for $n\le m$ when
restricted to $\sS_d$ if and only if, for any smooth, projective $X$ of
dimension $\le d$, the homomorphisms
\[H^*_\et(X,\bZ_l(0)^c\oo^L\alpha^*\bZ(n))\to H^*_\cont(X,\bZ_l(n))\]
are isomorphisms for all $n\le m$.
\end{prop}

\begin{pf} Let $K(n)$ be the cone of \eqref{eq3}: conjecture \ref{co1}
holds on $\sS_d$ if and only if, for any $X$ smooth of dimension
$\le d$, the groups $H^*(X,K(n))$ are $0$. The proof is similar to
that of theorem \ref{t7.1}. Note that $X\mapsto
H^*(X,K(n))$ defines a pure graded cohomology theory with values in
$\bQ_l$-vector spaces. The first claim follows from purity for continuous
\'etale cohomology (theorem \ref{t2.1}) and purity for \'etale motivic
cohomology \cite{suvo}, while the second one follows from lemma
\ref{l5.1}. Proposition \ref{p5.1} now follows once again from theorem
\ref{t3.2} b), noting that $K(n)$ has transfers.\end{pf}

\begin{thm}\label{p5.5} Let $X$ be a smooth, projective variety
over $\bF_p$. Assume that the groups $H^i_\et(X,\bZ_{(l)}(n))$ are
finitely generated
$\bZ_{(l)}$-modules for all $i\le 2n+1$. Then:\\
a) The Tate conjecture holds for $X$ in codimension $n$;\\
b) Let $K(n)$ be the cone of
\eqref{eq3}. Then $\bH^i_\et(X,K(n))=0$ for $i\neq 2n, 2n+1$ and there is
an exact sequence
\[0\to \bH^{2n}_\et(X,K(n))\to H^{2n}_\cont(X,\bQ_l(n))@>\cdot e>>
H^{2n+1}_\cont(X,\bQ_l(n))\to \bH^{2n+1}(X,K(n))\to 0.\]
\end{thm}

\prf By the assumptions and proposition \ref{p5.3}, $H^i_\et(X,\bZ_{(l)}(n))$ is
finite for
$i\neq 2n,2n+2$ and torsion for $i= 2n+2$. It follows that
\[H^i(X,\bZ(n))\otimes \bZ_l\iso H^i(X,\bZ(n)\otimes \bZ_l(0)^c)\]
for $i\neq 2n+1,2n+2$, hence that 
$\bH^i_\et(X,\bZ(n)\oo^L\bZ_l(0)^c)$ is a finitely generated
$\bZ_l$-module for $i\neq 2n,2n+1$.
So are the $H^i_\cont(X,\bZ_l(n))$ for all $i$. Since the
$\bH^i_\et(X,K(n))$ are uniquely divisible, it follows that they are $0$,
except perhaps for $i=2n,2n+1$. We have a commutative diagram with exact rows and
columns
\[\begin{CD}
&&H^{2n}_\cont(X,\bZ_l(n))\\
&&@V{A}VV\\
&&H^{2n}(X,K(n))\\
&&@V{B}VV\\
0\to H^{2n+1}(X,\bZ(n))\otimes\bZ_l @>>>H^{2n+1}(X,\bZ(n)\otimes \bZ_l(0)^c)@>>>
H^{2n}(X,\bQ(n))\otimes \bZ_l\\
&&@VVV\\
&&H^{2n+1}_\cont(X,\bZ_l(n))\\
&&@VVV\\
&&H^{2n+1}(X,K(n))\\
&&@VVV\\
&&H^{2n+2}(X,\bZ(n)\otimes \bZ_l(0)^c).
\end{CD}\]

Let $M=\Coker A$ and $N$ its torsion subgroup. Then $B(N)$ maps to $0$ in $
H^{2n}(X,\bQ(n))\otimes \bZ_l$, hence is contained in
$H^{2n+1}(X,\bZ(n))\otimes\bZ_l$. But this group is finite and $B(N)$ is divisible,
hence $B(N)=0$. This shows that $M$ is torsion-free and, since
$H^{2n}_\cont(X,\bZ_l(n))$ is finitely generated, this shows that $A=0$. We
first deduce that
\[H^{2n}(X,\bZ(n))\otimes \bQ_l\iso H^{2n}(X,\bZ(n)\otimes \bQ_l(0)^c)\iso  
H^{2n}_\cont(X,\bQ_l(n))\]
which gives a). Using this and tensoring the above vertical exact sequence by $\bQ$,
we get an exact sequence

\[0\to \bH^{2n}_\et(X,K(n))\to H^{2n}_\cont(X,\bQ_l(n))@>C>>
H^{2n+1}_\cont(X,\bQ_l(n))\to \bH^{2n+1}(X,K(n))\to 0\]
(note that $H^{2n+2}(X,\bZ(n)\otimes \bZ_l(0)^c)$ is torsion). Theorem \ref{p5.5}
results from this and the following lemma, which is an easy consequence of
corollary \ref{c2.2}:

\begin{lemma}\label{l8.3} The map $C$ is cup-product by $e$.\qed
\end{lemma}

\begin{cor}\label{c8.4}
Conjecture \ref{co1} holds for weights $\le d$ and smooth varieties
of dimension $\le d$ if and only if the condition of proposition \ref{p5.5} and
condition $S^n$ are verified for projective ones for all $n\le d$.\qed
\end{cor}

It is difficult in general to get more precise information on the cokernel
of the
$l$-adic cycle map
\[CH^n(X)\otimes\bZ_l\to H^{2n}_\cont(X,\bZ_l(n)).\]

We shall do so when the ground field is finite. (Compare \cite[remark
5.6 a)]{milne2}.)

\begin{prop}\label{p5.7} If conjecture \ref{co1} holds,
then, for any smooth, projective variety $X$ over $\bF_p$, there is an
exact sequence
\begin{multline*}
0\to H^{2n-1}_M(X,\bQ_l/\bZ_l(n))\to
CH^n(X)\otimes\bZ_l\oplus
H^{2n}_\cont(X,\bZ_l(n))_\tors\\
\to H^{2n}_\cont(X,\bZ_l(n))
\to
CH^n(X)\otimes\bQ_l/\bZ_l\to H^{2n}_\et(X,\bQ_l/\bZ_l(n)).\end{multline*}
 \end{prop}

\begin{pf} This follows readily from the commutative diagram of exact sequences
\[\begin{CD}
0&&0\\
@VVV @VVV\\
H^{2n-1}_M(X,\bQ_l/\bZ_l(n))@>>> H^{2n-1}_\et(X,\bQ_l/\bZ_l(n))\\
@VVV @VVV\\
H^{2n}(X,\bZ(n))\otimes\bZ_l@>>> H^{2n}_\cont(X,\bZ_l(n))\\
@VVV @VVV\\
H^{2n}(X,\bZ(n))\otimes\bQ_l\iso H^{2n}_\cont(X,\bQ_l(n))\\
@VVV @VVV\\
H^{2n}_M(X,\bQ_l/\bZ_l(n))@>>> H^{2n}_\et(X,\bQ_l/\bZ_l(n))\\
@VVV\\
0
\end{CD}\]

Here the top $0$s and the middle isomorphism come from proposition
\ref{c5.4} b), and the bottom $0$ comes from proposition \ref{p8.1} (iii).\end{pf}

\section{The one-dimensional case}\label{one}

In this section, we prove:

\begin{thm}\label{t9.1} Conjecture \ref{co7.1} is true when restricted to smooth
curves over
$\bF_p$.
\end{thm}

\begin{pf} We reduce as usual to the case of a smooth, projective curve $X/\bF_p$,
that we may suppose connected. For
$n\ge 2$, $\bQ_l(n)^c_{|X}=0$ by corollary \ref{c4.5}. On the other hand,
$(\sK_n^M)_X\otimes\bQ=0$: to see this, we reduce by theorem \ref{t5.1} to showing
that $K_n^M(\bF_p(X))\otimes \bQ=0$. This follows from the case $n=2$, which is a
consequence of Harder's theorem. Theorem \ref{t9.1} is proven for $n\ge 2$.

Assume now $n=1$. Let $k$ be the field of constants of $X$. We have $\sK_1^M=\bG_m$
and:
\begin{align*}
H^0_\Zar(X,\bG_m)\otimes\bQ&=k^*\otimes\bQ=0;\\
H^1_\Zar(X,\bG_m)\otimes\bQ&=\Pic(X)\otimes\bQ@>\deg>\sim> \bQ;\\
H^i_\Zar(X,\bG_m)\otimes\bQ&=0\quad\text{for}\quad i>1.
\end{align*}

On the other hand, we can compute $H^i_\cont(X,\bQ_l(1))$ with the
help of lemma \ref{l2.2} a) and the Kummer exact sequence. By proposition
\ref{p4.1}, it suffices to consider $i=2,3$. For $i=2$, we have short exact
sequences
\[\begin{CD}
0@>>> \Pic(X)/l^\nu@>>> H^2_\et(X,\mu_{l^\nu})@>>> {}_{l^\nu} Br(X)@>>>
0\\
0@>>> Br(X)/l^\nu @>>>  H^3_\et(X,\mu_{l^\nu})@>>> {}_{l^\nu}H^3_\et(X,\bG_m)@>>> 0
\end{CD}\]
and $Br(X)=0$, $H^3_\et(X,\bG_m)\iso \bQ/\bZ$ by class field theory.
This gives:
\begin{align*}
H^2_\cont(X,\bQ_l(1))&=\Pic(X)^{\hat{}}_l\otimes\bQ@>\deg>\sim>\bQ_l;\\
H^3_\cont(X,\bQ_l(1))&=T_l(\bQ/\bZ)\otimes\bQ\simeq\bQ_l.
\end{align*}

One checks easily from these computations that \eqref{eq33} induces isomorphisms
of cohomology groups.\end{pf}

\section{A conjecture in the Zariski topology}\label{zariski}

In this section as in section \ref{6}, we assume resolution of singularities.

Assume conjecture \ref{co1} holds. We can rewrite it as an exact triangle
\[\alpha^*\bZ(n)\otimes \bZ_l\to \bZ_l(n)^c\to\alpha^*\bZ(n)\otimes
\bQ_l[-1] \to \alpha^*\bZ(n)\otimes
\bZ_l[1].\]

Taking propositions \ref{p8.1} (i) and \ref{p5.-1} (i), (ii) into account, this can
be rewritten as another exact triangle
\[\alpha^*\bZ(n)\otimes \bZ_l\to \bZ_l(n)^c\to
\alpha^*\sK_n^M\otimes \bQ_l[-n-1]\to \alpha^*\bZ(n)\otimes
\bZ_l[1]\]
or an isomorphism
\[\alpha^*\bZ(n)\otimes
\bZ_l\iso\,\text{fibre}(\bZ_l(n)^c\to\alpha^*\sK_n^M\otimes \bQ_l[-n-1]).\]

Still by proposition \ref{p5.-1}, the right hand side can equally be written
$\tau_{\le n}\bZ_l(n)^c$. Applying $R\alpha_*$, we therefore get as a consequence
of conjecture
\ref{co1} a conjectural isomorphism
\[R\alpha_*\alpha^*\bZ(n)\otimes\bZ_l\iso R\alpha_*(\tau_{\le
n}\bZ_l(n)^c).\]

Recall now the Kato conjecture: for any field $F$ of characteristic $\neq l$, the
Galois symbol
\[K_n^M(F)/l\to H^n(F,\mu_l^{\otimes n})\]
is bijective, and the
Beilinson-Lichtenbaum conjecture \cite{suvo},
\cite{voem}
\[\bZ(n)\iso \tau_{\le n+1}R\alpha_*\alpha^*\bZ(n).\]

By the main result of \cite{suvo}, the Kato conjecture is equivalent to the
Beilinson-Lichtenbaum conjecture under resolution of singularities. Moreover, by the
main result of \cite{voem}, the Kato conjecture holds for $l=2$ (without assuming
resolution of singularities).

Assuming the Kato conjecture and conjecture \ref{co1}, we therefore get a new
conjectural isomorphism:
\[\bZ(n)\otimes\bZ_l\iso\tau_{\le n+1}R\alpha_*(\tau_{\le
n}\bZ_l(n)^c).\]

This is equivalent to an isomorphism $\bZ(n)\otimes\bZ_l\iso\tau_{\le 
n}R\alpha_*\bZ_l(n)^c$ plus an injection of Zariski sheaves
\[R^{n+1}\alpha_*\bZ_l(n)^c\Inj \sK_n^M\otimes\bQ_l\]
the latter being the same as asking for $R^{n+1}\alpha_*\bZ_l(n)^c$ to be
torsion-free. We can now state our conjecture in the Zariski topology:

\begin{conj}\label{co3}
a) The natural morphism in the derived category of abelian sheaves over
the big Zariski site of $\Spec\bF_p$
\[\bZ(n)\otimes\bZ_l\to \tau_{\le n}R\alpha_*\bZ_l(n)^c\]
is an isomorphism.\\
b) The Zariski sheaf $R^{n+1}\alpha_*\bZ_l(n)^c$ is torsion-free.
\end{conj}

\begin{rk} b) is closely related to `higher Hilbert 90'.
\end{rk}

We have proven part a) of the following (under resolution of singularities):

\begin{prop} \label{p5.6}Under resolution of singularities,\\
a) Conjecture
\ref{co1} and the Kato conjecture
imply conjecture \ref{co3}.\\
b) Conjecture \ref{co3} is equivalent to the conjunction of the Kato conjecture
and either
\begin{thlist}
\item finite generation of the $\bZ_{(l)}$-module $H^i(X,\bZ_{(l)}(n))$ for
all $i\le n+1$ and all smooth $X$
\item same as (ii), but only for smooth, projective $X$.
\end{thlist}
\end{prop}

\begin{pf} Assume conjecture \ref{co3}. Tensoring the isomorphism by $\bZ/l$, we
get the Beilin\-son-Lichten\-baum conjecture, hence the Kato conjecture. Moreover,
applying lemma
\ref{l2.2} b), we get the finite generation statement. Without
assumptions, the cone of
$R\alpha_*\alpha^*\bZ(n)\otimes\bZ_l\to R\alpha_*\bZ_l(n)^c$
has uniquely divisible cohomology sheaves. Therefore, the same holds in
degree $\le n$
for the cone of the truncations
$\tau_{\le n+1}R\alpha_*\alpha^*\bZ(n)\otimes\bZ_l\to 
\tau_{\le n+1}R\alpha_*\bZ_l(n)^c$. Assuming the Kato conjecture, hence the
Beilinson-Lichtenbaum conjecture, this means that the cohomology sheaves of
\[K(n)=\text{cone}(\bZ(n)\otimes \bZ_l\to \tau_{\le
n+1}R\alpha_*\bZ_l(n)^c)\]
are uniquely
divisible in degrees $\le n$ and that $R^{n+1}\alpha_*\bZ_l(n)^c$ is
torsion-free. Let
now $X$ be a smooth, projective variety over $\bF_p$.
Using lemma \ref{l2.2} b), the finite generation assumption implies that
$\bH^i(X,K(n))=0$ for $i\le n$. Applying theorem \ref{t3.2} b), we get the
same result for all smooth $X$. It follows that $\sH^i(K(n))=0$ for $i\le
n$, which is equivalent to the statement of conjecture \ref{co3}.\end{pf}

\begin{cor}\label{c8.1} Under resolution of singularities and the Kato conjecture,
finite generation of the groups $H^i(X,\bZ_{(l)}(n))$ for  $n>0$, $i\le 0$ and
smooth, projective varieties $X$ implies their vanishing for any smooth $X$
(motivic form of the Beilinson-Soul\'e vanishing conjecture).\qed
\end{cor} 

By proposition \ref{p5.6} a), conjecture \ref{co3} implies that the groups
$H^i(X,\bZ_{(l)}(n))$ are finitely generated over $\bZ_{(l)}$ for $i\le n+1$ and any
smooth variety $X$ over
$\bF_p$. What about $i>n+1$? Conjecture \ref{co1} implies that these groups have
finite rank, but not that they are finitely generated. In particular, it does not
imply the finite generation of the Chow groups of $X$. We shall now formulate a last
conjecture and show that, together with conjecture \ref{co1}, it implies this finite
generation.

\begin{conj}\label{co8.1} For any smooth scheme $X$ over a reasonable base
($\bF_p$, $\bZ[1/l]$, a separably closed field are reasonable), the groups 
\[H^i_\Zar(X,\sH^j_\et(\bZ/l^\nu(n)))\]
are finite for all $i,j,\nu,n$.
\end{conj}

\begin{thm} \label{t8.2} {\rm (Under resolution of singularities)} Assume that
the Kato conjecture holds, conjecture \ref{co1} holds and conjecture
\ref{co8.1} holds for $X$ projective and $\nu=1$,
$n=0$. Then the groups $H^i(X,\bZ_{(l)}(n))$ are finitely generated over $\bZ_{(l)}$
for all $i$ and any smooth $X/\bF_p$.
\end{thm}

\begin{pf} Once again, by theorem \ref{t3.2} it is sufficient to deal with the projective
case. Assume first that  $\Gamma(X,\bG_m)$ contains a primitive $l$-th root
of unity. Then, by the Beilinson-Lichtenbaum conjecture, we have for any $n$
an exact triangle on the small Zariski site of $X$:
\[\bZ/l(n-1)_M\to \bZ/l(n)_M\to \sH^n(\bZ/l)\to \bZ/l(n-1)_M[1]\]
where the index $M$ indicates motivic complexes.
From the assumption, we get inductively that $H^i_M(X,\bZ/l(n))$ is finite for all
$i,n$. If $\Gamma(X,\bG_m)$ does not contain a primitive $l$-th root
of unity, we reduce to this case by a classical transfer argument. By induction on
$\nu$ we then get the same result for $H^i_M(X,\bZ/l^\nu(n))$ and all $\nu\ge 1$.

The Soul\'e-Geisser motivic argument implies that $H^i_M(X,\bZ_{(l)}(n))$
has finite exponent for $i<2n$. The injection
\[H^i_M(X,\bZ_{(l)}(n))/l^\nu\Inj H^i_M(X,\bZ/l^\nu(n))\]
for $\nu$ sufficiently large then shows that this group is finite.

For $i=2n$, we use proposition \ref{p5.7}: it suffices to show that
$H_M^{2n-1}(X,\bQ_l/\bZ_l(n))$ is finite. We do this by showing that it has finite
exponent, just as above.
\end{pf}

\section{The case of rings of integers}\label{8}

In this section, let $S=\Spec \bZ[1/l]$ and $\sS$ be the category of regular
$S$-schemes. We don't have even a conjectural description of $\bZ_l(n)^c$ over
$\sS_\et$. We shall indicate the little knowledge we have, mostly for the
restriction of $\bZ_l(n)^c$ to the {\it small} \'etale site of $S$.

We begin with an elementary result which partially reduces the problem to
understanding what happens at the generic point of $S$. Unfortunately, we have to
state it in terms of \'etale homology.

\begin{prop}\label{p12.2} Let $\sX$ be a scheme of finite type over $S$,
$X$ its generic fibre and, for $p\neq l$, $\sX_p$ its closed fibre at $p$. Then, for
all
$n\in
\bZ$, there is a long exact sequence
\[\dots\to H_i^\cont(\sX/S,\bZ_l(n))\to \tilde
H_i^\cont(X/S,\bZ_l(n))\to\coprod_{p\neq l}H_{i-1}^\cont(\sX_p/\bF_p,\bZ_l(n))\to
H_{i-1}^\cont(\sX/S,\bZ_l(n))\to\dots\]
\end{prop}

\begin{pf} This follows from proposition
\ref{p2.2} a).\end{pf}

\begin{cor}\label{c12.2} Let $\sX$ be a scheme of finite type over $S$ and
$X$ its generic fibre. Then for all $n<0$ and
$i\in \bZ$, the natural map
\[H_i^\cont(\sX/S,\bQ_l(n))\to \tilde H_i^\cont(X/S,\bQ_l(n))\]
is an isomorphism.
\end{cor}

\begin{pf} This follows from proposition
\ref{p12.2} and corollary \ref{c4.6}.\end{pf}

Recall, for each prime $p$, the canonical generator of $H^1_\cont(\bF_p,\bZ_l)$
that we now denote by $e_p$. Given its importance in the theory in characteristic
$p$, the first question that comes to mind is: do the $e_p$ lift compatibly to
characteristic $0$? The answer is, predictably, yes but only up to a factor.

More precisely, the cyclotomic character gives a continuous homomorphism
\[e:\pi_1(S)\to \bZ_l^*\]
which is surjective by Gauss' theorem on the irreducibility of cyclotomic
polynomials. This defines a class in $H^1_\cont(S,\bZ_l^*)$. Consider the
continuous homomorphism
\begin{align*}
\ell:\bZ_l^*&\to \bZ_l\\
u&\mapsto 
\begin{cases}
{1\over l^2}\log(u^l)&\text{if $l>2$}\\
{1\over 8}\log(u^2)&\text{if $l=2$.}
\end{cases}
\end{align*}

This has kernel the roots of unity of $\bZ_l^*$ and coincides with $1/l$ of the
usual logarithm when restricted to $1+l\bZ_l$ ($l>2$) and with $1/4$ of the
usual logarithm when restricted to $1+4\bZ_2$
($l=2$). In all cases, it is surjective. Via
$\ell$, we associate to $e$ a class, still denoted by $e$:
\[e\in H^1_\cont(S,\bZ_l).\]

\begin{lemma}\label{l12.2} For all $p\neq l$, the image of $e$ in
$H^1_\cont(\bF_p,\bZ_l)$ is $\ell(p)e_p$.
\end{lemma}

\begin{pf} The absolute Frobenius automorphism at $p$ raises
roots of unity to the $p$-th power.
\end{pf}

\begin{defn}\label{d12.1} Let $\bR_0$ be ``the" field of real algebraic numbers
($\bR_0=\bR\cap \overline \bQ$ for a suitable embedding of $\overline \bQ$ into
$\bC$). For $n\in \bZ/2$, we denote by $\bZ(n)$ the $G_{\bR_0}$-module with
support $\bZ$ on which the action of $G_{\bR_0}$ is given by $\varepsilon^n$,
where $\varepsilon$ is its nontrivial character of order $2$. We set
\[A(n)=Ind_{G_{\bR_0}}^{G_\bQ}\bZ(n).\]
\end{defn}

\begin{thm}\label{t12.1} Let $j:\Spec \bQ\inj \Spec \bZ[1/2]$ be the inclusion
of the generic point. For
$n\ge 2$,\\ a) There is a non-canonical isomorphism
\[\bQ_l(n)^c_{|S}\simeq j_*A(n+1)\otimes \bQ_l[-1].\]
b) The natural morphism
\[Rj_*\alpha^*\bZ(n)_{|\bQ}\otimes \bQ_l\to \bQ_l(n)^c_{|S}\]
is an isomorphism.
\end{thm}

\begin{pf} a) Let $R$ be a ring of $S$-integers in a number field $F$ (we assume
$1/l\in R$). It is clear that $H^0_\cont(R,\bQ_l(n))=0$ and
$H^i_\cont(R,\bQ_l(n))=0$ for
$i\ge 3$ as well for reasons of cohomological dimension. By Soul\'e's main theorem
\cite{soule3}, we have
\begin{itemize}
\item $H^2(R,\bQ_l(n))=0$;
\item the Chern character induces an isomorphism 
\[K_{2n-1}(R)\otimes\bQ_l\iso H^1_\cont(R,\bQ_l(n)).\]
\end{itemize}

On the other hand, the localization exact sequence in algebraic $K$-theory shows
that
\[K_{2n-1}(R)\otimes\bQ\iso K_{2n-1}(F)\otimes\bQ.\]

Therefore, the Chern character induces an isomorphism
\begin{equation}\label{eq12.1}
j_*\alpha^*\sK_{2n-1}\otimes\bQ_l[-1]\iso \bQ_l(n)^c_{|S}.
\end{equation}

On the other hand, the Borel regulator can be interpreted as a
$Gal(\bC/\bR)$-equivariant homomorphism
\[K_{2n-1}(\bC)@>\rho_n>> \bR(n+1)\]
where $\bR(n):=\bZ(n)\otimes\bR$. By \cite{borel}, the composition
\[K_{2n-1}(F)@>>>\coprod_{v|\infty}
K_{2n-1}(F_v)@>(\rho_{2n-1})>>\coprod_{v|\infty}
\bR(n+1)^{G_v}\]
where $G_v=Gal(\bC/\bF_v)$ indices an isomorphism
\[K_{2n-1}(F)\otimes\bR\iso\coprod_{v|\infty}
\bR(n+1)^{G_v}\]
hence, in the limit, an isomorphism of $G_\bQ$-modules
\[K_{2n-1}(\overline \bQ)\otimes\bR\iso A(n+1)\otimes\bR.\]

By the theory of characters, this isomorphism implies an abstract isomorphism
\[K_{2n-1}(\overline \bQ)\otimes\bQ\simeq A(n+1)\otimes\bQ\]
hence an isomorphism
\[K_{2n-1}(\overline \bQ)\otimes\bQ_l\simeq A(n+1)\otimes\bQ_l\]
and a).

b) This follows from \eqref{eq12.1} and \cite[th. 4.1 and 4.2]{QL2}.
\end{pf}

\begin{rk}
Behind the noncanonical isomorphism of theorem \ref{t12.1} a) are of course
the determinants of the Borel regulators. This, and the prominent r\^ole of
archimedean places, shows how Arakelov geometry lurks in the background.
\end{rk}

Theorem \ref{t12.1} shows that the situation is markedly
different in characteristic $0$ and in characteristic $p$:

\begin{cor}\label{c12.1} For $n\ge 2$, cup-product by $e$
\[\sH^{i}(\bQ_l(n)^c_{|S})@>\cdot e>>\sH^{i+1}(\bQ_l(n)^c_{|S})\]
is $0$ for all $i$.\qed
\end{cor}

For $n<0$, contrary to
characteristic $p$, it is well-known that $\bQ_l(n)^c_{|S}\neq 0$ \cite{schneider}.
It is widely conjectured ({\it loc. cit.}):

\begin{conj}\label{co12.1} Theorem \ref{t12.1} a) extends to $n<0$.
\end{conj}

To describe $\bZ_l(1)^c_{|S}$, it is convenient to introduce more definitions. Let
$D$ be the discrete topological $G_\bQ$-module defined by
\[D=\colim_{[F:\bQ]<+\infty} Div(O_F)\]
and define the \'etale {\em sheaf of codivisors} $\Codiv$ as the
kernel of the surjective map
\[j_*D\to \Div.\]

The degree map on $D$ induces a map
\[\deg:\Codiv\to j_*\bZ.\]

We also define $\tilde A(0)$ as the kernel of the augmentation $A(0)\to\bZ$.
Finally, let $Br$ be the {\it Zariski} sheaf associated to the presheaf $U\mapsto
Br(U)$ on the small \'etale site of $S$.

\begin{prop}\label{p12.1} There are exact sequences
\begin{gather*}0\to j_*\tilde A(0)\otimes\bQ\to\bG_m\otimes\bQ\to
\Codiv\to 0\\
0\to Br\to \Codiv\otimes\bQ/\bZ@>\deg>> j_*\bQ/\bZ\to 0
\end{gather*}
except for the contribution of real places to the Brauer group.
\end{prop}

\begin{pf} This is a reformulation of well-known results. In the first sequence,
the left hand sheaf corresponds to $O_F^*\otimes\bQ$ for $F$ running through
number fields; its description follows from Dirichlet's units theorem (hence is
once again transcendental). The right hand side can be obtained similarly, using
the $S$-units theorem. The exact sequence concerning the Brauer group is a
reformulation of the Albert-Brauer-Hasse-Noether theorem.
\end{pf}

\begin{thm}\label{t12.2} a) $\sH^i(\bQ_l(1)^c_{|S})=0$ for $i\neq 1,2$.\\
b) The natural map
\[\bG_m[-1]\otimes\bZ_l@>>> \sH^1(\bZ_l(1)^c_{|S})\]
given by Kummer theory is an isomorphism.\\
c) There is an exact sequence
\[0\to \sH^2(\bQ_l(1)^c_{|S})\to \Codiv\otimes\bQ_l@>\deg>> j_*\bQ_l\to 0.\]
\end{thm}

\begin{pf} a) follows as in the proof of theorem \ref{t12.1} a). To see b),
consider $R\subset F$ as in the said proof. From the Kummer exact sequence, we get
short exact sequences
\[0\to R^*/(R^*)^{l^\nu}\to H^1(R,\mu_{l^\nu})\to {}_{l^\nu}\Pic(R)\to 0.\] 

Since all groups are finite, there are no $\lim^1$ and we get a resulting short
exact sequence
\[0\to \lim R^*/(R^*)^{l^\nu}\to H^1(R,\bZ_l(1))\to T_l(\Pic(R))\to 0.\] 

But $R^*$ is finitely generated and $\Pic(R)$ is finite (Dirichlet's theorem);
hence the left hand side is $R^*\otimes\bZ_l$ and the right hand side is $0$.

Finally, to see c), we use the similar exact sequences one degree higher:
\[0\to \Pic(R)/l^\nu\to H^2(R,\mu_{l^\nu})\to {}_{l^\nu}Br(R)\to 0\]
and proposition \ref{p12.1}.
\end{pf}

A nontrivial question is whether cup-product by $e$
\[\sH^1(\bQ_l(1)^c_{|S})@>\cdot e>>\sH^2(\bQ_l(1)^c_{|S})\]
is surjective. This seems related to the Leopoldt conjecture. Out of laziness we
don't examine this in more detail and just give a conjectural value for
$\bQ_l(0)^c_{|S}$, which is easily seen to be equivalent to the Leopoldt
conjecture.

\begin{conj} \label{co12.2} There is an exact triangle
\[j_*\bQ_l[-1]@>\cdot e>>\bQ_l(0)^c_{|S}@>>> j_*A(1)\otimes\bQ_l[-1]@>>>
j_*\bQ_l[0].\] In particular, $\sH^i(\bQ_l(0)^c_{|S})=0$ for $i\neq 1$.
\end{conj}

Further insight on the conjectural behaviour of $\bQ_l(n)^c$ in characteristic $0$
can be found in \cite{jannsen5} and \cite[\S 13]{jannsen2}. These two references,
together with theorems \ref{t12.1} and \ref{t12.2}, at least suggest:

\begin{conj} \label{co12.3} For $n\ge 0$, the map
\[\alpha^*\bZ(n)\otimes\bZ_l\to \tau_{\le n}j^*\bZ_l(n)_S\]
restricted to smooth varieties over $\bQ$ is an isomorphism.
\end{conj}

The case $n=0$ is trivial and the case $n=1$ is easily checked as in the proof of
theorem \ref{t12.2} b) (use the finite generation of $\Gamma(X,\sO_X^*)$ and of
$Pic(X)$ for $X$ regular of finite type over $\Spec \bZ$, which follow from
Dirichlet, Mordell-Weil and N\'eron-Severi). On the other hand, I have no idea of
the nature of $\tau_{> n}j^*\bZ_l(n)_S$.

\appendix
\section{The Bass conjecture implies the rational Beilinson-Soul\'e conjecture}

Let $n,i\ge 0$ be two integers; for any regular scheme
$X$, denote by $K_n(X)^{(i)}$ the $i$-th eigenspace of the Adams operations acting
on $K_n(X)$ \cite{soule}. For any prime number $l$ and any abelian group $A$,
denote by $A_{(l)}$ the localisation $A\otimes \bZ_{(l)}$ of $A$ at $l$. In this
appendix, we prove:

\begin{thm}\label{tA.1} Assume that $K_n(X)^{(i)}$ is finitely generated when
$X$ is of finite type over $\bZ$.\\
a) If $1\le i\le n/2$, then for any $d\ge 0$ there exists an effectively computable
integer $M(n,d)$ such that $M(n,d)K_n(X)^{(i)}_{(2)}=0$ for any regular scheme of
dimension $\le d$.\\
b) One can choose $M(n,d)$ such that, if $n=2i-1>1$, for any regular scheme
$X$ of dimension $\le d$, in the commutative diagram
\[\begin{CD}
K_{2i-1}(X)^{(i)}_{(2)}@>>> K_{2i-1}(\kappa(X))^{(i)}_{(2)}\\
@AAA @AAA\\
K_{2i-1}(X_0)^{(i)}_{(2)}@>>> K_{2i-1}(\kappa(X_0))^{(i)}_{(2)}
\end{CD}\]
all maps have kernel and cokernel killed by $M(n,d)$. Here $X_0$ is the scheme of
constants of $X$, \ie the normalisation of $\Spec\bZ$ in $X$, and $\kappa(X)$
denotes the total ring of functions of $X$.
\end{thm}

The proof goes along the same lines as the one of corollary \ref{c8.1}, but is
more involved. The game is to dodge resolution of singularities. The method goes
back to \cite{knote}.

The proof will be divided into a series of lemmas.

\begin{lemma}\label{lA.1} For all $n\ge 0$ there is an
effectively computable constant $M(n)$ such that, if $F$ is a field of
characteristic
$0$, 
\[M(n)K_n(F)^{(i)}\]
is uniquely $2$-divisible for $i\le n/2$ and, for $n=2i-1>1$, the natural map
\[K_n(F_0)^{(i)}\2\to K_n(F)^{(i)}\2\]
multiplied by $M(n)$ is injective with uniquely divisible cokernel.
\end{lemma}

\begin{pf} We have the Bloch-Lichtenbaum-Friedlander-Suslin-Voevodsky spectral
sequence
\[E_2^{p,q}=H^{p-q}(F,\bZ(-q))\Rightarrow K_{-p-q}(F)\quad (p,q\le 0)\]
(\cite{bl}, \cite{frivoe}, \cite{vtr}, see the discussion in \cite[2.1]{QL2}). By
Soul\'e \cite{soule4}, for all $k>0$ the Adams operation $\psi^k$ acts on this
spectral sequence, by multiplication by $k^{-q}$ on $E_2^{p,q}$ and converging to
the standard action on $K_{-p-q}(F)$. It follows that, for all $r\ge 2$, the
differentials $d_r$ are killed by 
\[w_{r-1}:=gcd_{k\ge 2} k^N(k^{r-1}-1)\]
where $N$ is large with respect to $r$ (\cf \cite[p. 498]{soule}). It follows
that, for any $p,q\le 0$, the group
\[w_1w_2\dots w_{-p}E_2^{p,q}\]
consists of universal cycles, while universal boundaries of degree $(p,q)$ are
killed by
\[w_1w_2\dots w_{-q}.\]

To summarise, there is a chain of inclusions
\[w_1\dots w_{-p}w_1\dots w_{-q}E_2^{p,q}\subseteq w_1\dots
w_{-q}E_\infty^{p,q}\subseteq w_1\dots w_{-q}E_2^{p,q}.\]

On the other hand, let $(F^iK_n(F))_{i\ge 0}$ be the filtration on $K_n(F)$
induced by the spectral sequence. The Adams operation $\psi^k$ respects this
filtration and acts on $\gr^iK_n(F)=E_\infty^{i-n,-i}$ by $k^i$. Let $\varphi$ be
the composition
\[K_n(F)^{(i)}\Inj K_n(F)\Surj K_n(F)/F^{i+1}K_n(F).\]

Then $w_1\dots w_{n-i}\Ker\varphi=0$ and $w_1\dots w_{i-1}\IM\varphi\subseteq
\gr^iK_n(F)$. Conversely, by the argument of \cite[p. 499]{soule}, this subgroup
has index $w_1\dots w_{n-i}$. In other words, there is a chain of inclusions
\[w_1\dots w_{n-i}w_1\dots w_{i-1}K_n(F)^{(i)}\subseteq w_1\dots
w_{n-i}E_\infty^{i-n,-i}\subseteq w_1\dots w_{i-1}\overline{K_n(F)}^{(i)}\]
where $\overline{K_n(F)}^{(i)}$ is a quotient of $K_n(F)^{(i)}$ by a
subgroup of exponent $w_1\dots w_{n-i}$.

Let $M(n)=(w_1\dots w_n)^4$. We see from the above that there exists for all $i$ a
chain of homomorphisms
\[H^{2i-n}(F,\bZ(n))@>>> ? @<<< K_n(F)^{(i)}\]
with kernel and cokernel killed by $M(n)$. 

By \cite[th. 3.1 a)]{QL2}, $H^{2i-n}(F,\bZ(n))$ is uniquely $2$-divisible for
$2i-n\le 0$, and by {\it loc. cit.}, th. 7.1, $H^1(F_0,\bZ\2(n))\to
H^1(F,\bZ\2(n))$ is injective with uniquely divisible cokernel (this relies on
\cite{suvo} and
\cite{voem}). The result follows.
\end{pf}

\begin{lemma}\label{lA.2} Lemma \ref{lA.1} also holds for $F$ of characteristic
$p>0$.
\end{lemma}

\begin{pf} We distinguish two cases:

1) $p=2$. By results of Geisser and Levine \cite{gele}, $\Coker(K_n^M(F)\to
K_n(F))$ is uniquely $2$-divisible. The result then follows from \cite[p. 498, cor.
1]{soule} (and its proof).

2) $p>2$. By the long homotopy exact sequence, it is enough to show that, for all
$i\le n/2$ and $\nu\ge 1$, $M(n)K_n(F,\bZ/2^\nu)^{(i)}$ is $0$ if $i< n/2$ and
that $M(n)K_n(F_0,\bZ/2^\nu)^{(i)}\to M(n)K_n(F,\bZ/2^\nu)^{(i)}$ is bijective
(\resp injective) for $n=2i-1$ (\resp $2i-2$). We may assume
$F$ perfect (transfer argument). Let $E$ be the field of fractions of the Witt
vectors on $F$. We have split short exact sequences
\[0\to K_n(F,\bZ/2^\nu)\to K_n(E,\bZ/2^\nu)\to K_{n-1}(F,\bZ/2^\nu)\to 0\]
which are respected by the Adams operations, and similarly for $F_0$. The result
follows from this, up to increasing $M(n)$ a bit.\end{pf}

\begin{lemma}\label{lA.3} For any $n,d$, there is a constant $M(n,d)$ such
that,  for any regular scheme $X$ of dimension $d$,
\[M(n,d)K_n(X)^{(i)}\]
is uniquely $2$-divisible for $i\le n/2$ and, for $n=2i-1>1$, in the commutative
diagram
\[\begin{CD}
K_{2i-1}(X)^{(i)}_{(2)}@>>> K_{2i-1}(\kappa(X))^{(i)}_{(2)}\\
@AAA @AAA\\
K_{2i-1}(X_0)^{(i)}_{(2)}@>>> K_{2i-1}(\kappa(X_0))^{(i)}_{(2)}
\end{CD}\]
with natural maps multiplied by $M(n,d)$, all maps are injective with uniquely
divisible cokernel.
\end{lemma}

\begin{pf} To do this, we reduce to the field case by using the Quillen spectral
sequence
\cite{quillen} (for $X$ and $X_0$)
\[E_1^{p,q}=\coprod_{x\in X^{(p)}}K_{-p-q}(\kappa(x))\Rightarrow K_{-p-q}(X),\]
Soul\'e's theorem that the Adams operations act on this spectral sequence \cite[p.
521, th. 4 (i), (ii)]{soule} and an argument similar to that in step 1 (\cf also
\cite{knote}).
\end{pf}

\begin{lemma}\label{lA.4} Theorem \ref{tA.1} holds for $X$ of finite type over
$\bZ$.
\end{lemma}

\begin{pf}This follows immediately from lemma \ref{lA.3}.
\end{pf}

\begin{lemma}\label{lA.5} Theorem \ref{tA.1} holds for $X=\Spec F$, $F$ a
field.
\end{lemma}

\begin{pf}We may assume $F$ finitely generated. Write $F=\kappa(X)$ for $X$ regular
of finite type over $\bZ$. The result follows from lemma \ref{lA.4}, a direct limit
argument and  lemma \ref{lA.2}.
\end{pf}

\noindent{\it Proof of theorem \ref{tA.1}}. Follows from step 5 and a new
application of the Quillen spectral sequence.\qed

\begin{cor}\label{cA.1} Under the assumption of theorem \ref{tA.1}, for any
smooth scheme $X$ over a field $F$ of characteristic $0$,\\
a) $H^{2i-n}(X,\bZ\2(i))=0$ for $i>0$ and $2i-n\le 0$.\\
b) In the diagram
\[\begin{CD}
H^1(X,\bZ\2(i))@>>> H^1(F(X),\bZ\2(i))\\
&&@AAA\\
&&H^1(F_0,\bZ\2(i))
\end{CD}\]
the two maps are isomorphisms.The same holds for $F$ of characteristic $>0$, under
resolution of singularities.
\end{cor}

\begin{pf} The case $X=\Spec F$ follows from theorem \ref{tA.1} by applying the
argument in the proof of lemma \ref{tA.1} backwards. The general case follows from
the coniveau spectral sequence plus purity for motivic cohomology.
\end{pf}

By the same method as in the proof of theorem \ref{tA.1}, one can further show that
the Bass conjecture implies the Parshin conjecture for smooth, projective varieties
$X$ over $\bF_p$ ($K_i(X)\otimes \bQ=0$ for $i>0$), hence recover the consequences
Geisser deduces from this conjecture in \cite[\S\S 3.3 and 3.4]{geisser} (see
corollaries \ref{c7.4} and \ref{p7.4} here).

\end{document}